\documentclass[a4paper]{amsart}

\usepackage{amsmath,amssymb,amsthm}
\usepackage{graphicx}
\usepackage[english]{babel}
\usepackage[utf8]{inputenc}
\usepackage[T1]{fontenc}
\usepackage{mathtools}
\usepackage[margin=1cm]{caption}
\usepackage{bbm}
\usepackage[margin=3cm]{geometry}
\usepackage{xcolor}
\usepackage{enumitem}
\usepackage{hyperref}
\allowdisplaybreaks

\newcommand{\N}{\mathbb{N}}
\newcommand{\Z}{\mathbb{Z}}
\newcommand{\Q}{\mathbb{Q}}
\newcommand{\R}{\mathbb{R}}

\newcommand{\E}{\mathbb{E}}

\newcommand{\e}{{\mathrm{e}}}

\newcommand{\PP}{\mathbb{P}}

\newcommand{\F}{\mathcal{F}}
\newcommand{\I}{\mathbbm{1}}

\newcommand{\g}{\mathbf{g}}
\newcommand{\A}{\mathbf{a}}
\newcommand{\Ev}{\mathcal{E}}
\newcommand{\eps}{\varepsilon}

\newcommand{\ahom}{\bar{\a}}

\newcommand{\Sav}{S} 
\newcommand{\Inc}{Z} 
\newcommand{\Bla}{L} 
\newcommand{\Ero}{\Ev^{\mathrm{rgh}}} 
\newcommand{\Ebd}{\Ev^{\mathrm{bd}}} 
\newcommand{\Ecp}{\Ev^{\mathrm{cpl}}} 
\newcommand{\Eup}{\Ev^{\mathrm{upbr}}} 
\newcommand{\Elw}{\Ev^{\mathrm{lwbr}}} 
\newcommand{\com}{\complement}
\newcommand{\GGamma}{\Upsilon}

\newcommand{\Bup}{\mathcal{B}^{\mathrm{upbr}}}
\newcommand{\Blwp}{\mathcal{B}^{\mathrm{lwbr,+}}}
\newcommand{\Blwm}{\mathcal{B}^{\mathrm{lwbr,-}}}
\newcommand{\kmax}{k_{\mathrm{max}}}
\renewcommand{\ahom}{{\overline{\mathbf{a}}}}

\newcommand{\Ind}{X}
\newcommand{\Err}{Y}

\DeclareMathOperator{\dist}{dist}

\DeclareMathOperator{\osc}{osc}

\DeclareMathOperator{\Cov}{Cov}
\DeclareMathOperator{\Var}{Var}

\DeclareMathOperator{\Law}{Law}

\numberwithin{equation}{section}
\numberwithin{figure}{section}

\newtheorem{theorem}{Theorem}[section]

\newtheorem{conjecture}[theorem]{Conjecture}

\newtheorem{lemma}[theorem]{Lemma}

\theoremstyle{remark}
\newtheorem{remark}[theorem]{Remark}
\theoremstyle{definition}

\allowdisplaybreaks
\makeatletter
\newcommand{\vast}{\bBigg@{4}}
\newcommand{\Vast}{\bBigg@{5}}
\makeatother

\begin{document}

\title{Tightness of the maximum of Ginzburg-Landau fields}
\author{Florian Schweiger}
\address[Florian Schweiger]{Section de mathématiques, Université de Genève, 7-9 rue du Conseil 
G\'{e}n\'{e}ral, 1205 Gen\`{e}ve, Switzerland.}
\email{florian.schweiger@unige.ch}
\author{Wei Wu}
\address[Wei Wu]{Department of Mathematics,  NYU Shanghai and and NYU-ECNU Math Institute,  Shanghai, China, 2001224.}
\email{wei.wu@nyu.edu}
\author{Ofer Zeitouni}
\address[Ofer Zeitouni]{Department of Mathematics, Weizmann Institute of Science, Rehovot 76100, Israel.}\email{ofer.zeitouni@weizmann.ac.il}

\date{\today}
\maketitle
\begin{abstract}
  We consider the discrete Ginzburg-Landau field with 
  potential satisfying a uniform convexity condition,
  in the critical dimension
  $d=2$, and prove that its maximum over boxes of sidelength $N$, centered
  by an explicit $N$-dependent centering, is tight.
\end{abstract}

\section{Introduction and main results}
Let $V\in C^{2}(\R)$ be symmetric around 0 and assume that there are constants $0<c_-\le c_+<\infty$ such that
\[c_-\le V''(x) \le c_{+}\]
Let $Q_N:=(-N,N)^2\cap\Z^2$ and $\partial^+ Q_N$ be the (outer) boundary of $Q_N$, that is the set of vertices in $\Z^2\setminus Q_N$ that are connected to $Q_N$ by some edge.
Then the Ginzburg-Landau Gibbs probability
measure on $Q_N$ with zero boundary condition is given by
\begin{equation}\label{e:defGLmeasure}
\PP^{Q_N,0}(d\phi)=\frac{1}{Z^{Q_N,0}}\exp \left(-\sum_{x\in Q_N\cup \partial^+ Q_N}\sum_{y\sim x}V\left(
\phi(y)-\phi(x) \right) \right) \prod_{x\in Q_N}d\phi (x)\prod_{x\in \partial^+Q_N}\delta
_{0}(\phi(x))
\end{equation}
where we extend the definition so that $\phi(x) = 0$ for $x\in \Z^2\setminus Q_N$, and $Z^{Q_N,0}$ is a suitable normalization constant. We denote by $\E^{Q_N,0}$ and $\Var^{Q_N,0}$
the expectation and variance with respect to the measure $\PP^{Q_N,0}$.

The Ginzburg-Landau measure  is a natural generalization of the discrete Gaussian free field, which corresponds to the special choice $V(x)=\frac{x^2}{2}$. The two-dimensional DGFF is a fundamental example of a log-correlated Gaussian field. The Ginzburg-Landau model, on the other hand, is in general not a Gaussian model, but it is still log-correlated. In fact, there is a constant $\g=\g(V)$ (the so-called effective stiffness) such that
\begin{equation}\label{e:covarianceGL}
\Cov^{Q_N,0}(\phi(x),\phi(y))\approx \g\log\frac{N}{1+|x-y|}.
\end{equation}
See \cite{NS97} for the corresponding infinite volume statement and \cite[Theorem 1]{M11} for the finite volume case. In these references it is also proved that as a random distribution
(that is, after integration against macroscopic test functions), the Ginzburg-Landau field converges to the (continuous) DGFF.

In view of this, the Ginzburg-Landau measure  is expected to
behave also in terms of its extremes like a multiple of the DGFF. It is known that the centered maximum of the DGFF
 converges in distribution to a randomly shifted Gumbel distribution \cite{BDZ16}, and so it is plausible that the same holds true for the Ginzburg-Landau  field. More precisely, definining
\begin{equation}\label{e:defmN}
m_N:=2\log N-\frac{3}{4}\log\log N,
\end{equation}
we expect that
\begin{equation}\label{e:tightness_maximum}
\max_{x\in Q_N}\phi(x)=\sqrt{\g}m_N+O(1).
\end{equation}
Previously it was shown that the first-order term is correct \cite{BW20} and that there is a (deterministic) subsequence along which the maximum centered by its mean
is tight \cite{WZ19}.

\subsection{Main result}
In this paper we prove the tightness of the centered maximum, that is
\eqref{e:tightness_maximum} under a mild additional regularity assumption on $V$.  Stated precisely, our main result is the following.
\begin{theorem}\label{t:main}
\label{main}
Let the sequence of measures $\PP^{Q_N,0}$ be as defined in \eqref{e:defGLmeasure},  and suppose that $V''$ is Lipschitz-continuous. Let $\phi_N$ be a sample from $\PP^{Q_N,0}$. Then there is a deterministic constant $\g=\g(V)$, such that with $m_N$ as in \eqref{e:defmN},  for every $\eps >0$,  there exists $C_\eps<\infty$,  such that for all $N\in\N$,
\begin{equation}\label{e:main_precise}
\PP^{Q_N,0}\left(\left|\max_{x\in Q_N}\phi(x)-\sqrt{\g}m_N\right|\le C_\eps\right) >
1-\eps.
\end{equation}
\end{theorem}
The constant $\g$ will turn out to be the same as in \eqref{e:covarianceGL} (and in fact as part of our proof we will establish a precise version of \eqref{e:covarianceGL}). The assumption that $V''$ is Lipschitz-continuous is needed as we will use results from \cite{M11,AW22,W22}. Presumably this assumption
could be relaxed.

\subsection{Earlier work}
Motivated by statistical mechanics, the Ginzburg-Landau Gibbs measure, as a
generalization of the Gaussian free field,
has been the focus of extensive studies for several decades. Important
milestones include the introduction of the Brascamp-Lieb inequalities
\cite{BL76}, the introduction of the Helffer-Sj\"{o}strand
dynamical representation
\cite{HS}, and its use by Naddaf-Spencer
to derive CLTs and decay of correlation \cite{NS97}. These in turn were
used in the construction of infinite volume limits,
hydrodynamics, and large deviations. See \cite{DGI00,FS97,GOS01} and the
lecture notes \cite{Fu05}; see also \cite{Sh05}.

In recent years and with the advent of quantitative homogenization,
refined CLT's, hydrodynamic limits
and large deviation estimates were obtained.
Of particular relevance to the Ginzburg-Landau measure are \cite{AW22,Da22,Da22a,AD23} and the local CLT in \cite{W22}.

The interest in the maximum of the Ginzburg-Landau field is partially
motivated by the entropic repulsion question, see \cite{DG00}, which in turn
is motivated by analogous results in the Gaussian setup \cite{BDG01}.
The latter paper motivated in parts the study of the maximum of the Gaussian
free field in the critical dimension $2$, culminating with the
proof of convergence in \cite{BDZ16} and its (Gaussian) generalization
in \cite{DRZ17}. For the Ginzburg-Landau model, as mentioned above, the
best results prior to this work are the evaluation of the leading order
asymptotics in \cite{BW20}, and the proof of sequential tightness in
\cite{WZ19}. We emphasize that the proof in 
\cite{WZ19} is not quantitative, does not seem to generalize to full tightness,
and further provides no information on
the centering term $m_N$ beyond the leading order.

\subsection{Outline of the proof}
When studying the extrema of log-correlated fields, the standard approach in the literature is to relate the extrema of the field to the extrema of a branching structure such as a branching random walk. This branching makes then first- and second-moment calculations feasible. 
We refer the reader to \cite{Bi20,Ze16} for a pedagogical exposition in the Gaussian setup (where various arguments are simpler) and to  \cite{BRZ20} for a sample application in the non-Gaussian setup.

In the case of the Ginzburg-Landau field, this asymptotic branching structures will be given in terms of a decomposition in terms of harmonic averages
introduced in \cite{BW20}. We fix the sequence of radii $r_k=N\e^{-k}$ as in \eqref{eq-rk}. Now for any point $x\in Q_N$, we would like to decompose
\[\phi(x)=\sum_{k=0}^K
(\Gamma_{Q_{r_{k+1}}(x)}(x,\phi)-\Gamma_{Q_{r_{k}}(x)}(x,\phi)),\] where $Q_{r_k}(x)$ is the box of side $r_k$ centered
at $x$ and $\Gamma_{Q_{r_{k}}(x)}(x,\phi)$ is a weighted average of the values of $\phi$ on $\partial^+\Q_{r_k}(x)$ as defined in \eqref{eq-harmonic}. For this random walk approximation to be useful, we would like the increments $\Gamma_{Q_{r_{k+1}}(x)}(x,\phi)-\Gamma_{Q_{r_{k}}(x)}(x,\phi)$ to be independent and with some good control on their distribution. For the DGFF, these increments turn out to be i.i.d Gaussians (the independence follows immediately from the Gauss-Markov property). However for the
Ginzburg-Landau field, these increments are neither independent nor do we have control over their distribution.

The latter problem can be solved by modifying our definition of the increments: Inside of considering $\Gamma_{Q_{r_{k}}(x)}(x,\phi)$, we consider $\Sav_k(x,\phi)$, defined as a radially smoothed out average of $\Gamma_{Q_{r}(x)}(x,\phi)$ for $r\sim r_k$ so that we can use the CLT for the Ginzburg-Landau field. In view of the quantitative CLT from \cite{AW22,W22} we do not need to smooth out over length $r_k$, but only over $r_k^{1-\eps}$ for some small $\eps>0$. This gain is very important for our goal, as it ensures that various error terms (in particular the boundary layer to be discussed shortly) are summable.

The former problem, that the increments $\Inc_k(x,\phi):=\Sav_{k+1}(x,\phi)-\Sav_k(x,\phi)$ need not be independent, is much more troublesome.
 However, due to 
Miller's coupling \cite{M11} reviewed in  Section \ref{subsec-miller}, this independence is 
``almost'' true,  as long as one is ready to make various concessions: First, the precise notion of independence is obtained by introducing independent variables and coupling with them on a ``good event'' $\Ecp_{k,x}$.  Second, this coupling is only possible if the boundary values of $\phi$ on $\partial^+Q_{r_k}(x)$ are not too rough. Third and last, we need at least mesoscopic separation between the different increments, and so need to introduce boundary layers separating them. These issues force us to introduce more ``good events'' $\Ero_{k,x}$ and $\Ebd_{k,x}$ that we need to occur. When these do not occur,
 we need to resort to using rather wasteful Brascamp-Lieb estimates to retain any control over the field, and so we need to use that points where a good event at an early scale fails are rare enough not to affect the maximum (a similar difficulty appeared in \cite{SZ24}).

Our construction of the random walk approximation is inspired by \cite{BW20}. However,  since we require a much higher precision, our construction is different. In particular, since the probability that the good coupling event $\Ecp_{k,x}$ fails is relatively large (of order $r_k^{-\delta}$,  for some small $\delta>0$), we need to continue with the coupling at finer scales even if one coupling has failed. The construction present in this paper couples all boundary conditions at once,  couples the increments for a large set of representative points in $Q_N$ altogether,  and proceeds iteratively. We explain the construction in details in 
Section \ref{sec-RWapp}.

With the random walk approximation in place, the proof proceeds differently for the upper and lower bounds. For the upper bound, in order to obtain
the correct centering term $m_N$  of \eqref{e:defmN} and in particular the $\log\log N$ term there, it is crucial to proceed iteratively and show that partial sums of approximate harmonic averages
do not exceed a certain ``barrier''; this is the barrier event 
$\Eup_x$ introduced in \eqref{eq-upperbarrier}. The main step in introducing 
the barrier is to prove some continuity of the harmonic averages at intermediate scales, see Lemma \ref{l:fluctuations_uncond}, whose proof builds on a chaining argument. (We note that for $k$ large, the chaining has to be done conditionally, see Lemma \ref{l:fluctuations_cond}, 
which adds an extra technical complication.) Another complication in
handling the barrier is that the estimates on the probability
of the good coupling event 
$\Ecp_{x,k}$ are not strong enough to show that it holds for all relevant $x$ and all $k$. Therefore, when introducing the barrier, some couplings will fail, and we handle those increments, as in \cite{BW20}, by a Brascamp-Lieb argument. We
also make use of the strong local CLT from \cite{W22}, in the form of
Lemma \ref{l.meantilt}. The argument that allows one to insert the barrier
is based on the ``good'' coupling event, and uses the barrier estimates
given in 
Section \ref{s:ballot_thms}; those in turn are versions of the classical ballot theorem, but given for walks with independent but not identically distributed
variables. 
Once the barrier is in place, the 
upper bound follows by a simple union bound.

For the lower bound, we proceed by the classical first-second moment
argument; we only consider points $x$ for which the coupling is consistently succesful at all scales, and such that an appropriate barrier condition 
$\Elw_x$, see
\eqref{eq-Elw},
is satisfied. (This barrier event is more restrictive than the one used in 
the upper bound.)
The proof then uses barrier estimates (a lower bound for the
first moment, and an upper bound for the second moment). The details
are provided in Section \ref{s:lower}.

\subsection{A conjecture}
In view of Theorem \ref{t:main}, it is natural to expect convergence
of the maximum to a randomly shifted Gumbel. Formally, we expect the following.
\begin{conjecture}
  \label{conj-conv}
  With notation as in Theorem \ref{t:main}, 
  it holds that for  $c^*=2/\sqrt{\g}$ and some random variable $Z>0$,
  \begin{equation}\label{econj:main_precise}
\lim_{N\to\infty}\PP^{Q_N,0}\left(\max_{x\in Q_N}\phi(x)-\sqrt{\g}m_N\leq t\right)
=\E\left( \e^{- Z\e^{-c^*t}}\right),\quad\forall t\in\R.
\end{equation}
\end{conjecture}
To obtain Conjecture \ref{conj-conv} would require a more precise control than ours on the bad coupling events, and a more careful analysis of the small scale 
differences of harmonic averages, and in particular non-centered ones; see
\cite{BDZ16} for a a simpler setup where this program is implemented.

\subsection*{Acknowledgement}
This work was partially supported by ISF grant number 421/20. FS was supported by the Foreign Postdoctoral Fellowship Program of the Israel
Academy of Sciences and Humanities, the NCCR SwissMAP, the Swiss FNS, and the Simons collaboration on localization of waves.
The work of WW was partially supported by MOST grant 2021YFA1002700,  NSFC grant 20220903 NYTP and SMEC grant 2010000080.
OZ thanks Cheuk Ting Li for a useful discussion concerning \cite{AL19} and the coupling
argument in Lemma \ref{l:simple_coupling}.

\section{Set-up and important tools}

\subsection{Notation}
Throughout the paper we will work under the standing assumption that $V''$ is Lipschitz-continuous and satisfies
\[c_-\le V''(x) \le c_{+}.\]

We first generalize \eqref{e:defGLmeasure} to non-zero boundary values and general domains. So let $D\subset\Z^2$, and $f\colon\partial^+D\to\R$. Then the
Ginzburg-Landau measure on $D$ with Dirichlet boundary condition $f$ is given by
\begin{equation}\label{e:defGLmeasure_f}
\PP^{D,f}(d\phi)=\frac{1}{Z^{D,f}}\exp \left(-\sum_{x\in D\cup \partial^+ D}\sum_{y\sim x}V\left(
\phi(y)-\phi(x) \right) \right) \prod_{x\in D}d\phi (x)\prod_{x\in \partial^+D}\delta
_{0}(\phi(x)-f(x)).
\end{equation}
We continue to denote a generic random variable distributed according to $\PP^{D,f}$ by $\phi$. When there is risk of confusion, we add superscripts.

For a function $g\colon D\to\R$ we define
\[\langle\phi,g\rangle:=\sum_{x\in D}\phi(x)g(x).\]

We let $Q_R(x)=x+(-R,R)^2\cap\Z^2$ denote the lattice cube of sidelength $2R$ centered at $x$. When $x=0$ we generally omit it.
As we will often encounter the complement $(Q_R(x))^\com$ of a cube $Q_R(x)$, we use the slightly shorter notation $Q_R^\com(x)$ instead. 

For a set $D'\subset D$, we let $\F_{D'}=\sigma(\{\phi(x)\colon x\in D'\})$.  For $f: D \to \R$,  we denote by $f\restriction_{D'}$ the restriction of $f$ on $D'$.

The letters $C$ and $c$ denote generic constants (that might in particular depend on $V(\cdot)$), and whose precise value might change from line to line.

\subsection{Brascamp-Lieb inequalities and rough upper bounds}
Our assumption $V''\ge c_-$ ensures that the measure $\PP^{Q_N,0}$ is strictly log-concave.
Thus, the Brascamp-Lieb inequalities \cite{BL76} apply to this measure and give bounds on the (polynomial and exponential) moments of $\PP^{Q_N,0}$.

\begin{lemma}[Brascamp-Lieb inequalities] \label{l:BL}
Let $\mathbb{E}_{\mathrm{DGFF}}^{D,f}$ and $\Var_{\mathrm{DGFF}}^{D,f}$ denote the expectation and variance with respect to the discrete GFF measure (i.e. \eqref{e:defGLmeasure} with $V(x)=\frac{x^2}{2}$). Then for any $g\colon D\to\R$ we have
\begin{align}
\Var^{D,f}\langle \phi ,g\rangle \quad & \le \frac{1}{c_-} \Var_{\mathrm{DGFF}}^{D,f}\langle \phi ,f\rangle ,
\label{e:BLvarbound} \\
\E^{D,f}\left( \exp \left(\langle \phi ,g\rangle -
\E^{D,f}\langle \phi ,g\rangle \right) \right) \quad
& \le\exp \left( \frac{1}{2c_-}\Var_{\mathrm{DGFF}}^{D,f}\langle \phi ,f\rangle\right).  \label{e:BLexpmombound}
\end{align}
\end{lemma}

The exponential moment bound \eqref{e:BLexpmombound} directly implies a tail bound on $\phi(x)$ for $x\in Q_N$.

\begin{lemma}
\label{l:BLtailbound} There is $C>0$ such that for any $N\in\N$ and any $x\in Q_N$
\begin{equation}
\PP^{Q_N,0}\left( \phi(x)\ge t\right) \le C\exp\left(-\frac{c_-\sqrt{\pi}}{2\sqrt{2}\log\dist(x,\partial^+ Q_N)}t^2\right).\label{e:BLtailbound}
\end{equation}
\end{lemma}
This is essentially \cite[Lemma 2.2]{BW20}. However, in the statement of that Lemma a prefactor $C$ and a $\log$ in the denominator are missing, and so we give the corrected result with its proof. 
\begin{proof}
It is well-known that
\[\left|\Var_{\mathrm{DGFF}}^{Q_N,0}\phi_N(x)-\sqrt{\frac{2}{\pi}}\log\dist(x,\partial^+Q_N)\right|\le C\]
for any $N$ and any $x\in Q_N$. Combining this with \eqref{e:BLexpmombound} for $g=s\delta_x$ for some $s\in\R$ we find
\[\PP^{Q_N,0}(\phi_N(x)\ge t)\le \e^{-st}\E^{Q_N,0}(\e^{s\phi_N(x)})\le \exp\left(-st+\frac{s^2}{2c_-}\sqrt{\frac{2}{\pi}}\log\dist(v,\partial^+Q_N)+C\right).\]
Choosing $s=\frac{c_-t}{\sqrt{\frac{2}{\pi}}\log\dist(x,\partial^+Q_N)}$ yields the claimed result.
\end{proof}

\subsection{Miller's coupling}
\label{subsec-miller}
We next state an extremely useful coupling result due to Miller \cite{M11}.

For $D\subset\Z^2$ and $r\ge0$ we let $D^r=\left\{ v\in D\colon\dist(v,\partial
D)>r\right\} $. 

\begin{theorem}[{\cite[Theorem 1.2]{M11}}]
\label{t:miller_coupling} Let $\Lambda>0$ be a constant. Then there are constants $c_\Lambda$, $\alpha_\Lambda$, $\beta_\Lambda$ with the following property.

Let $D\subset \mathbb{Z}^{2}$,  $R:= \max \{|x-y|: x,y\in D\}$,  and let $r\ge c_\Lambda R^{1-\alpha_\Lambda}$.
Fix some $f\colon\partial^+ D\rightarrow \R$ such that $\osc_{x\in \partial D}|f|\le 2\Lambda |\log R|^\Lambda$. There exists a coupling of $\PP^{D,f}$ and $\PP^{D,0}$ such that
$\phi^{D,0} $ and $\phi^{D,f}$ are distributed according to $\PP^{D,0}$ and $\PP^{D,f}$,
respectively, and further, if
 $\hat{\phi}\colon D^r\to\R$ is discrete harmonic with $\hat{\phi}
\restriction_{\partial^+ D^{r}}=\left(\phi^{D,f}-\phi^{D,0}\right)\restriction_{\partial^+D^{r}}$, then
\[
\PP\left( \phi^{D,f}=\phi^{D,0}+\hat{\phi} \text{ in }D^{r}\right) \ge1-\frac{c'_\Lambda}{R^{\beta_\Lambda}}.
\]
\end{theorem}
We note that in \cite[Theorem 1.2]{M11} the assumption on $f$ is that $\max_{x\in \partial D}|f|\le \Lambda |\log R|^\Lambda$. But the statement of the theorem is invariant under adding a constant to $f$, so the version stated here directly follows.

We also have the following estimate on $\E^{D,f}\phi$.

\begin{theorem}[{\cite[Theorem 1.3]{M11}}]
\label{t:miller_mean} Let $\Lambda>0$ be a constant, and let $c_\Lambda$, $\alpha_\Lambda$, $\beta_\Lambda$, $D$ and $f$ be as in Theorem \ref{t:miller_coupling}. Then, if $r\ge c_\Lambda R^{1-\alpha_\Lambda}$ and if $\hat\phi\colon D^r\to\R$ is discrete harmonic function  with $\hat{\phi}\restriction_{\partial^+ D^{r}}=\E\phi ^{f}\restriction_{\partial^+ D^{r}}$, then
\[
\max_{x\in D^{r}}\left|\E^{D,f}\phi(x) -\hat{\phi}(x)\right|\le \frac{c'_\Lambda}{R^{\beta_\Lambda}}.
\]
\end{theorem}

For our purposes we need a slightly modified version of Theorem \ref{t:miller_coupling}, where we make two additional claims. The first is that if we are only interested in coupling the two fields in an annulus instead of the whole domain, then the event that the coupling succeeds can be chosen to be measurable with respect to that annulus. The second is that we can obtain the coupling not just for some fixed boundary data $f$, but simultaneously for all of them, 
in a measurable way. 

To be precise, we claim the following.
\begin{theorem}
\label{t:miller_coupling_annulus} Let $\Lambda>0$ be a constant. Then there are constants $c_\Lambda$, $\alpha_\Lambda$, $\beta_\Lambda$ with the following property.

Let $D\subset \mathbb{Z}^{2}$, and let $r\ge c_\Lambda R^{1-\alpha_\Lambda}$. Moreover let
\[F=\left\{f\colon\partial^+ D\rightarrow \R\middle|\osc_{x\in \partial ^+D}|f|\le 2\Lambda |\log R|^\Lambda\right\}\]
and let $f\in F$. Let $D'\subset D^r$. Then there exists a coupling 
$\PP^{D,D',f,0}$ of \,  $\PP^{D,f}$ and $\PP^{D,0}$
and an event $\Ev_f\in \F_{D\setminus D'}\times\F_{D\setminus D'}$  such that
 $\phi^{D,0} $ and $\phi^{D,f}$ are  distributed
according to $\PP^{D,0}$ and $\PP^{D,f}$,
respectively, and the following properties hold.
%
If $\hat{\phi}\colon D^r\to\R$ is discrete harmonic with $\hat{\phi}
\restriction_{\partial^+ D^{r}}=\left(\phi^{D,f}-\phi^{D,0}\right)\restriction_{\partial^+D^{r}}$, then on the event $\Ev_f^\com$ we have
\begin{equation}\label{e:miller_coupling_annulus1}
\phi^{D,f}=\phi^{D,0} +\hat{\phi}\text{ in }D^{r}\setminus D'.
\end{equation}
Moreover,
\begin{equation}\label{e:miller_coupling_annulus2}
\PP^{D,D',f,0}(\Ev_f) \le\frac{C_\Lambda}{R^{\beta_\Lambda}}.
\end{equation}

Finally this coupling is regular in the following sense: The mapping from $F$ to probability measures on $(\R^{D})^2$ defined by $f\mapsto \PP^{D,D',f,0}$ is a Markov kernel and the set
\begin{equation}\label{e:miller_coupling_annulus3}
  \bigcup_{f\in F}(\{f\}\times\Ev_f)
  \subset\R^{\partial^+D}\times (\R^{D})^2
\end{equation}
is (Borel-)measurable.
\end{theorem}

\begin{proof}
Consider first some fixed boundary data $f$. Theorem \ref{t:miller_coupling} implies that
\[d_{TV}\left(\Law \left(\phi^{D,f}\restriction_{D^r}\right),\Law\left(\left(\phi^{D,0}+\hat\phi\right)\restriction_{D^r}\right)\right)\le \frac{c'_\Lambda}{R^{\beta_\Lambda}}.\]
Integrating out some coordinates, we obtain in particular that
\[d_{TV}\left(\Law \left(\phi^{D,f}\restriction_{D^r\setminus D'}\right),\Law\left(\left(\phi^{D,0}+\hat\phi\right)\restriction_{D^r\setminus D'}\right)\right)\le \frac{c'_\Lambda}{R^{\beta_\Lambda}}.\]
and this in turn means that $\phi^{D,f}$ and $\phi^{D,0} +\hat{\phi}$ can be coupled so that they agree on $D^{r}\setminus D'$ with probability at least $1-\frac{c'_\Lambda}{R^{\beta_\Lambda}}$. In particular, if we define $\Ev_f^\com$ as the event that \eqref{e:miller_coupling_annulus1} holds, then $\Ev_f$ is $\F_{D\setminus D'}\times\F_{D\setminus D'}$-measurable and \eqref{e:miller_coupling_annulus2} holds trivially.

It remains to prove the regularity of the coupling. The fact that $f\mapsto\Law\left(\phi^{D,f}(\cdot),\phi^{D,0}(\cdot)\right)$ defines a Markov kernel follows from the general Lemma \ref{l:coupling_measurability} below. Finally, the measurability of \eqref {e:miller_coupling_annulus3} is clear, because the set $\Ev_f=\left\{(\phi^{D,f},\phi^{D,0})\colon \phi^{D,f}=\phi^{D,0} +\hat{\phi}\text{ in }D^{r}\setminus D'\right\}\subset(\R^{D})^2$ does not depend on $f$ (and is closed and therefore measurable) 
\end{proof}

We state the lemma on coupling of families of probability measures that we used above.
\begin{lemma}\label{l:coupling_measurability}
Let $F$ be a measurable space, and $\nu^f$ for $f\in F$ be a family of probability measures on some Polish space $E$ such that $f\mapsto\nu^f(\cdot)$ is a Markov kernel. Let $\bar\nu$ be another probability measure on $E$, and suppose that
\[\sup_{f\in F}d_{TV}(\nu^f,\bar \nu)\le \varepsilon.\]
Then if $X^f$ and $\bar X$ are distributed according to $\nu^f$ and $\bar \nu$, respectively, there is a coupling $\PP^f$ of $\Law(X^f),\Law(\bar X)$ such that $\PP^f(X^f\neq \bar X)\le\varepsilon$ and such that
$f\mapsto\PP^f$ is a Markov kernel.
\end{lemma}
This lemma is most likely well-known, but we could not locate a reference, and so we sketch a proof. 
\begin{proof}
For a fixed $f$ this is the standard coupling lemma (see e.g. \cite[Theorem 5.2]{Lin92}). In the following we follow the proof given there while verifying at each step the measurability.

Note that $\nu^f$ and $\bar \nu$ are absolutely continuous with respect to $\nu^f+\bar \nu$. Let $g^f$ and $\bar g^f$ denote the respective densities. Importantly, both $g^f$ and $\bar g^f$ are jointly measurable functions $F\times E\to\R$ (as follows e.g. from \cite[Theorem 58]{DM82}).

Furthermore let $d\sigma^f=(g^f\wedge\bar g)d(\nu^f+\bar \nu)$. This is a subprobability measure on $E$, and it is a sub-Markov kernel in the sense that $f\mapsto\sigma^f(\Ev)$ is measurable for any event $\Ev$.

In particular the set $T=\left\{f\colon \int d\sigma^f=1\right\}$ is measurable. On $T$ we have $\nu^f=\bar \nu$ and so we can take e.g. the diagonal map to obtain a trivial coupling. On $T^\com$ we instead consider the measure
\[\PP^f=\frac{1}{1-\int d\sigma^f}(\nu^f-\sigma^f)\times(\bar\nu-\sigma^f)+\beta_*\sigma^f\]
on $E\times E$, where $\beta_*\sigma^f$ is the pushforward of $\sigma^f$ under the diagonal embedding $\beta\colon E\to E\times E$. Then $f\mapsto \PP^f(\cdot)$ is a Markov kernel, and as in \cite[Theorem 5.2]{Lin92} one verifies that it yields a coupling with the desired properties.
\end{proof}

\subsection{The Helffer-Sj\"ostrand representation} \label{sec:HS}

In this section we summarize the observation
of \cite{NS97} (which  in turn was inspired by the works~\cite{HS}), that the variance of linear functionals of the Ginzburg-Landau field can be written in terms of the Helffer-Sj\"ostrand operator. The Helffer-Sj\"ostrand representation is crucial for the proof of
a quantitative central limit theorem for the statistics of $\phi$,  stated in the section below. It is also used in the proof of Miller's coupling, see
\cite{M11}.

The finite volume measure $\PP^{Q_N, 0}$ is invariant under the Langevin-dynamics
\begin{equation}
\label{e.dynamics}
\left\{
\begin{aligned}
& d\phi_t(x)
= \sum_{y\sim x} V'( \phi_t(y)-\phi_t(x)) \,dt  + \sqrt{2} \,dB_t(x), && x\in Q_N,
\\ &
\phi_t(x) = 0, && x \in \partial^+ Q_N,
\end{aligned}
\right.
\end{equation}
where $\{ B_t(x) \,:\, x\in Q_N \}$ is a family of independent Brownian motions.
Let $\omega_x: Q_N \to \R$ be defined by $\omega_x(y) = 1_{x=y}$.
The infinitesimal generator of this process
is the operator~$\Delta_\phi$ defined by
\begin{align*}
\Delta_\phi F (\phi)
: =
\sum_{x\in Q_N} \partial_x^2 F(\phi)
-
\sum_{x\in Q_N}
 \sum_{y\sim x}
V'(\phi(y)-\phi(x) )\partial_xF(\phi),
\end{align*}
where
\begin{equation*}
\partial_x F (\phi):= \lim_{h\to 0} \frac1h \left( F(\phi+h\omega_x) - F(\phi) \right).
\end{equation*}
We notice that the formal adjoint of~$\partial_x$ with respect to $\PP^{Q_N,0}$, which we denote as $\partial_x^*$, is given by
\begin{equation*}
\partial_x^* F := -\partial_xF
+
\sum_{y\sim x}
V'(\phi(y)-\phi(x) ) F(\phi).
\end{equation*}

We let~$\mathcal{E}(Q_N)$ denote the set of directed edges on $Q_N$, equipped with a lexicographical ordering  $\ll$ of $\Z^2$.
We define the gradient operator $\nabla$ acting on functions
$h: Q_N\to \R$ by 
\begin{equation} 
  \nabla h: \mathcal{E}(Q_N)\to \R, 
  \quad \nabla h(e) =h(y)-h(x) \; \mbox{\rm if $e=(x,y)$ with $x\ll y$.}
\end{equation}
We also define the discrete version of the negative of the divergence operator, for functions~$f:\mathcal{E}(Q_N)\to \R$ by
\begin{equation}
\left( \nabla^*\cdot f\right)(x)
:=
\sum_{y\sim x, \, y\ll x} f(y,x) 
- \sum_{y\sim x, \, x\ll y} f(x,y), 
\quad x\in Q_N.
\end{equation}

Define the Helffer-Sj\"ostrand operator $\mathcal L:= -\Delta_\phi + \nabla^*\cdot V'' \nabla$, such that for functions 
$\hat F: Q_N\times \R^{Q_N}\to \R$,
$$
\mathcal L \hat F(x,\phi) :=  -\Delta_\phi  \hat F + \sum_{y\sim x} V''(\phi(y) -\phi(x)) (\hat F(y,\phi)- \hat F(x,\phi)).
$$
 Notice that
$\mathcal L$ is the generator for the Markov process $(X_t,\phi_t)$, where $\phi_t$
is the Langevin dynamics \eqref{e.dynamics} and $X_t$ is a continuous time random walk in $Q_N$
with jump rates $V''(\phi_t(y)- \phi_t(x))$, stopped when hitting the boundary.

The following representation of the variance is obtained by Naddaf and Spencer
in \cite{NS97} .

\begin{theorem}[Helffer-Sj\"ostrand representation]
\label{p.HS}
For all $F$ such that
\begin{equation*}
\E^{Q_N,0} \biggl(F(\phi)^2 + \sum_{x \in Q_N} (\partial_x F(\phi))^2 \biggr) < \infty,
\end{equation*}
we have
\begin{equation*}
\Var^{Q_N,0}(F(\phi)) = \left\langle \partial F,  \mathcal L^{-1}  \partial F  \right\rangle,
\end{equation*}
where $\left\langle \partial F, \mathcal L^{-1}  \partial F  \right\rangle := \sum_{x,y \in Q_N} \E^{Q_N,0}(\partial_x F \mathcal L^{-1} _{xy} \partial_y F)$.
\end{theorem}

If we consider a linear statistics of $\phi$, and take $F(\phi) =  \sum_{x \in Q_N} f(x) \phi(x)$ for some test
function $f$, then the above proposition implies
\begin{equation}
\label{e.HSlin}
\Var^{Q_N,0}\left(\sum_{x \in Q_N} f(x) \phi(x)\right) =  \left\langle f, \mathcal L^{-1} f  \right\rangle.
\end{equation}
\subsection{Quantitative CLT}

On large scales we expect the Ginzburg-Landau model to behave like a DGFF. One way to formalize this is via a central limit theorem (CLT). In a qualitative form, such a CLT is due to \cite{NS97,M11}. For our purposes, however, it is crucial to have a quantitative version of the CLT. Such a result was recently shown in \cite{W22}, building upon the recent breakthrough on quantitative stochastic homogenization of the Helffer-Sj\"{o}strand equation in \cite{AW22}.  

Define the volume-normalized versions of the $L^2$ norm,  for every $U\subset \Z^2$ and $f: \mathcal E(U) \to \R$, by 
\begin{equation*}
\left\| f \right\|_{\underline{L}^2(U)}
:=
\left(
\frac1{|U|}
\sum_{e\in \mathcal E(U) } | f(e) |^2 \right)^{\frac12}.
\end{equation*}

\begin{theorem}[{\cite[Theorem 4.1]{W22}}]
\label{t:qclt}
Let $R>0$, $\alpha\in(0,\frac 12)$ and  consider the measure $\PP^{Q_R,0}$.
Let $g\colon Q_R \to \R$ and suppose that there is a H\"{o}lder continuous function $\tilde g\in C^{1,\alpha} ([-1,1]^2)$ such that $g(x)=\tilde g\left(\frac{x}{R}\right)$ for all $x\in Q_R$. Moreover,  we assume that $g$ has mean zero: $\sum_{x\in Q_R}g(x) =0$.  Let $f: \mathcal E(Q_R) \to \R$ be such that $g  = \nabla^* \cdot f$.

Then for any $0<\kappa\le1$ there are constants $\overline{\A}= \overline{\A}(V)>0$, $C$ and $\gamma>0$ (depending on $V$ and $\kappa$ only) such that
\begin{equation}\label{e:qclt}
\left|\frac{1}{R^{4}}\Var^{Q_R,0}\langle \phi,g\rangle-\frac{1}{\overline{\A}}\int_{(-1,1)^2\times (-1,1)^2}\tilde g(x)G^{-\Delta}(x,y)\tilde g(y)dxdy\right|
\le
C\frac{\| \tilde g\|_{C^{1,\kappa}}^2+ \| f\|_{\underline {L}^2(Q_R)}^2}{R^{\gamma}},
\end{equation}
where $G^{-\Delta}$ denotes the Green's function of the Laplacian on $[-1,1]^2$.
\end{theorem}
Note that \cite[Theorem 4.1]{W22} is stated only for gradient observables. However, one can recover our result easily from the one in \cite{W22} by using that every function with mean 0 can be written as the divergence of some vector field $f$.    Since $f(u,v) = \sum_{y\in Q_R} g(y) (G_{Q_R} (u,y) -G_{Q_R}(v,y))$,  where $G_{Q_R}$ is the Dirichlet Green's function in $Q_R$,  and $g(x)=\tilde g\left(\frac{x}{R}\right)$, we see that $f$ is also a slow varying function as in \cite[Theorem 4.1]{W22}.  The scaling changes from $R^{-2}$ to $R^{-4}$ since we study the linear functional of $\phi$,  instead of the gradient variables.

In our application the functions $\tilde g$ for different $R$ will depend on $R$, but they will be close in $H^{-1}$ to a deterministic function. This will allow to replace the $R$-dependent integral term in \eqref{e:qclt} by a constant term up to a small error.

The next lemma estimates the derivative of the moment generating function,  which will be applied to estimate the mean under the measure with a linear tilt.

\begin{lemma}
\label{l.meantilt}
Under the  assumptions and notation  of Theorem
\ref{t:qclt},
 define for $t\in\R$
\[G(t) := \E^{Q_R,0} \exp\left(t R^{-2}\langle \phi, g\rangle\right).\]
Then, with $\overline{\A}= \overline{\A}(V)>0$  the constant from Theorem \ref{t:qclt},  for any $0<\kappa\le1$ there are constants $C<\infty$, $c_1>0$ and $\gamma>0$ (depending on $V$ and $\kappa$ only) such that for every $|t| \le c_1 \sqrt{\log R}$,
\begin{align}
\label{e.meantilt}
 \left|\frac{d}{dt}G(t) - t\sigma^2 G(t)\right|
 \leq
 C\frac{\| \tilde g\|_{C^{1,\kappa}}^2+  \| f\|_{\underline {L}^2(Q_R)}^2}{R^{\gamma/4}},
  \end{align}
  where
  \begin{align*}
  \sigma^2:= \frac{1}{\overline{\A}}\int_{(-1,1)^2\times (-1,1)^2} \tilde g(x)G^{-\Delta}(x,y) \tilde g(y)dxdy.
   \end{align*}
\end{lemma}
\begin{proof}
Since $g$ has mean zero,  we may write $g= \nabla^* \cdot f$ for some $f:\mathcal E(Q_R )\to \R$.  By the Brascamp-Lieb inequality  we have for every $t\in\R$,  the bound $G(t) \le \e^{Ct^2}$,  for some $ C=C(f)$.
Let $u$ be the solution of
\begin{equation}
\label{e.HS}
\left\{
\begin{aligned}
&
-\Delta_\phi u + \nabla^*\cdot \A \nabla u
 =  \nabla^* \cdot f_R
\quad &\mbox{in} & \ Q_R\times \Omega_0(Q_R)
\\ &
u_R = 0& \mbox{on} & \ \partial Q_R \times\Omega_0(Q_R),
\end{aligned}
\right.
\end{equation}
where $\A(e) = V''(\nabla \phi(e))$,  and
\begin{equation*}
f_R(e):= R^{-2} f\left(e \right),
\end{equation*}
satisfies $\partial_x R^{-2} \langle  \phi,  g \rangle = \nabla^*\cdot f_R$.
Notice that $u$ satisfies
\begin{equation}
\label{e.vfirstorder}
\sum_{x\in Q_R} \partial_{x}^* u(x,\cdot)
= R^{-2} \langle \phi,  g \rangle.
\end{equation}
This can be justified by applying $\partial_x$ to the equation above, which yields \eqref{e.HS}, and use the uniqueness of the solution to \eqref{e.HS}.

Using~\eqref{e.vfirstorder} and  integration by parts, we compute
\begin{align*}
&G'\left( t\right)
=
\E^{Q_R,0}\left(
R^{-2} \langle  \phi,  g \rangle
\exp\left( tR^{-2} \langle  \phi,  g \rangle \right)
\right)
\\ &
=
\E^{Q_R,0}\left(
\sum_{x\in Q_R} \partial_{x}^* u(x,\cdot)
\exp\left( tR^{-2} \langle  \phi,  g \rangle \right)
\right)
=
t \E^{Q_R,0}\left(
\sum_{x\in Q_R}  u(x,\cdot) \nabla^*\cdot f_R(x)
\exp\left( tR^{-2} \langle  \phi,  g \rangle\right)
\right).
\end{align*}
We apply \cite[Theorem 4.4]{W22} to conclude that
\begin{align*}
\E^{Q_R,0} \left(
\left|  \sum_{x\in Q_R}  u(x,\cdot) \nabla^*\cdot f_R(x)- \int  \bar u(x) \tilde g(x)\,dx \right|^2
\right)
\leq
CR^{-\gamma}\left( \| \tilde g\|_{C^{1,\alpha}}^2+ \| f\|_{\underline {L}^2(Q_R)}^2\right),
\end{align*}
where $\bar u$ is the solution to the constant coefficient equation
\begin{equation}
\label{e.homog}
\left\{
\begin{aligned}
&
-\nabla^* \cdot \ahom \nabla \bar u = \tilde g \quad &\mbox{in}& \ [-1,1]^2
\\ &
\bar u = 0& \mbox{on} & \ \partial ([-1,1]^2).
\end{aligned}
\right.
\end{equation}
Therefore an integration by parts yields
\begin{equation*}
\int  \bar u(x) \tilde g(x)\,dx
=
 \frac{1}{\overline{\A}}\int_{[-1,1]^2\times [-1,1]^2}   \tilde g(x)G^{-\Delta}(x,y) \tilde g(y) dx dy
 := \sigma^2.
 \end{equation*}
 Applying the Brascamp-Lieb inequality we obtain that, for every~$t\leq c_1\sqrt{\log R}$ with $c_1$ sufficiently small, there exists $C<\infty$ such that
\begin{align*}
\left| G'\left( t\right)  - t \sigma^2 G(t)  \right|
&
\leq
t \left|  \E^{Q_R,0}\left(
\left( \sum_{e}  u(e,\cdot) \nabla^*\cdot f_R(e) -\sigma^2 \right)
\exp\left( tR^{-2} \langle  \phi,  g \rangle \right)
\right) \right|
\\ &
\leq
t \E^{Q_R,0}\left(
\left| \sum_{e}  u(e,\cdot) \nabla^*\cdot f_R(e) -\sigma^2 \right|^2 \right)^{1/2}
\E^{Q_R,0}\left( \exp\left( 2tR^{-2} \langle  \phi,  g \rangle\right)
\right)^{1/2}
\\ &
\leq
C R^{-\gamma/2}\left( \| \tilde g\|_{C^{1,\kappa}}^2+  \| f\|_{\underline {L}^2(Q_R)}^2\right) t \exp\left( Ct^2 \right)
\\ &
\leq C\left( \| \tilde g\|_{C^{1,\kappa}}^2+  \| f\|_{\underline {L}^2(Q_R)}^2\right) R^{-\gamma/4}.
\end{align*}
\end{proof}

By integrating \eqref{e.meantilt} over $t$,   we conclude the following quantitative CLT.

\begin{lemma}\label{l:qclt_exponentialmoment}  
With assumptions and notation as in Theorem \ref{t:qclt},
there are constants $C<\infty$, $c_1>0$ and $\gamma>0$ (depending on $V$ and $\kappa$ only) such that for every $|t| \le c_1 \log^{1/2} R$
\begin{align}
 \left|G(t) - \exp\left(\frac 12 \sigma^2 t^2\right)\right|
 \leq
 C|t| \frac{ \| \tilde g\|_{C^{1,\kappa}}^2+  \| f\|_{\underline {L}^2(Q_R)}^2}{R^{\gamma}},
  \end{align}
  where
  \begin{align*}
 \sigma^2:= \frac{1}{\overline{\A}}\int_{[-1,1]^2\times [-1,1]^2} \tilde g(x)G^{-\Delta}(x,y) \tilde g(y)dxdy.
   \end{align*}
\end{lemma}

\section{Random walk approximation}
\label{sec-RWapp}
The goal of this section is to define an approximation of $\phi(x)$ for $x\in Q_N$ via a random walk with approximately independent and approximately Gaussian increments. These increments will be given as differences of square averages of $\phi$ (a small modification of the circle averages used in \cite{BW20}), and applying Theorem \ref{t:miller_coupling_annulus} and Theorem \ref{t:qclt} iteratively from large to small scales we will prove that the increments have the desired properties.

\subsection{Square averages}

In \cite{BW20,W22} the Ginzburg-Landau field was analyzed via a circle average process. We use the same strategy, but for technical reasons (namely the fact that Theorem \ref{t:qclt} is stated only for squares) it is advantageous for us to work with square averages instead.

We begin by defining the constant $\g$ that appears in Theorem \ref{t:main}. In view of the quantitative CLT in Theorem \ref{t:qclt}, $\g$ should be proportional to ${1}/{\ahom}$. We define the proportionality constant in terms of (continuous) square averages. Namely, for $r>0$ let $\bar{Q}_r=(-r,r)^2$ be the continuous square of sidelength $2r$ centered at 0. Consider the continuous harmonic measure $\bar{a}_{\bar{Q}_r}(x,\cdot)$ on $\bar{Q}_r$ (that is, for a Borel-set $A\subset\partial\bar{Q}_r$, $\bar{a}_{\bar{Q}_r}(x,A)$ is the probability that the first point where Brownian motion started from $x\in \bar{Q}_r$ hits $\partial\bar{Q}_r$ is in in $A$.). Then define the constant
\begin{equation}\label{e:def_c_*}
c_*:=\int_{[-1,1]^2\times[-1,1]^2} G^{-\Delta}(x,y)d\left(\bar{a}_{\bar{Q}_{1/\e}}-\bar{a}_{\bar{Q}_1}\right)(0,x)d\left(\bar{a}_{\bar{Q}_{1/\e}}-\bar{a}_{\bar{Q}_1}\right)(0,y).
\end{equation}
One can check that $0<c_*<\infty$, and we define
\begin{equation}\label{e:def_g}
\g:=\frac{c_*}{\ahom}.
\end{equation}

\begin{remark}
The precise value of $c_*$ is not-important for us, and so we have chosen the definition \eqref{e:def_c_*} so that the variance estimates for the square averages (in particular \eqref{e:sharpmomentbound}) are as easy as possible.
However, let us point out that one can compute $c_*$ explicitly, and in fact $c_*={1}/{(2\pi)}$. The easiest way to see this is to consider the special case of the discrete Gaussian free field (i.e. the case $V(x)={x^2}/{2}$): In that case $\ahom=1$, while it is well-known that
\[\Var_{\mathrm{DGFF}}^{Q_N}(\phi(x))=\frac{1}{2\pi}\log\dist(x,\partial^+Q_N)+O(1),\]
which together with \eqref{e:sharpmomentbound} below means that $\g=1/(2\pi)$, and hence necessarily $c_*=1/(2\pi)$.
\end{remark}

Having defined $\g$, we can now introduce the discrete square average process that will be crucial for our argument.
For $D\subset\Z^2$ and $x\in D$ let $a_D(x,\cdot)$ be the (discrete) harmonic measure of $x$. Formally, let $S^x$ denote simple random walk starting at $x$, and $\tau _{\partial^+ D}=\inf \left\{
t>0\colon S^x_t \in \partial^+ D\right\} $ the hitting time of $\partial^+ D$. Define
\begin{equation}
\label{eq-a}
a_{D}(x,y)=\PP\left( S^x_{\tau_{\partial^+ D}} =y\right) .
\end{equation}
We also define the weighted average
\begin{equation}
  \label{eq-harmonic}
  \Gamma_D(x,\phi)=\sum_{y\in \partial^+D}a_D(x,y)\phi(y).
\end{equation}

In the following we ignore all rounding effects.

For our purposes Theorems \ref{t:miller_coupling}, \ref{t:miller_mean} and \ref{t:miller_coupling_annulus} hold for every $\Lambda>1$. So for concreteness we fix $\Lambda=2$, and drop it from our notation from now on. With this in mind, let
\begin{equation}\label{e:def_omega}
\omega=\frac{\min(1,\alpha,\beta,\gamma)}{8},
\end{equation}

 with the exponents $\alpha, \beta$ from Theorem \ref{t:miller_coupling_annulus} and $\gamma$ from Theorem \ref{t:qclt}.

We define the following sequences of radii $\left\{ r_k\right\}
_{k=1}^{\infty }$, $\left\{ r_{k,+}\right\} _{k=0}^{\infty }$ and $\left\{r_{k,-}\right\} _{k=0}^{\infty }$ by
\begin{equation}
  \label{eq-rk}
  r_{k} =N\e^{-k},\quad 
r_{k,+} =r_k+r_k^{1-\omega},\quad
r_{k,-} =r_k-r_k^{1-\omega}.
\end{equation}

We want to construct smoothed out square averages associated with the radii $r_{k,+}$ and $r_{k,-}$. To do so, fix some smooth function $\eta\in C^\infty(\R)$ that is non-negative, supported in $(-1,1)$ and has integral 1, and let $\eta_\rho(r)=\frac{1}{\rho}\eta\left(\frac{r}{\rho}\right)$.
Then $\eta_\rho$ still has integral 1, and we have
\begin{equation}\label{e:propertyeta}
\sum_{r\in\N}\eta_\rho(r+a)=1+O\left(\frac{1}{\rho}\right)
\end{equation}
for any $a\in\R$.

We then define
\begin{align}
\Sav_{k,+}(x,\phi)&=\frac{\sum_{r\in\N}\eta_{r_k^{1-\omega}}(r-r_{k,+})\Gamma_{Q_r(x)}(x,\phi)}{\sum_{r\in\N}\eta_{r_k^{1-\omega}}(r-r_{k,+})}\label{e:def_squareav+},\\
\Sav_{k,-}(x,\phi)&=\frac{\sum_{r\in\N}\eta_{r_k^{1-\omega}}(r-r_{k,-})\Gamma_{Q_r(x)}(x,\phi)}{\sum_{r\in\N}\eta_{r_k^{1-\omega}}(r-r_{k,-})}\label{e:def_squareav-}.
\end{align}
Thus \eqref{e:def_squareav+} is a radially smoothed-out version of a square average of $\phi$, where the radial smoothing takes places over a strip of width $r_k^{1-\omega}$. The denominator ensures that if $\phi$ is constant then $\Sav_{k,+}(x,\phi)$ is equal to that same constant. Moreover the denominator is very close to 1, as in view of \eqref{e:propertyeta} we have
\begin{equation}\label{e:num_close_to_one}
\sum_{r\in\N}\eta_{r_k^{1-\omega}}(r-r_{k,+})=1+O\left(\frac{1}{r_k^{1-\omega}}\right).
\end{equation}
Analogous statements hold true for \eqref{e:def_squareav-}.

Moreover, note that if $k$ is so small that $r_k>\dist(x,\partial^+Q_N)$ then $S_{k,\pm}$ are not well-defined. This is not an issue for us as we will use $S_{k,\pm}$ only for large enough indices $k$.

Using $\Sav_{k,+}$ and $\Sav_{k,-}$ we define the increment process
\[\Inc_k(x,\phi)=\Sav_{k+1,+}(x,\phi)-\Sav_{k,-}(x,\phi),\]
and the boundary layer process
\[\Bla_k(x,\phi)=\Sav_{k,-}(x,\phi)-\Sav_{k,+}(x,\phi).\]

With these definitions, we can write a telescoping sum
\begin{equation}\label{e:telescope}
\phi(x)=\sum_{k=k_0}^{k_\infty-1}\Inc_k(x,\phi)+\sum_{k=k_0+1}^{k_\infty-1}\Bla_k(x,\phi)+\left(\phi(x)-\Sav_{k_\infty,+}(x,\phi)\right)+\Sav_{k_0,-}(x,\phi)
\end{equation}
for some parameters $\log\frac{N}{\dist(x,\partial^+Q_N)}\le k_0\le k_\infty\le\log N$ to be chosen later. Note that both $k_0$ and $k_\infty$ may depend on $N$ and on $x$.

The summands of the first sum will be the main terms in our approximation, and we will argue using Theorem \ref{t:miller_coupling_annulus} and Theorem \ref{t:qclt} that they are approximately Gaussian and approximately independent. The summands of the second sum are typically very small so that we will be able to neglect them. Finally, the third and fourth summand will be 
estimated
using the Brascamp-Lieb inequality from Lemma \ref{l:BL}.

We define in this context some events that characterize a typical behaviour of the field. We will think of these events as good events and show that they occur with large enough probability.

We begin with the event that the field is not too rough at some scale (in the sense that its oscillation is controlled). The definition is made in such a way that we can apply Theorem \ref{t:miller_coupling_annulus} whenever the field is non-rough. We define
\[\Ero_{k,x}=\left\{\phi\colon \osc_{y\in\partial^+Q_{r_k}(x)}\phi(y)\le 4(\log r_k)^2\right\}.\]
Next we consider the event that a boundary layer term is not much larger than typical. Note that
\[\Var \Bla_k(x,\phi)\approx\g\log\frac{r_{k,+}}{r_{k,-}}\approx\frac{2\g}{r_k^{\omega}},\]
 and so typically $\Bla_k(x,\phi)$ is of order $\frac{1}{r_k^{\omega/2}}$.
In view of this we define
\[\Ebd_{k,x}=\left\{\phi\colon \left|\Bla_k(x,\phi)\right|\le\frac{1}{r_k^{\omega/4}}\right\}.\]

We need to quantify that these events are very likely, even conditionally  on the field outside of the relevant square.
\begin{lemma}\label{l:badevents}
If $x\in Q_N$ and $\log\frac{N}{\dist(x,\partial^+Q_N)}\le k\le\log N$, then
\begin{align}
\PP^{Q_N,0}\left((\Ero_{k,x})^\com\right)&\le C\exp\left(-c(\log r_k)^{3}\right),\label{e:badevents1}\\
\PP^{Q_N,0}\left((\Ebd_{k,x})^\com\right)&\le C\exp\left(-cr_k^{\omega/2}\right).\label{e:badevents2}
\end{align}
Moreover, there is a constant $C'$ such that if $\log\frac{N}{\dist(x,\partial^+Q_N)}+1\le k\le\log N-C'$, then
\begin{align}
\PP^{Q_N,0}\left((\Ero_{k,x})^\com\mid\F_{Q_{r_{k-1}}^\com(x)}\right)\I_{\Ero_{k-1,x}}&\le C\exp\left(-c(\log r_k)^{3}\right),\label{e:badevents3}\\
\PP^{Q_N,0}\left((\Ebd_{k,x})^\com\mid\F_{Q_{r_{k-1}}^\com(x)}\right)\I_{\Ero_{k-1,x}}&\le C\exp\left(-cr_k^{\omega/2}\right).\label{e:badevents4}
\end{align}
\end{lemma}
\begin{proof}
We have that
\begin{align*}
\PP^{Q_N,0}\left((\Ero_{k,x})^\com\right)&\le \sum_{y,y'\in\partial^+Q_{r_k}(x)}\PP^{Q_N,0}\left((\phi(y)-\phi(y')>2(\log r_k)^2\right)\\
&\le Cr_k^2\max_{y,y'\in\partial^+Q_{r_k}(x)}\PP^{Q_N,0}\left((\phi(y)-\phi(y')>2(\log r_k)^2\right).
\end{align*}
Using the Brascamp-Lieb inequality just like in the proof of Lemma \ref{l:BLtailbound}, we find that
\begin{align*}
\PP^{Q_N,0}\left((\phi(y)-\phi(y')>4(\log r_k)^2\right)&\le C\exp\left(-c\frac{(\log r_k)^{4}}{\log r_k}\right) \le C\exp\left(-c(\log r_k)^{3}\right).
\end{align*}
Together with the previous estimate this implies \eqref{e:badevents1}. The estimate \eqref{e:badevents2} follows from an analogous application of the Brascamp-Lieb inequality.

The arguments for \eqref{e:badevents3} and \eqref{e:badevents4} are more involved, as this time we cannot just use the Brascamp-Lieb inequality, but also need to control the conditional expectation of the quantity in question. For that purpose we need to use that on $\Ero_{k-1,x}$ this conditional expectation is well-behaved.

We begin with \eqref{e:badevents3}. The Helffer-Sj\"ostrand representation implies that $\E^{Q_N,0}\left(\phi(\cdot)\mid\F_{Q_{r_{k-1}}^\com(x)}\right)$ assumes its extrema over $Q_{r_{k-1}}(x)$ on the boundary. On the event $\Ero_{k-1,x}$ we know that \[\osc_{y\in\partial^+Q_{r_{k-1}}(x)}\phi(y)\le 4(\log r_{k-1})^2,\] 
and hence also
\begin{equation}\label{e:badevents5}
\osc_{y\in Q_{r_{k-1}}(x)}\E^{Q_N,0}\left(\phi(y)\mid\F_{Q_{r_{k-1}}^\com(x)}\right)\le 4(\log r_{k-1})^2.
\end{equation}
Next, recall that by Theorem \ref{t:miller_mean}
\[\E^{Q_N,0}\left(\phi(y)\mid\F_{Q_{r_{k-1}}^\com(x)}\right)=\hat\phi(y)+O\left(\frac{1}{r_{k-1}^\omega}\right)\text{ on }Q_{r_{k-1,-}}(x),\]
where $\hat\phi$ is the harmonic function that agrees with $\E^{Q_N,0}\left(\phi(\cdot)\mid\F_{Q_{r_{k-1}}^\com(x)}\right)$ on $\partial^+Q_{r_{k-1,-}}(x)$. From \eqref{e:badevents5} we conclude
\[\osc_{y\in Q_{r_{k-1,-}}(x)}\hat\phi\le 4(\log r_{k-1})^2,\]
and the Harnack inequality for discrete harmonic functions (see e.g. \cite[Theorem 6.3.9]{LL10} implies that
\[\osc_{y\in Q_{r_{k}}(x)}\hat\phi\le 4(1-\theta)(\log r_{k-1})^2\]
for a constant $\theta>0$. We conclude that
\[
\osc_{y\in Q_{r_k}(x)}\E^{Q_N,0}\left(\phi(y)\mid\F_{Q_{r_{k-1}}^\com(x)}\right)\le 4(1-\theta)(\log r_{k-1})^2+\frac{C}{r_{k-1}^\omega},\]
and note that if $r_k$ is large enough (i.e. if we choose $C'$ large enough) then the right-hand side is bounded by $4\left(1-\frac{\theta}{2}\right)(\log r_k)^2$.
So if $(\Ero_{k,x})^\com$ occurs then we must have
\[\osc_{y\in Q_{r_k}(x)}\left(\phi(y)-\E^{Q_N,0}\left(\phi(y)\mid\F_{Q_{r_{k-1}}^\com(x)}\right)\right)\ge 2\theta(\log r_k)^2.\]
However the probability of the latter event can be estimates using Brascamp-Lieb and a union bound in just the same way as \eqref{e:badevents1}.

Finally, we turn to \eqref{e:badevents4}. For this we need to control the conditional expectation of $\Bla_{k,x}$. For that purpose we use once again Theorem \ref{t:miller_mean}. First of all, $\Bla_k(x,\cdot)$ is a linear functional and hence
\[\E\left(\Bla_k(x,\phi)\mid\F_{Q_{r_{k-1}}^\com(x)}\right)=\Bla_k\left(x,\E\left(\phi\mid\F_{Q_{r_{k-1}}^\com(x)}\right)\right).\]
Now
Theorem \ref{t:miller_mean} implies that on the event $\Ero_{k,x}$, $\E\left(\phi(\cdot)\mid\F_{Q_{r_{k-1}}^\com(x)}\right)$ is equal to a harmonic function (up to an error uniformly bounded by $\frac{C}{r_{k-1}^{\omega}}$), while $\Bla_k(x,\cdot)$ vanishes for harmonic functions. Hence
\[\left|\E^{Q_N,0}\left(\Bla_k(x,\phi)\mid\F_{Q_{r_{k-1}}^\com(x)}\right)\right|\I_{\Ero_{k-1,x}}\le \frac{C}{r_k^\omega}.\]
This means that if $(\Ebd_{k,x})^\com$ occurs then we must have
\[\Bla_k(x,\phi)-\E^{Q_N,0}\left(\Bla_k(x,\phi)\mid\F_{Q_{r_{k-1}}^\com(x)}\right)\ge\frac{1}{r_k^{\omega/4}}-\frac{C}{r_k^\omega}\ge \frac{1}{2r_k^{\omega/4}}\]
(where the last inequality holds as soon as $C'$ and hence $r_k$ are large enough). The probability of that event can again be estimated using Brascamp-Lieb (as for \eqref{e:badevents2}).
\end{proof}

\subsection{A single coupling step}
We will use Theorem \ref{t:miller_coupling_annulus} iteratively to couple the random variables $\Inc_k(x,\phi)$ to independent random variables. For simplicity we state first the properties of a single such coupling step. 
\begin{lemma}\label{l:simple_coupling}
Let $x\in Q_N$, and $\log\frac{N}{\dist(x,\partial^+Q_N)}\le k\le\log N-C$. 
Then one can construct on a (possibly extended) probability
space with probability measure $\tilde\PP^{Q_N,0}$ a random
variable $\Ind_k(x)$ and an event $\Ecp_{k,x}$, with the following properties:
\begin{itemize}
\item $\Ind_k(x)$ is independent of $\F_{Q_{r_k}^\com(x)}$, and equal in law to $\Inc_k(x,\phi^0)$, where $\phi^0$ is distributed according to $\PP^{Q_{r_k}(x),0}$.
\item The event $\Ecp_{k,x}$ is $\sigma\left(\F_{Q_{r_{k+1}}^\com(x)},\Ind_k(x)\right)$-measurable.
\item On the event $\Ero_{k,x}\cap \Ecp_{k,x}$ we have $\Inc_k(x,\phi)=\Ind_k(x)$.
\item We have
\[\tilde\PP^{Q_N,0}_f\left((\Ecp_{k,x})^\com\mid\F_{Q_{r_k}^\com(x)}\right)\le \frac{C}{r_k^{\beta}}\]
almost surely on the event $\Ero_{k,x}$.
\end{itemize}
\end{lemma}
\begin{proof}
  Let $\phi^{Q_N,0}\restriction_{\partial^+Q_{r_k}(x)}=f$, and note that the conditional law of $\phi^{Q_N,0}$ given $\phi^{Q_N,0}\restriction_{Q_{r_k}^\com(x)}$  is equal to $\phi^{Q_{r_k(x)},f}$.
 We have that $\phi^{Q_N,0}\in\Ero_{k,x}$ if and only if $\osc_{x\in \partial^+ Q_{r_k(x)}}|f|\le 4|\log R|^2$.

We claim that,  if $\phi^{Q_N,0}\in\Ero_{k,x}$,  then we may apply Theorem \ref{t:miller_coupling_annulus} to couple $\phi^{Q_{r_k}(x),f}$ and $\phi^{Q_{r_k}(x),0}$ on $Q_{r_k}(x)\setminus Q_{r_{k+1}}(x)$.  Indeed,  applying Theorem \ref{t:miller_coupling_annulus} on $Q_{r_k}(x)\setminus Q_{r_{k+1}}(x)$ with $r:=r_k^{1-\omega}$  yields a specific family of probability measures $\PP^f$ on $(\R^{Q_{r_k}(x)})^2$ such that the first marginal is the law of $\phi^{Q_{r_k}(x),f}$ and the second marginal is the law of $\phi^{Q_{r_k}(x),0}$. In particular, the law of the second marginal does not depend on $f$. 
This means that the Markov kernel $f\mapsto\PP^f(\cdot)$ actually induces a coupling of the family $\Law\left(\phi^{Q_{r_k}(x),f}\restriction_{Q_{r_k}(x)\setminus Q_{r_{k+1}}(x)}\right)$ to the fixed $\Law\left(\phi^{Q_{r_k}(x),0}\restriction_{Q_{r_k}(x)\setminus Q_{r_{k+1}}(x)}\right)$. In particular 
we can fix, once and for all,
a realization of $\phi^0$ distributed according to $\mu^{Q_{r_k}(x),0}$ (which is independent of $\F_{Q_N\setminus Q_{r_k}(x)}$),
and construct for each $f$, the random variable  $\phi^{D,f}$ distributed 
according to $\phi^{Q_{r_k(x)},f}$, in a measurable way in $f$. By the Markov property of the field, $\phi^{Q_N,0}$ can be obtained by first sampling $\phi^{Q_N,0}\restriction_{Q_{r_k}^\com(x)}$, 
  and then resampling $\phi^{Q_N,0}\restriction_{Q_{r_k}(x)}$ according to $\PP^{Q_{r_k},f}$. So we can construct the coupling by starting with the realization $\phi^0$, then sampling independent of it $\phi^{Q_N,0}\restriction_{Q_{r_k}^\com(x)}$, and finally (re-)sampling $\phi^{Q_N,0}\restriction_{Q_{r_k}(x)}$ using the coupling induced by $\PP^f$.
  
Moreover, according to Theorem \ref{t:miller_coupling_annulus}, the coupling has the property that there is an event $\Ecp_{k,x}:=\Ev_k$ on which we have
\begin{equation}\label{e:simple_coupling1}
\phi^{Q_{r_k},f}=\phi^0+\hat{\phi}
\end{equation}
in $Q_{r_k-r_k^{1-\omega}}(x)\setminus Q_{r_{k+1}(x)}$.
Now note that the linear functional $\phi\mapsto\Inc_k(x,\phi)=\Sav_{k+1,+}(x,\phi)-\Sav_{k,-}(x,\phi)$ can be represented as $\langle\phi,g\rangle$ for some function $g$, and inspecting the definitions reveals that this $g$ is supported in
\[Q_{r_{k,-}+r_k^{1-\omega}}\setminus Q_{r_{k+1,+}-r_k^{1-\omega}}\subset Q_{r_k-r_k^{1-\omega}}(x)\setminus Q_{r_{k+1}(x)}.\]
That is, the support of $g$ is contained in the set of sites on which \eqref{e:simple_coupling1} holds. Thus, on the event $\Ero_{k,x}\cap \Ecp_{k,x}$ we have
\[\Inc_k(x,\phi)=\Inc_k(x,\phi^0)+\Inc_k(x,\hat\phi).\]
Now our definition of square averages was made in such a way that for each harmonic function $h$ we have
\[h(x)=\Sav_{k+1,+}(x,h)=\Sav_{k,-}(x,h).\]
This means that $\Inc_k(x,\hat\phi)$ vanishes, and so we actually have $\Inc_k(x,\phi)=\Inc_k(x,\phi^0)$ on the event $\Ero_{k,x}\cap \Ecp_{k,x}$. So if we define $\Ind_k(x)=\Inc_k(x,\phi^0)$ then $\Ind_k(x)$ has all the properties claimed in the lemma.
\end{proof}

\begin{remark}
  \label{rem-2011}
  In the proof we have used Theorem 
  \ref{t:miller_coupling_annulus} which provides one coupling for each fixed boundary value $f$. The resampling argument ensured that nonetheless we obtain one coupling between $\phi^{Q_N,0}$ and $\phi^0$ (and not a family of couplings indexed by $f$).
  
An alternative approach would be to directly couple all the $\phi^{Q_{r_k},f}$ \emph{simultaneously}. Such a coupling can be obtained, using
  the pairwise, multi-marginal coupling provided by 
  \cite[Theorem 3]{AL19}, which allows one to not only construct a random variable $\phi^0$ as claimed in the proof but in fact construct at once 
  the family 
  $(\phi^0,\phi^f)$ (first with $f$ in some countable set, and then for all $f$ by continuity), 
 with only a multiplicative error in the pairwise 
  variation distance.
\end{remark}

\subsection{Iterated coupling}

As already mentioned, we can use the coupling result from Lemma \ref{l:simple_coupling} iteratively to couple the main term in \eqref{e:telescope} to the sum of independent random variables. We state this rather technical result in detail.

\begin{lemma}\label{l:iterated_coupling}
There is a constant $C$ and an extension $\tilde\PP^{Q_N,0}$ of $\PP^{Q_N,0}$,  such that the random variable $\phi$ is distributed according to $\PP^{Q_N,0}$,  and for any $x\in Q_N$, $\log\frac{N}{\dist(x,\partial^+Q_N)}\le k\le\log N-C$, there exist random variables $\Ind_k(x)$, $\Err_k(x)$ and events $\Ecp_{k,x}$
with the following properties:
\begin{itemize}
\item For any $x\in Q_N$ and any $\log\frac{N}{\dist(x,\partial^+Q_N)}\le k\le\log N-C$ we have
\begin{equation}\label{e:iterated_coupling1}
\I_{\Ero_{k,x}}\Inc_k(x,\phi)=\I_{\Ero_{k,x}}\left(\Ind_k(x)+\I_{(\Ecp_{k,x})^\com}\Err_k(x)\right).
\end{equation}
Moreover, for any $\log\frac{N}{\dist(x,\partial^+Q_N)}\le k_0\le k_\infty\le\log N-C$ we have
\begin{equation}\label{e:telescope2}
\begin{split}
\phi(x)&=\sum_{k=k_0}^{k_\infty-1}\I_{\Ero_{k,x}}\left(\Ind_k(x)+\I_{(\Ecp_{k,x})^\com}\Err_k(x)\right)\\
&\qquad+\sum_{k=k_0}^{k_\infty-1}\I_{(\Ero_{k,x})^\com}\Inc_k(x,\phi)+\sum_{k=k_0+1}^{k_\infty-1}\Bla_k(x,\phi)+\left(\phi(x)-\Sav_{k_\infty,+}(x,\phi)\right)+\Sav_{k_0,-}(x,\phi).
\end{split}
\end{equation}
\item For any $k$ the random variable $\Ind_k(x)$ is equal in law to $\Inc_k(x,\phi^0)$, where $\phi^0$ is distributed according to $\PP^{Q_{r_k}(x),0}$. Moreover, it is independent of
\[\tilde \F_{Q_{r_k}^\com(x)}:=\sigma\left(\F_{Q_{r_k}^\com(x)},\{\Ind_j(y)\colon Q_{r_j}(y)\cap Q_{r_k}(x)=\varnothing\}\right).\]
\item The event $\Ecp_{k,x}$ is $\tilde\F_{Q_{r_{k+1}}^\com(x)}$-measurable.
\item We have
\[\tilde\PP^{Q_N,0}\left((\Ecp_{k,x})^\com\mid\tilde\F_{Q_{r_k}^\com(x)}\right)\le \frac{C}{R^{\beta}}.\]
almost surely on the event $\Ero_{k,x}$.
\item Finally, Lemma \ref{l:badevents} continues to apply even when conditioning on $\tilde\F$ instead of $\F$. That is, we have
\begin{align}
\tilde\PP^{Q_N,0}\left((\Ero_{k,x})^\com\mid\tilde\F_{Q_{r_{k-1}}^\com(x)}\right)\I_{\Ero_{k-1,x}}&\le C\exp\left(-c(\log r_k)^{3}\right),\label{e:badevents3tilde}\\
\tilde\PP^{Q_N,0}\left((\Ebd_{k,x})^\com\mid\tilde\F_{Q_{r_{k-1}}^\com(x)}\right)\I_{\Ero_{k-1,x}}&\le C\exp\left(-cr_k^{\omega/2}\right).\label{e:badevents4tilde}
\end{align}
\end{itemize}
\end{lemma}
\begin{proof}
Let $\prec$ denote the lexicographic ordering on $\Z^2$. 
We will construct the coupling inductively, and order the $\Ind_k(x)$ first according to increasing $k$ and then among equal $k$ according to $\prec$.

We will construct the $\Ind_k(x)$ so that they are independent of
\[\tilde{\tilde \F}_{Q_{r_k}^\com(x)}:=\sigma\left(\F_{Q_{r_k}^\com(x)},\{\Ind_j(y)\colon j<k\text{ or } j=k\text{ and }y\prec x,Q_{r_j}(y)\cap Q_{r_k}(x)=\varnothing\}\right).\]
This is a smaller $\sigma$-algebra than $\tilde \F_{Q_{r_k}^\com(x)}$, but the full statement then follows from the easy fact that if in a finite sequence of random variables each is independent of the previous random variables, then they are in fact all independent.

So take some $k,x$ and suppose that the $\Ind_j(y)$ with $j<k$ or $j=k$ and $y\prec x$ have already been constructed. We want to use Lemma \ref{l:simple_coupling} to obtain on a (possible enlarged) probability space a new random variable $\Ind_k(x)$ and an event $\Ecp_{k,x}$ that the coupling succeeds. There is a small subtlely, as the lemma only gives that $\Ind_k(x)$ and $\Ecp_{k,x}$ are independent of $\F_{Q_{r_k}^\com(x)}$, not $\tilde{\tilde\F}_{Q_{r_k}^\com(x)}$. However, an inspection of the proof of the lemma shows that this can be rectified easily: First we can add $\phi^0$ independently of $\tilde{\tilde\F}_{Q_{r_k}^\com(x)}$, then resample $\phi^{Q_N,0}\restriction_{Q_{r_k}^\com(x)}$ using the coupling, and finally resampling those $X_j(y)$ for which $\{\Ind_j(y)\colon Q_{r_j}(y)\cap Q_{r_k}(x)\neq\varnothing\}$. In this manner we can couple the new random variable $\Ind_k(x)$ to the field and the $X_j(y)$ already constructed in such a manner that it is independent of $\tilde{\tilde\F}_{Q_{r_k}^\com(x)}$, as claimed. The same holds also for $\Ecp_{k,x}$, and so this completes the construction of $\Ind_k(x)$.

Having constructed the $\Ind_k(x)$ with the right independence properties, the other properties then follow from Lemma \ref{l:simple_coupling}.
\end{proof}

From now on we will work with the fixed probability space $\tilde\PP^{Q_N,0}$. We denote expectation and variance with respect to it by $\tilde\E^{Q_N,0}$, $\tilde\Var^{Q_N,0}$.  We also define,  for $x,y\in\Z^2$,   $|x-y|_{\infty} := \max\{ |x_1-y_1|,  |x_2-y_2|\} $.

Of course we need a finer understanding of $\Ind_k$ and $\Err_k$ in order  to use \eqref{e:telescope2} at its full strength. We will describe the properties of these random variables by means of estimates on their moment generating function.

\begin{lemma}\label{l:exponential_moment_bounds}
Consider the setting of Lemma \ref{l:iterated_coupling}. Let $x\in Q_N$, let $\log\frac{N}{\dist(x,\partial^+Q_N)}\le k_0\le k\le k_\infty\le \log N-C$
 and let $\lambda_*>0$ be arbitrary. Then we have the following estimates on the exponential moments.
\begin{itemize}
\item
The $\Ind_k(x)$ have approximately Gaussian distribution in the sense that for any $|\lambda|\le\lambda_*$ we have
\begin{equation}\label{e:exponential_moment_bounds1}
\left|\log\tilde\E^{Q_N,0}\exp\left(\lambda\Ind_k(x)\right)-\frac{\lambda^2\g}{2}\right|\le\frac{C_{\lambda_*}}{r_k^{\omega/2}}.
\end{equation}
where $\g$ is the constant defined in \eqref{e:def_g}.
\item For any $|\lambda|\le\lambda_*$ we have
\begin{equation}\label{e:exponential_moment_bounds2}\left|\log\tilde\E^{Q_N,0}\left(\exp\left(\lambda\I_{(\Ecp_{k,x})^\com}\Err_k(x)\right)\middle|\tilde\F_{Q_{r_k}^\com(x)}\right)-1\right|\I_{{\Ero_{k,x}}}\le \frac{C_{\lambda_*}}{r_k^{\omega}}.
\end{equation}
\item For any $|\lambda|\le\lambda_*$ we have
\begin{equation}\label{e:exponential_moment_bounds3}
\left|\log\tilde\E^{Q_N,0}\left(\exp\left(\lambda(\phi(x)-\Sav_{k_\infty,-}(x,\phi))\right)\middle|\tilde\F_{Q_{r_{k_\infty}}^\com(x)}\right)\right|\I_{{\Ero_{k_\infty,x}}}\le C_{\lambda_*}(\log N-k_\infty).
\end{equation}
\item For any $|\lambda|\le\lambda_*$ we have
\begin{equation}\label{e:exponential_moment_bounds4}\log\tilde\E^{Q_N,0}\exp\left(\lambda\left(\Sav_{k_0,+}(x,\phi)\right)\right)\le C_{\lambda_*}\left(k_0-\log\frac{N}{\dist(x,\partial^+Q_N)}\right).
\end{equation}
\end{itemize}
\end{lemma}

\begin{proof}
We begin with the proof of \eqref{e:exponential_moment_bounds1}. According to Lemma \ref{l:iterated_coupling}, $\Ind_k(x)$ is equal in law to $\Inc_k(x,\phi^0)$, where $\phi^0$ is distributed according to $\PP^{Q_{r_k}(x),0}$.   
Notice that we may write
$$
\Inc_k(x,\phi^0) = S_{k+1,+}(x,\phi^0) -  S_{k,-}(x,\phi^0) 
= \frac 1{r_k}  \sum_{y\in Q_{r_k}(x)} \phi^0(y) g(y),
$$
where for $r= |y-x|_\infty$,
\begin{equation}
\label{e.g}
\frac 1{r_k}g(y) = \frac{\eta_{r_{k+1}^{1-\omega}}(|y-x|_\infty-r_{k+1,+})a_{Q_r(x)}(x,y)}{\sum_{\tilde r\in\N}\eta_{r_{k+1}^{1-\omega}}(\tilde r-r_{k+1,+})}
-
\frac{\eta_{r_k^{1-\omega}}(|y-x|_\infty-r_{k,-})a_{Q_r(x)}(x,y)}{\sum_{\tilde r\in\N}\eta_{r_k^{1-\omega}}(\tilde r-r_{k,-})}.
\end{equation}
If we consider $g(r_k(\cdot-x))$ as a function on $\frac{1}{r_k}\Z^2\cap[-1,1]^2$, we can extend it to a function $\hat g$ on $[-1,1]^2$, for example by subdiving each face of $\frac{1}{r_k}\Z^2$ along a diagonal and defining $\hat g$ so that it is piecewise affine on each triangle.

Recall the continuous harmonic measure introduced at the beginning of Section \ref{sec-RWapp}. It is clear from our construction that $\hat g\to \bar{a}_{\bar{Q}_{1/\e}}(0,\cdot)-\bar{a}_{\bar{Q}_1}(0,\cdot)$ weakly in the space of (signed) measures. To quantify the speed of convergenge, note that standard estimates for harmonic measure (e.g. \cite[Proposition 8.1.4]{LL10}) imply that $ra_{Q_r(x)}(x,r(\cdot+y))$ agrees with the density of $\bar{a}_{\bar{Q}_1}(0,\cdot)$ up to an error of order $\frac{1}{r^2}$ (which is negligible for our purposes). Because of the two cut-off functions, we only need to focus on radii $r$ with $|r-r_{k+1}|\le r_{k+1}^{1-\omega}$ or $|r-r_{k}|\le r_{k}^{1-\omega}$. From this one can show the quantitative estimate
\begin{equation}\label{e:conv_H-1}
\left\|\hat g-\bar{a}_{Q_*}(0,\cdot)\right\|_{H^{-1}((-1,1)^2)}\le \frac{C}{r_k^{\omega/2}}
\end{equation}
Moreover, we know that $a$ is Lipschitz-continuous in the sense that for $y,y'\in Q_{r_k}(x)$,  $y\sim y'$, if we  denote by $r= |y-x|_\infty$ and $r'= |y'-x|_\infty$, we have
\[
|a_{Q_r(x)}(x,y) - a_{Q_{r'}(x)}(x,y')|\le
\frac C {r}.
\]
(this follows for example from the explicit formula in \cite[Proposition 8.13]{LL10}). From this we can easily estimate that
\begin{equation}\label{e:hatg_holder}
\left\|\hat g\right\|_{C^{0,1}((-1,1)^2)}\le Cr_k^{2\omega}
\end{equation}

We  want to apply Theorem \ref{t:qclt} to $g,\hat g$. The only issue with that it that in general $\hat g$ is only Lipschitz-continuous, and not $C^{1,\kappa}$ for any $\kappa>0$. However,  applying Lemma \ref{l.approx} below for some $\kappa>0$ to be chosen shortly, we see that there exists $\tilde g\colon[-1,1]^2\to\R$ such that $\tilde g = \bar g$ on $\frac{1}{r_k} \Z^2 \cap [-1,1]^2$. Moreover, using \eqref{e:hatg_holder} we have that
\begin{equation}\label{e:bound_H-1}
\left\|\hat g-\tilde g\right\|_{H^{-1}((-1,1)^2)}\le\left\|\hat g-\tilde g\right\|_{L^\infty((-1,1)^2)}\le \frac{C}{r_k}\left\|\hat g\right\|_{C^{0,1}((-1,1)^2)}\le \frac{C}{r_k^{1-2\omega}}
\end{equation}
and 
\begin{equation}\label{e:tildeg_holder}\left\|\tilde g\right\|_{C^{1,\kappa}((-1,1)^2)}\le Cr_k^\kappa\left\|\bar g\right\|_{C^{0,1}((-1,1)^2)}\le Cr_k^{2\omega+\kappa}.
\end{equation}

We also notice that $g$ has mean zero,  therefore $g = \nabla^*\cdot f$,  such that
\[
f(e) = \sum_{y\in Q_{r_k}(x)} g(y) \nabla G_{Q_{r_k}(x)}(e,y),
\]
where $G_{Q_{r_k}(x)}$ is the Dirichlet Green's function in $Q_{r_k}(x)$ with zero boundary condition.  
Notice that  $g(y)$ is a weighted average of $r_k a_{Q_r(x)} (x,y)$,  for $r$ such that $|r-r_{k+1,+}|\le r_{k+1}^{1-\omega}$ or $|r-r_{k,-}|\le r_{k}^{1-\omega}$.  For every fixed $r$,  we  apply the
Cauchy-Schwarz which yields
 \begin{multline*}
 \sum_{y\in \partial^+ Q_{r}(x)} r_k a_{Q_r(x)} (x,y) \nabla G_{Q_{r_k}(x)}(e,y)
 \le
 C \sum_{y\in \partial^+ Q_{r}(x)}| \nabla G_{Q_{r_k}(x)}(e,y)|\\
 \le
 C  \biggl(\sum_{y\in Q_{r_k}(x)} (\nabla G_{Q_{r_k}(x)}(e,y))^2\biggr)^{1/2} 
 =
 C G_{Q_{r_k}}(y,y)^{1/2}
 \le
 C (\log r_k)^{1/2}.
 \end{multline*}
 Averaging over $r$ we conclude that
 $\| f\|_{L^\infty} \le C (\log r_k)^{1/2}.$  This implies $\| f\|_{\underline L^2(Q_{r_k}(x))} \le C (\log r_k)^{1/2}.$
Therefore we may apply Lemma \ref{l:qclt_exponentialmoment} to $\frac 1{r_k}  \sum_{y\in Q_{r_k}(x)} \phi^0(y) \tilde g(\frac {y-x}{r_k})$ and  conclude that for $r_k$ large enough, 
\[
\left|\log\tilde\E^{Q_N,0}\exp\left(\lambda\Ind_k(x)\right)-\frac{\lambda^2\sigma^2}{2}\right|\le\frac{C_{\lambda_*}(\| \tilde g\|_{C^{1,\kappa}((-1,1)^2}^2 + \|  f\|^2_{\underline L^2(Q_{r_k}(x))})}{r_k^{\gamma}}
\le
\frac{C_{\lambda_*} }{r_k^{\gamma -4\omega-2\kappa}}.
\]
where
\[\sigma^2= \frac{1}{\ahom}\int_{(-1,1)^2\times(-1,1)^2} \tilde g(x)G^{-\Delta}(x,y) \tilde g(y)dxdy.\]

Our choice of $\omega$ in \eqref{e:def_omega} ensures that, if we take $\kappa=\omega$, then $\gamma -4\omega-2\kappa>\omega$. It remains to verify that $\sigma^2$ is close enough to $\g$. To see this, we can combine \eqref{e:conv_H-1}, \eqref{e:bound_H-1} and the definition of $c_*$ in \eqref{e:def_c_*} and compute that
\begin{align*}
\sigma^2&=\frac{1}{\ahom}\int_{(-1,1)^2\times(-1,1)^2} G^{-\Delta}(x,y)d\left(\bar{a}_{\bar{Q}_{1/\e}}-\bar{a}_{\bar{Q}_1}\right)(0,x)d\left(\bar{a}_{\bar{Q}_{1/\e}}-\bar{a}_{\bar{Q}_1}\right)(0,y)+O\left(\frac{1}{r_k^{\omega/2}}\right)\\
&=\frac{c_*}{\ahom}+O\left(\frac{1}{r_k^{\omega/2}}\right)=\g+O\left(\frac{1}{r_k^{\omega/2}}\right)
\end{align*}
This completes the proof of \eqref{e:exponential_moment_bounds1}.

For \eqref{e:exponential_moment_bounds2}, we need to use the Brascamp-Lieb inequality. We know from \eqref{e:iterated_coupling1} that
on the event ${\Ero_{k,x}}$ we have
\[\I_{(\Ecp_{k,x})^\com}\Err_k(x)=\I_{(\Ecp_{k,x})^\com}\left(\Inc_k(x,\phi)-\Ind_k(x)\right).\]
It is straightforward to bound the exponential moment of $\Inc_k(x,\phi)$ and $\Ind_k(x)$, and the difficulty is in obtaining the decaying factor $\frac{1}{r_k^{\omega}}$. For that purpose we will use H\"{o}lder's inequality. Indeed, we have
\begin{equation}\label{e:exponential_moment_bounds5}
\begin{split}
&\left|\tilde\E^{Q_N,0}\left(\exp\left(\lambda\I_{(\Ecp_{k,x})^\com}\Err_k(x)\right)\middle|\tilde\F_{Q_{r_k}^\com(x)}\right)-1\right|\\
&
=\tilde\E^{Q_N,0}\left(\I_{(\Ecp_{k,x})^\com}\exp\left(\lambda\left(\Inc_k(x,\phi)-\Ind_k(x)\right)\right)\middle|\tilde\F_{Q_{r_k}^\com(x)}\right)\\
&
\leq \! \left(\tilde\E^{Q_N,0}\left(\exp\left(4\lambda\Inc_k(x,\phi)\right)\middle|\F_{Q_{r_k}^\com(x)}\right)\right)^{1/4}\!\left(\tilde\E^{Q_N,0}\exp\left(4\lambda\Ind_k(x)\right)\right)^{1/4}\!\left(\tilde\PP^{Q_N,0}\left((\Ecp_{k,x})^\com\middle|\tilde\F_{Q_{r_k}^\com(x)}\right)\right)^{1/2}.
\end{split}
\end{equation}
The second factor can be estimated using the Brascamp-Lieb inequality, and Lemma \ref{l:simple_coupling} allows to control the third one. Regarding the first factor, one needs to be more careful, as the conditional expectation of $\Inc_k$ given $\F_{Q_{r_k}^\com(x)}$ need not be 0. However, Theorem \ref{t:miller_mean} implies that on the event $\Ero_{k,x}$ the conditional expectation of $\phi$ given $\F_{Q_{r_k}^\com(x)}$ is close to harmonic, and so
\[\left|\tilde\E^{Q_N,0}\left(\Inc_k(x,\phi)\middle|\tilde\F_{Q_{r_k}^\com(x)}\right)\right|\le\frac{C}{r_k^{\beta}}.\]
On the other hand, by Brascamp-Lieb we have
\[\tilde\E^{Q_N,0}\left(\exp\left(4\lambda\Inc_k(x,\phi)-4\lambda\tilde\E^{Q_N,0}\left(\Inc_k(x,\phi)\middle|\tilde\F_{Q_{r_k}^\com(x)}\right)\right)\middle|\tilde\F_{Q_{r_k}^\com(x)}\right)\le \exp(C\lambda^2).\]
Combining this with the previous estimate, we can bound
\[\tilde\E^{Q_N,0}\left(\exp\left(4\lambda\Inc_k(x,\phi)\right)\middle|\tilde\F_{Q_{r_k}^\com(x)}\right)\le \exp\left(C\lambda^2+\frac{C\lambda}{r_k^{\beta}}\right)\le C_{\lambda_*}.\]
This controls the first factor in \eqref{e:exponential_moment_bounds5}.

Altogether, we obtain
\begin{align*}
\left|\tilde\E^{Q_N,0}\left(\exp\left(\lambda\I_{(\Ecp_{k,x})^\com}\Err_k(x)\right)\middle|\F_{Q_{r_k}^\com(x)}\right)-1\right|&\le\left(C_{\lambda_*}\right)^{1/4}\left(\exp(4C\lambda^2\right)^{1/4}\left(\frac{C}{r_k^{\beta}}\right)^{1/2}
\le\frac{C_{\lambda_*}}{r_k^{\omega}},
\end{align*}
which immediately implies \eqref{e:exponential_moment_bounds2}.

Regarding \eqref{e:exponential_moment_bounds3}, we can argue similarly. Theorem \ref{t:miller_mean} implies that
\[\left|\tilde\E^{Q_N,0}\left(\lambda(\phi(x)-\Sav_{k_\infty,-}(x,\phi))\middle|\tilde\F_{Q_{r_{k_\infty}}^\com(x)}\right)\right|\I_{{\Ero_{k_\infty,x}}}\le\frac{C}{r_{k_\infty}^\beta},\]
while from the Brascamp-Lieb inequality we obtain
\begin{align*}
&\tilde\E^{Q_N,0}\left(\exp\left(\lambda(\phi(x)-\Sav_{k_\infty,-}(x,\phi))-\lambda\tilde\E^{Q_N,0}\left(\phi(x)-\Sav_{k_\infty,-}(x,\phi)\middle|\tilde\F_{Q_{r_{k_\infty}}^\com(x)}\right)\right)\middle|\tilde\F_{Q_{r_{k_\infty}}^\com(x)}\right)\\
&\le \exp(C\lambda^2\log r_{k_\infty}).\end{align*}
As before, combining the last two displays we find that
\[\log\tilde\E^{Q_N,0}\left(\exp\left(\lambda(\phi(x)-\Sav_{k_\infty,-}(x,\phi))\right)\middle|\tilde\F_{Q_{r_{k_\infty}}^\com(x)}\right)\I_{{\Ero_{k_\infty,x}}}\le C_{\lambda_*}\log r_{k_\infty}=C_{\lambda_*}(\log N-k_\infty),\]
which is the upper bound in \eqref{e:exponential_moment_bounds3}. The lower bound follows from the observation that
\begin{align*}
&\tilde\E^{Q_N,0}\left(\exp\left(\lambda(\phi(x)-\Sav_{k_\infty,-}(x,\phi))\right)\middle|\tilde\F_{Q_{r_{k_\infty}}^\com(x)}\right)\\
&\ge\left(\tilde\E^{Q_N,0}\left(\exp\left(-\lambda(\phi(x)-\Sav_{k_\infty,+}(x,\phi))\right)\middle|\tilde\F_{Q_{r_{k_\infty}}^\com(x)}\right)\right)^{-1},\end{align*}
and the upper bound (applied for $-\lambda$ instead of $\lambda$).

The argument for \eqref{e:exponential_moment_bounds4} follow once again from the Brascamp-Lieb inequality, we omit further details.
\end{proof}

\begin{lemma}
\label{l.approx}
Let $f\colon(-1,1)^2 \to \R$ be a Lipschitz continuous function,  $\eps>0$ be a fixed constant.  Then there exists a function $g\colon(-1,1)^2 \to \R$ such that $g(x) = f(x)$ for $x\in \eps\Z^2 \cap (-1,1)^2$,
\[\|f-g\|_{L^\infty((-1,1)^2)}\le C\eps\|f\|_{C^{0,1}((-1,1)^2)}\]
and 
\[\| g\|_{C^{1,\alpha}((-1,1)^2)} \le C\eps^{-\alpha} \|f\|_{C^{0,1}((-1,1)^2)}.\]
\end{lemma}
\begin{proof}
Let $\varphi$ be a standard nonnegative mollifier,  $\varphi_\eps (x) = \eps^{-2} \varphi(x/\eps)$, and set $\Phi_\eps(f) =  \varphi_\eps * f$. Then there exists $C>0$,  such that $\|\Phi_\eps (f)\|_{C^{1,\alpha}} \le C \eps^{-\alpha}  \|f\|_{C^{0,1}}$.  For every $x\in \eps\Z^2 \cap (-1,1)^2$, 
\[
| \Phi_\eps (f)(x) - f(x) | \le C  \|f\|_{C^{0,1}((-1,1)^2)}\int |y-x|  \varphi_\eps (x-y) \, dy
\le
C \eps \|f\|_{C^{0,1}((-1,1)^2)}.
\]
Let $h_x: [0,1]^2 \to \R$ be a symmetric smooth cutoff function such that $h_x(x) = f(x) - \Phi_\eps (f)(x)$ and $h_x (y) = 0$ for $y \in [0,1]^2\setminus B(x,\eps/2)$.  It is clear that one may choose an $h_x$ such that $\|h_x\|_{C^{0,1}((-1,1)^2)} \le C \eps^{-\alpha} \|f\|_{C^{0,1}((-1,1)^2)}$ and $\|h_x\|_{L^\infty((-1,1)^2)} \le C\eps \|f\|_{C^{0,1}((-1,1)^2)}$.  Take $H = \sum_{x\in \eps\Z^2 \cap (-1,1)^2} h_x$,  and define $g = \Phi_\eps (f) + H$.  Then $g(x) = f(x)$ for $x\in \eps\Z^2 \cap [0,1]^2$, and the two estimates follow immediately.
\end{proof}

Lemma \ref{l:exponential_moment_bounds} allows us to control the exponential moments of the individual summands in \eqref{e:telescope2}. Applying these estimates iteratively after suitable conditioning, we can obtain bounds on the exponential moments of the square averages $\Sav_{k,\pm}$.

\begin{lemma}\label{l:sharpmomentbound}
Let $N\in\N$, $x\in Q_N$.
For any $\lambda_*>0$ and any $|\lambda|\le\lambda_*$, we have
\begin{equation}
\left|\log\tilde\E^{Q_N,0}(\exp(\lambda\phi(x))-\frac{\lambda^2\g\log\dist(x,\partial Q_N)}{2}\right|\le C_{\lambda_*}.\label{e:sharpmomentbound}
\end{equation}
Moreover, if $\log\frac{N}{\dist(x,\partial^+Q_N)}\le k< k'\le\log N-C$, we also have
\begin{align}
\left|\log\tilde\E^{Q_N,0}\left(\exp\left(\lambda\Sav_{k',+}(x,\phi)\right)\right)-\frac{\lambda^2\g k'}{2}\right|&\le \frac{C_{\lambda_*}}{r_{k'}^{\omega/4}},
\label{e:sharpmomentbound_Sk}\\
\left|\log\tilde\E^{Q_N,0}\left(\exp\left(\lambda(\phi(x)-\Sav_{k,-}(x,\phi))\right)\middle|\tilde\F_{Q_{r_{k}}^\com(x)}\right)-\frac{\lambda^2\g (\log N-k)}{2}\right|\I_{{\Ero_{k,x}}}&\le C_{\lambda_*},\label{e:sharpmomentbound_phi_Sk}\\
\left|\log\tilde\E^{Q_N,0}\left(\exp\left(\lambda(\Sav_{k',+}(x,\phi)-\Sav_{k,-}(x,\phi))\right)\middle|\tilde\F_{Q_{r_{k}}^\com(x)}\right)-\frac{\lambda^2\g(k'-k)}{2}\right|\I_{{\Ero_{k,x}}}&\le \frac{C_{\lambda_*}}{r_{k'}^{\omega/4}}.\label{e:sharpmomentbound_Sk_Sk'}
\end{align}
\end{lemma}

Before we turn to the proof of these estimates, let us point out that \eqref{e:sharpmomentbound} allows us to improve upon Lemma \ref{l:BLtailbound} by replacing $c_-$ with the correct prefactor. The result obtained thereby is very similar to \cite[Theorem 1.4]{BW20}. However, with our more refined approach (and the use of a quantitative instead of a qualitative CLT) we obtain a much better bound for the error term.
\begin{lemma}\label{l:quitesharptailbound}
For any $N\in\N$, any $x\in Q_N$ and any $t\le C\log\dist(x,\partial Q_N)$ we have
\begin{equation}\label{e:quitesharptailbound}
\tilde\PP^{Q_N,0}\left( \phi(x)\ge t\right) \le C\exp\left(-\frac{1}{2\g\log\dist(x,\partial Q_N)}t^2\right).
\end{equation}
\end{lemma}
\begin{proof}
This follows from \eqref{e:sharpmomentbound},  and an application of the exponential Chebyshev inequality similar to the proof of Lemma \ref{l:BLtailbound}.
\end{proof}

\begin{proof}[Proof of Lemma \ref{l:sharpmomentbound}]
The argument for all four estimates is very similar. We give the proof for one of them, \eqref{e:sharpmomentbound_phi_Sk}, in detail, and then briefly explain how to explain the other estimates as well.

We choose $k_0=\log\frac{N}{\dist(x,\partial^+Q_N)}$ and $k_\infty=\log N-C$. Our assumption on $x$ ensures that $k_\infty\ge k\ge k_0$. To show \eqref{e:sharpmomentbound_phi_Sk}, we will establish the slightly stronger statement
\begin{equation}\label{e:sharpmomentbound1}
\begin{split}
&\left|\log\tilde\E^{Q_N,0}\left(\exp\left(\lambda(\phi(x)-\Sav_{k,-}(x,\phi))\right)\middle|\tilde\F_{Q_{r_{k}}^\com(x)}\right)-\frac{\lambda^2\g (k_\infty-k)}{2}\right|\I_{{\Ero_{k,x}}}\nonumber\\
&\qquad \le C_{\lambda_*}\left(\log N- r_{k_\infty}+\sum_{j=k+1}^{k_\infty}\frac{1}{r_j^{\omega/4}}\right),
\end{split}
\end{equation}
via reverse induction on $k$. The case $k=k_\infty$ is given by \eqref{e:exponential_moment_bounds3}.

Assume now that \eqref{e:sharpmomentbound1} holds for some $k$. We can decompose
\begin{equation}\label{e:sharpmomentbound2}
\begin{split}
&\tilde\E^{Q_N,0}\left(\exp\left(\lambda(\phi(x)-\Sav_{k-1,-}(x,\phi))\right)\middle|\tilde\F_{Q_{r_{k-1}}^\com(x)}\right)\I_{{\Ero_{k-1,x}}}\\
&=\tilde\E^{Q_N,0}\left(\I_{{\Ero_{k,x}}\cap\Ebd_{k,x}}\exp\left(\lambda(\phi(x)-\Sav_{k-1,-}(x,\phi))\right)\middle|\tilde\F_{Q_{r_{k-1}}^\com(x)}\right)\I_{{\Ero_{k-1,x}}}\\
&\qquad+\tilde\E^{Q_N,0}\left(\I_{({\Ero_{k,x}}\cap\Ebd_{k,x})^\com}\exp\left(\lambda(\phi(x)-\Sav_{k-1,-}(x,\phi))\right)\middle|\tilde\F_{Q_{r_{k-1}}^\com(x)}\right)\I_{{\Ero_{k-1,x}}},
\end{split}
\end{equation}
and we will estimate the two summands on the right-hand side separately. The second summand is an error term that can be bounded with a similar argument as in the proof of \eqref{e:exponential_moment_bounds2}. Indeed, by Cauchy-Schwarz we have
\begin{align*}
&\left|\tilde\E^{Q_N,0}\left(\I_{({\Ero_{k,x}}\cap\Ebd_{k,x})^\com}\exp\left(\lambda(\phi(x)-\Sav_{k-1,-}(x,\phi))\right)\middle|\tilde\F_{Q_{r_{k-1}}^\com(x)}\right)\right|\\
&\le\left(\tilde\E^{Q_N,0}\left(\exp\left(2\lambda(\phi(x)-\Sav_{k-1,-}(x,\phi))\right)\middle|\tilde\F_{Q_{r_{k-1}}^\com(x)}\right)\right)^{1/2}\left(\tilde\PP^{Q_N,0}\left(({\Ero_{k,x}}\cap\Ebd_{k,x})^\com\middle|\tilde\F_{Q_{r_{k-1}}^\com(x)}\right)\right)^{1/2}.
\end{align*}
On the event $\Ero_{k-1,x}$ the exponential moment can be controlled using Theorem \ref{t:miller_mean} and the Brascamp-Lieb inequality, and the probability using Lemma \ref{l:badevents}. We find that
\begin{equation}\label{e:sharpmomentbound3}
\begin{split}
&\left|\tilde\E^{Q_N,0}\left(\I_{({\Ero_{k,x}}\cap\Ebd_{k,x})^\com}\exp\left(\lambda(\phi(x)-\Sav_{k-1,-}(x,\phi))\right)\middle|\tilde\F_{Q_{r_{k-1}}^\com(x)}\right)\right|\I_{{\Ero_{k-1,x}}}\\
&\le (\exp(C\lambda^2\log r_{k-1}))^{1/2}\left(\exp(-c(\log r_{k-1})^{3/2}\right)^{1/2}\\
&\le C_{\lambda_*}\exp(-c(\log r_{k-1})^{3/2}).
\end{split}
\end{equation}

The main challenge is to bound the first summand in \eqref{e:sharpmomentbound2}. By definition we have
\[\phi(x)-\Sav_{k-1,-}(x,\phi)=\phi(x)-\Sav_{k,-}(x,\phi)+\Bla_k(x,\phi)+\Inc_k(x,\phi).\]
On the event $\Ebd_{k,x}$ we have a uniform bound on $\Bla_k(x,\phi)$, and so
\begin{align*}
&\tilde\E^{Q_N,0}\left(\I_{{\Ero_{k,x}}\cap\Ebd_{k,x}}\exp\left(\lambda(\phi(x)-\Sav_{k-1,-}(x,\phi))\right)\middle|\tilde\F_{Q_{r_{k-1}}^\com(x)}\right)\\
&=\tilde\E^{Q_N,0}\left(\I_{{\Ero_{k,x}}\cap\Ebd_{k,x}}\exp\left(
\lambda\left(\phi(x)-\Sav_{k,-}(x,\phi)+\Inc_k(x,\phi)+O\left(\frac{1}{r_k^{\omega/4}}\right)\right)\right)\middle|\tilde\F_{Q_{r_{k-1}}^\com(x)}\right)\\
&=\exp\left(O_{\lambda_*}\left(\frac{1}{4r_k^{\omega/4}}\right)\right)\tilde\E^{Q_N,0}\left(\I_{{\Ero_{k,x}}\cap\Ebd_{k,x}}\exp\left(\lambda\left(\phi(x)-\Sav_{k,-}(x,\phi)+\Inc_k(x,\phi)\right)\right)\middle|\tilde\F_{Q_{r_{k-1}}^\com(x)}\right).
\end{align*}
The increment $\Inc_k(x,\phi)$ and the indicator $\I_{\Ero_{k,x}}$ are measurable with respect to $\tilde\F_{Q_{r_{k}}^\com(x)}$, and so we can rewrite this as
\begin{equation}\label{e:sharpmomentbound4}
\begin{split}
&\tilde\E^{Q_N,0}\left(\I_{{\Ero_{k,x}}\cap\Ebd_{k,x}}\exp\left(\lambda(\phi(x)-\Sav_{k-1,-}(x,\phi))\right)\middle|\tilde\F_{Q_{r_{k-1}}^\com(x)}\right)
=\exp\left(O_{\lambda_*}\left(\frac{1}{4r_k^{\omega/4}}\right)\right)\\
&\times \tilde\E^{Q_N,0}\left(\tilde\E^{Q_N,0}\left(\I_{\Ebd_{k,x}}\exp\left(\lambda\left(\phi(x)-\Sav_{k,-}(x,\phi)\right)\right)\middle|\tilde\F_{Q_{r_{k}}^\com(x)}\right)\I_{\Ero_{k,x}}\exp\left(\lambda\Inc_k(x,\phi)\right)\middle|\tilde\F_{Q_{r_{k-1}}^\com(x)}\right).
\end{split}
\end{equation}

In order to apply our inductive assumption, we need to get rid of $\I_{\Ebd_{k,x}}$ in \eqref{e:sharpmomentbound4}. For that purpose observe that the same argument that lead to \eqref{e:sharpmomentbound3} also shows that
\[\left|\tilde\E^{Q_N,0}\left(\I_{(\Ebd_{k,x})^\com}\exp\left(\lambda\left(\phi(x)-\Sav_{k,-}(x,\phi)\right)\right)\middle|\tilde\F_{Q_{r_{k}}^\com(x)}\right)\right|\I_{\Ero_{k,x}}\le C_{\lambda_*}\exp\left(-cr_k^{\omega/2}\right).\]
Combining this with our inductive assumption
\begin{align*}
  &\tilde\E^{Q_N,0}\left(\exp\left(\lambda(\phi(x)-\Sav_{k,-}(x,\phi))\right)\middle|\tilde\F_{Q_{r_{k}}^\com(x)}\right)\I_{{\Ero_{k,x}}}\\
  &\qquad =\exp\left(\frac{\lambda^2\g (k_\infty-k)}{2}+O_{\lambda_*}\left(\log N- r_{k_\infty}+\sum_{j=k+1}^{k_\infty}\frac{1}{r_j^{\omega/4}}\right)\right),\end{align*}
we obtain that on the event $\Ero_{k,x}$ we have
\begin{align*}
&\tilde\E^{Q_N,0}\left(\I_{\Ebd_{k,x}}\exp\left(\lambda\left(\phi(x)-\Sav_{k,-}(x,\phi)\right)\right)\middle|\tilde\F_{Q_{r_{k}}^\com(x)}\right)\\
&=\exp\left(\frac{\lambda^2\g (k_\infty-k)}{2}+O_{\lambda_*}\left(\log N - r_{k_\infty}+\sum_{j=k+1}^{k_\infty}\frac{1}{r_j^{\omega/4}}\right)\right)+O_{\lambda_*}\left(\exp\left(-cr_k^{\omega/2}\right)\right)\\
&=\exp\left(\frac{\lambda^2\g (k_\infty-k)}{2}+O_{\lambda_*}\left(\log N -  r_{k_\infty}+\sum_{j=k+1}^{k_\infty}\frac{1}{r_j^{\omega/4}}+\frac{1}{4r_k^{\omega/4}}\right)\right).
\end{align*}

When we insert this into \eqref{e:sharpmomentbound4}, we see that, with
\[A:=\frac{\lambda^2\g (k_\infty-k)}{2}+O_{\lambda_*}\left(\log N- r_{k_\infty}+\sum_{j=k+1}^{k_\infty}\frac{1}{r_j^{\omega/4}}+\frac{1}{4r_k^{\omega/4}}\right),\]
\begin{equation}\label{e:sharpmomentbound5}
\begin{split}
&\tilde\E^{Q_N,0}\left(\I_{{\Ero_{k,x}}\cap\Ebd_{k,x}}\exp\left(\lambda(\phi(x)-\Sav_{k-1,-}(x,\phi))\right)\middle|\tilde\F_{Q_{r_{k-1}}^\com(x)}\right)\\
&=\exp\left(O_{\lambda_*}\left(\frac{1}{4r_k^{\omega/4}}\right)\right)
\tilde\E^{Q_N,0}\left(\exp (A)
\I_{\Ero_{k,x}}\exp\left(\lambda\Inc_k(x,\phi)\right)\middle|\tilde\F_{Q_{r_{k-1}}^\com(x)}\right)\\
&=\exp(A)
\tilde\E^{Q_N,0}\left(\I_{\Ero_{k,x}}\exp\left(\lambda\Inc_k(x,\phi)\right)\middle|\tilde\F_{Q_{r_{k-1}}^\com(x)}\right).
\end{split}
\end{equation}

We will use the results from Lemma \ref{l:exponential_moment_bounds} in order to estimate the exponential moment of $\Inc_{k-1}(x,\phi)=\Ind_{k-1}(x)+\I_{(\Ecp_{k-1,x})^\com}\Err_{k-1,x}(x)$ conditional on $\tilde\F_{Q_{r_k-1}^\com(x)}$. From the lemma we know that
\begin{equation}\label{e:sharpmomentbound6}
\tilde\E^{Q_N,0}\left(\exp\left(\lambda\Ind_{k-1}(x,\phi) \right)\middle|\tilde\F_{Q_{r_{k-1}}^\com(x)}\right)=\exp\left(\frac{\lambda^2\g}{2}+O_{\lambda_*}\left(\frac{1}{r_{k-1}^{\omega/2}}\right)\right).
\end{equation}
Combining this with \eqref{e:exponential_moment_bounds2} and H\"{o}lder's inequality, we also have that
\begin{align*}
&\tilde\E^{Q_N,0}\left(\exp\left(\lambda\Ind_{k-1}(x,\phi) \right)\left(\exp\left(\lambda\I_{(\Ecp_{{k-1},x})^\com}\Err_{{k-1},x}\right)-1\right)\middle|\tilde\F_{Q_{r_{k-1}}^\com(x)}\right)\\
&\qquad\le \left(\tilde\E^{Q_N,0}\left(\exp\left(2\lambda\Ind_{k-1}(x,\phi) \right)\middle|\tilde\F_{Q_{r_{k-1}}^\com(x)}\right)\right)^{1/2}\\
&\qquad \qquad \times \left(\tilde\E^{Q_N,0}\left(\left(\exp\left(\lambda\I_{(\Ecp_{k,x})^\com}\Err_{k-1}(x)\right)-1\right)^2\middle|\tilde\F_{Q_{r_{k-1}}^\com(x)}\right)\right)^{1/2}\\
&\qquad\le  \frac{C_{\lambda_*}}{r_{k-1}^{\omega/2}}.
\end{align*}
In combination with \eqref{e:sharpmomentbound6} this means that on the event $\Ero_{k-1,x}$ we have
\begin{equation}\label{e:sharpmomentbound7}
\begin{split}
\tilde\E^{Q_N,0}\left(\exp\left(\lambda\Inc_{k-1}(x,\phi)\right)\middle|\tilde\F_{Q_{r_{k-1}}^\com(x)}\right)&=\tilde\E^{Q_N,0}\left(\exp\left(\lambda\left(\Ind_{k-1}(x) +\I_{(\Ecp_{{k-1},x})^\com}\Err_k(x)\right)\right)\middle|\tilde\F_{Q_{r_{k-1}}^\com(x)}\right)\\
&=\exp\left(\frac{\lambda^2\g}{2}+O_{\lambda_*}\left(\frac{1}{r_{k-1}^{\omega/2}}\right)\right).
\end{split}
\end{equation}
Using once more the argument that lead to \eqref{e:sharpmomentbound3}, we also see
\[\left|\tilde\E^{Q_N,0}\left(\I_{(\Ero_{k,x})^\com}\exp\left(\lambda\Inc_k(x,\phi)\right)\middle|\tilde\F_{Q_{r_{k-1}}^\com(x)}\right)\right|\I_{{\Ero_{k-1,x}}}\le C_{\lambda_*}\exp(-c(\log r_{k-1})^{3/2}),\]
and so \eqref{e:sharpmomentbound7} implies that
\[\tilde\E^{Q_N,0}\left(\I_{\Ero_{k,x}}\exp\left(\lambda\Inc_k(x,\phi)\right)\middle|\tilde\F_{Q_{r_{k-1}}^\com(x)}\right)\I_{{\Ero_{k-1,x}}}=\exp\left(\frac{\lambda^2\g}{2}+O_{\lambda_*}\left(\frac{1}{4r_k^{\omega/4}}\right)\right).\]

Returning to \eqref{e:sharpmomentbound5}, this means that
\begin{equation}\label{e:sharpmomentbound8}
\begin{split}
&\tilde\E^{Q_N,0}\left(\I_{{\Ero_{k,x}}\cap\Ebd_{k,x}}\exp\left(\lambda(\phi(x)-\Sav_{k-1,-}(x,\phi))\right)\middle|\tilde\F_{Q_{r_{k-1}}^\com(x)}\right)\I_{{\Ero_{k-1,x}}}\\
&=\exp\left(\frac{\lambda^2\g (k_\infty-k+1)}{2}+O_{\lambda_*}\left(\log N- r_{k_\infty}+\sum_{j=k+1}^{k_\infty}\frac{1}{r_j^{\omega/4}}+\frac{3}{4r_k^{\omega/4}}\right)\right)\I_{{\Ero_{k-1,x}}}.
\end{split}
\end{equation}
We can now combine \eqref{e:sharpmomentbound2}, \eqref{e:sharpmomentbound3} and \eqref{e:sharpmomentbound8}, and find that
\begin{align*}
&\tilde\E^{Q_N,0}\left(\exp\left(\lambda(\phi(x)-\Sav_{k-1,-}(x,\phi))\right)\middle|\tilde\F_{Q_{r_{k-1}}^\com(x)}\right)\I_{{\Ero_{k-1,x}}}\\
&=\exp\left(\frac{\lambda^2\g (k_\infty-k+1)}{2}+O_{\lambda_*}\left(\log N- r_{k_\infty}+\sum_{j=k+1}^{k_\infty}\frac{1}{r_j^{\omega/4}}+\frac{3}{4r_k^{\omega/4}}\right)\right)\I_{{\Ero_{k-1,x}}}\\
&\qquad \qquad +O_{\lambda_*}\left(\exp(-c(\log r_{k-1})^{3/2})\right)\\
&=\exp\left(\frac{\lambda^2\g (k_\infty-k+1)}{2}+O_{\lambda_*}\left(\log N- r_{k_\infty}+\sum_{j=k}^{k_\infty}\frac{1}{r_j^{\omega/4}}\right)\right)\I_{{\Ero_{k-1,x}}}.
\end{align*}
This quickly implies \eqref{e:sharpmomentbound1} for $k-1$, thereby completing our induction and the proof of \eqref{e:sharpmomentbound_phi_Sk}.

The argument for \eqref{e:sharpmomentbound} is very similar. First of all, we can assume that $\log\dist(x,\partial Q_N)\ge C$ as otherwise the estimate follows from Brascamp-Lieb. Now if $\log\dist(x,\partial Q_N)\ge C$ then we can choose $k_0=\log\frac{N}{\dist(x,\partial Q_N)}\le k_\infty$. Once we know \eqref{e:sharpmomentbound_phi_Sk} for $k=k_0$, all that remains to be done is one conditioning step similar to the inductive step we just did, but using \eqref{e:exponential_moment_bounds4} instead of \eqref{e:exponential_moment_bounds1} and \eqref{e:exponential_moment_bounds2}.

The arguments for \eqref{e:sharpmomentbound_Sk} and \eqref{e:sharpmomentbound_Sk_Sk'} are even simpler. We prove
\[
\left|\log\tilde\E^{Q_N,0}\left(\exp\left(\lambda(\Sav_{k',+}(x,\phi)-\Sav_{k,-}(x,\phi))\right)\middle|\tilde\F_{Q_{r_{k}}^\com(x)}\right)-\frac{\lambda^2\g (k'-k)}{2}\right|\I_{{\Ero_{k,x}}}\le C_{\lambda_*}\sum_{j=k+1}^{k'}\frac{1}{r_j^{\omega/4}}\]
by reverse induction on $k$.
Here the case $k=k'-1$ follows directly from \eqref{e:sharpmomentbound7}, and the inductive step proceeds using exactly the same argument as before. This directly implies \eqref{e:sharpmomentbound_Sk_Sk'}, while for \eqref{e:sharpmomentbound_Sk} we again need to do one final inductive step using \eqref{e:exponential_moment_bounds4}.
\end{proof}

\section{Upper bound on the maximum}\label{s:upper}

\subsection{Preliminaries}
Our goal is to prove that with probability tending to 1 as $\Gamma\to\infty$ there is no $x$ such that $\phi(x)\ge\sqrt{\g}m_N+\Gamma$. In fact, we will prove a slightly stronger result involving not just $\phi(x)$ but also the square averages $\Sav_k(x,\phi)$.

For that purpose we introduce some notation. Throughout we fix $n=\log N$.
 Let $\Delta$ be a constant that is larger than the constants $C$ in Lemmas \ref{l:iterated_coupling}, \ref{l:exponential_moment_bounds} and 
 \ref{l:sharpmomentbound}, and fix it throughout this section. Let $\Gamma,\ell$ be two other constants that will be fixed later. For $x\in Q_N$ we define $k_0(x)=\left\lceil\log\frac{N}{\dist(x,\partial^+Q_N)}\right\rceil$ and $k_\infty=n-\Delta$. Moreover we set $k_\infty'=k_\infty-\ell$.
We consider the barrier event
\begin{equation}
  \label{eq-upperbarrier}
\Eup_x=\left\{\phi\colon \Sav_{j,+}(x,\phi)\le\Bup(j)+\Gamma\ \forall j\in\{k_0(x),\ldots,k_\infty'\}\right\}
\end{equation}
where
\[\Bup(j):=\frac{\sqrt{\g}(m_N-2\Delta)}{k_\infty}j+(j\wedge(k_\infty-j))^{2/5},\]
and our goal is to prove that the probability of the event
\[\bigcap_{x\in Q_N}\Eup_x\cap\{\phi\colon \phi(x)\le \sqrt{\g}m_N+\Gamma\} \]
tends to 1 as $\Gamma\to\infty$, uniformly in $n$. We will prove this arguing scale-by-scale, beginning with $k=1$. In fact, we will use two different arguments. For $k\le k_\infty-(\log n)^3$, say, there is a rather easy argument that we will give in the next section. The remaining case $k_\infty-(\log n)^3<k\le k_\infty'$ is more complicated, and it will be treated in the section thereafter.

We  introduce some more notation. We define for $k\in\{k_0(x),\ldots k_\infty\}$
\[
\Eup_{\le k,x}=\left\{\phi\colon \Sav_{j,+}(x,\phi)\le\Bup(j)+\Gamma\ \forall j\in\{k_0(x),\ldots,k\wedge k_\infty'\}\right\}.
\]
We also define the events
\begin{align*}
&\Ecp_x=\bigcap_{k_0(x)\le j\le k_\infty-1}\Ecp_{j,x},\qquad
\Ero_x=\bigcap_{k_0(x)\le j\le k_\infty}\Ero_{j,x},\qquad
\Ebd_x=\bigcap_{k_0(x)\le j\le k_\infty}\Ebd_{j,x},
\end{align*}
and, for $k\in\{k_0(x),\ldots k_\infty\}$, the events
\begin{align*}
&\Ecp_{\le k,x}=\bigcap_{k_0(x)\le j\le k}\Ecp_{j,x},\qquad
\Ero_{\le k,x}=\bigcap_{k_0(x)\le j\le k}\Ero_{j,x},\qquad
\Ebd_{\le k,x}=\bigcap_{k_0(x)\le j\le k}\Ebd_{j,x}.
\end{align*}

For technical reasons, it will be convenient for us to know that the maximum of $\phi$ is unlikely to occur very close to the boundary. In fact, we have

\begin{lemma}\label{l:pointsnearbdry}
For any $a>0$ we have
\[\tilde\PP^{Q_N,0}\left(\max_{\substack{x\in Q_N\\\dist(x,\partial^+Q_N)\le \frac{N}{n^{5/2}}}}\phi(x)\ge \sqrt{\g}m_N\right)\le\frac{1}{n}.\]
\end{lemma}
\begin{proof}
Using the sharp tail bound from Lemma \ref{l:quitesharptailbound} as well as a union bound we obtain
\[\tilde\PP^{Q_N,0}\left(\max_{\substack{x\in Q_N\\\dist(x,\partial^+Q_N)\le N/n^{5/2}}}\phi(x)\ge m_N\right)\le CN\frac{N}{n^{5/2}}\exp\left(-2n+\frac32\log n\right)\le\frac{1}{n}.
\]
\end{proof}

In view of this lemma, it will suffice to bound the maximum of $\phi$ over points $x$ where $k_0(x)\le \frac52\log n\ll\frac{n}{2}$.

\subsection{Early barrier}
In this subsection we will prove that the event $\Eup_{\le k_\infty-(\log n)^3,x}$ occurs with probability tending to 1 as $\Gamma\to\infty$, uniformly in all sufficiently large $n$. As our barrier $(j\wedge(k_\infty-j))^{2/5}$ is curved upwards, in the regime $j\le k_\infty-(\log n)^3$ it is bounded below by $c\log j$, 
with $c$ large at least for all  $j$ large. This gives us enough room to argue that the barrier is unlikely to be crossed.

For some fixed $k\le k_\infty-(\log n)^3$ we cannot simply use a union bound to estimate the probability of crossing the barrier at $k$, as the factor $O(N^2)$ is way too large. Instead we will consider a set $B_k$ of representatives and argue that
 the $\Sav_{k,+}(x,\phi)$ for $x\in B_k$ are likely to stay well below the barrier. If we can also show that the maximum of $\Sav_{k,+}(\cdot,\phi)$ over $Q_N$ is likely to be not much larger than the maximum over $B_k$, we can obtain the desired lower bound on the probability of $\Eup_{\le k,x}$.

Regarding the choice of the representatives, note that $\Sav_{k,+}(x,\phi)$ fluctuates on scale $r_k$. Because of this, a natural choice would be to take $B_k=r_k\Z^2\cap Q_N$. However, with this choice the maximal error from the fluctuations of $\Sav_{k,+}(x,\phi)$ would be of order $\sqrt{k}$, which is too large for our purposes. Instead we will choose them as $B_k=\frac{r_k}{k}\Z^2\cap Q_N$. This makes the error from the fluctuations significantly smaller, while the error from bounding the maximum over the representatives increases only moderately.

We now make this rigorous. We begin by stating an estimate on the fluctuations.

\begin{lemma}\label{l:fluctuations_uncond}
There is a constant $C$ such that for each $\Gamma>0$ and any $1\le k\le k_\infty$ we have
\begin{equation}\label{e:fluctuations_uncond}
\begin{split}
&\tilde\PP^{Q_N,0}\left(\exists x,x'\in Q_N\colon k_0(x)\wedge k_0(x')\ge k, |x-x'|\le \frac{r_k}{k}, \left|\Sav_{k,+}(x,\phi)-\Sav_{k,+}(x',\phi)\right|\ge \frac{\Gamma}{2}\right)\\
&\le C\exp\left(2k+2\log k-\frac{\Gamma^2 k}{C}\right).
\end{split}
\end{equation}
\end{lemma}
Note in particular that for $\Gamma$ large enough the right-hand side becomes summable in $k$.
\begin{proof}
This is essentially a chaining argument. Let
\begin{equation}\label{e:def_Ak}
A_{k}=\{x\in Q_N\colon k_0(x)\ge k\}=Q_{N-r_k}.
\end{equation}
By the Brascamp-Lieb inequality we have that 
\begin{equation}
\label{eq-Dec17a}
\tilde\E^{Q_N,0}\left(\exp\left(\lambda(\Sav_{k,+}(x,\phi)-\Sav_{k,+}(x',\phi))\right)\right)\le \exp\left(C\lambda^2\Var_{\mathrm{DGFF}}(\Sav_{k,+}(x,\phi)-\Sav_{k,+}(x',\phi)) \right),
\end{equation}
for any $x,x'\in A_k$.
Note that $\Sav_{k,+}(x,\phi)-\Sav_{k,+}(x',\phi) = \langle \phi,  g_+ - g_-\rangle$,  where for $r= |y-x|_\infty$ and $r'= |y-x'|_\infty$ ,
\[
 g_+(y) = \frac{\eta_{r_{k}^{1-\omega}}(|y-x|_\infty-r_{k,+})a_{Q_r(x)}(x,y)}{\sum_{r\in\N}\eta_{r_{k}^{1-\omega}}(r-r_{k,+})}
\]
and
\[
 g_-(y) = \frac{\eta_{r_{k}^{1-\omega}}(|y-x'|_\infty-r_{k,+})a_{Q_{r'}(x')}(x',y)}{\sum_{r\in\N}\eta_{r_{k}^{1-\omega}}(r-r_{k,+})}.
\]
We claim that
\begin{equation}\label{e.twocenter}
\begin{split}
\Var_{\mathrm{DGFF}}(\Sav_{k,+}(x,\phi)-\Sav_{k,+}(x',\phi))
&=
\sum_{y,z\in Q_N} (g_+(z) - g_-(z)) G_{Q_N}(y,z) (g_+(y) - g_-(y))\\
&\leq
C\biggl(\frac{|x-x'|}{r_k}\wedge1\biggr).
\end{split}
\end{equation}
which together with \eqref{eq-Dec17a} implies
\begin{equation}\label{e.twocenter_expmom}
\tilde\E^{Q_N,0}\left(\exp\left(\lambda(\Sav_{k,+}(x,\phi)-\Sav_{k,+}(x',\phi))\right)\right)\le \exp\left(C\lambda^2\left(\frac{|x-x'|}{r_k}\wedge1\right)\right).
\end{equation}

To see the claim \eqref{e.twocenter} ,  it suffices to check the case $|x-x'| < \delta r_k$,  for some small $\delta>0$.  We first sum over $z$ (resp. $y$) at the boundary $\partial^+ Q_r(x)$ (resp.  $\partial^+ Q_s(x)$).  Denote by $r' = |x'-z|_\infty$ and $s' = |x'-y|_\infty$, we need to estimate
\begin{equation*}
\sum_{z\in \partial^+ Q_r(x),  y\in \partial^+ Q_s(x)}  (a_{Q_r(x)}(x,z) -  a_{Q_{r'}(x')}(x',z) )G_{Q_N}(y,z)(a_{Q_s(x)}(x,y) -  a_{Q_{s'}(x')}(x',y) ),
\end{equation*}
where $a_B$ is as in \eqref{eq-a}.
We apply the Harnack inequality to conclude that for $r> 4|x-x'|$,
\begin{equation*}
|a_{Q_r(x)}(x,z) -  a_{Q_{r'}(x')}(x',z)|
\le
\frac Cr \biggl(\frac{|x-x'|}{r}\wedge1 \biggr).
\end{equation*}
Together with the fact that
\begin{equation*}
\sum_{z\in \partial^+ Q_r(x),  y\in \partial^+ Q_s(x)} \frac 1{rs} G_{Q_N}(y,z)
< C,
\end{equation*}
we obtain,  for each $r,s\in [r_{k,+} - r_k^{1-\omega}, r_{k,+} + r_k^{1-\omega} ]$,
\begin{align*}
&\sum_{z\in \partial^+ Q_r(x),  y\in \partial^+ Q_s(x)}  (a_{Q_r(x)}(x,z) -  a_{Q_{r'}(x')}(x',z) )G_{Q_N}(y,z)(a_{Q_s(x)}(x,y) -  a_{Q_{s'}(x')}(x',y) )\\
&\qquad\qquad 
\leq
C  \biggl(\frac{|x-x'|}{r_k}\wedge1 \biggr).
\end{align*}
Summing over $r,s$ we conclude \eqref{e.twocenter}.

Applying \eqref{e.twocenter_expmom} and
the exponential Markov inequality
 thus implies the Gaussian tail bounds
\begin{equation}\label{e:fluctuations_uncond1}
\tilde\PP^{Q_N,0}\left(\left|\Sav_{k,+}(x,\phi)-\Sav_{k,+}(x',\phi)\right|\ge s\right)\le \exp\left(-\frac{s^2r_k}{C|x-x'|}\right).
\end{equation}
We can bound
\begin{align*}
&\tilde\PP^{Q_N,0}\left(\exists x,x'\in A_k\colon |x-x'|\le \frac{r_k}{k}, \left|\Sav_{k,+}(x,\phi)-\Sav_{k,+}(x',\phi)\right|\ge \frac{\Gamma}{2}\right)\\
&\le\sum_{l=0}^{\log_2 (r_k/k)}\sum_{\substack{x_l\in A_k\cap 2^l\Z^d\\x_{l+1}\in A_k\cap 2^{l+1}\Z^d)\\|x_l-x_{l+1}|\le 2^l}}\tilde\PP^{Q_N,0}\left(\left|\Sav_{k,+}(x_l,\phi)-\Sav_{k,+}(x_{l+1},\phi)\right|\ge \frac{2^{l/4}k^{1/4}\Gamma}{32r_k^{1/4}}\right)\\
&\le\sum_{l=0}^{\log_2 (r_k/k)} C\frac{N^2}{2^{2l}}\exp\left(-\frac{k^{1/2}r_k^{1/2}\Gamma^2}{C2^{l/2}}\right)
\le C\frac{N^2k^2}{r_k^2}\exp\left(-\frac{k\Gamma^2}{C}\right)\\
&\le C\exp\left(2k+2\log k-\frac{\Gamma^2k}{C}\right),
\end{align*}
which is \eqref{e:fluctuations_uncond}.
\end{proof}

Using Lemma \ref{l:fluctuations_uncond}, it is now relatively straightforward to bound the probability of $(\Eup_{\le k_\infty-(\log n)^3,x})^\com$.

\begin{lemma}\label{l:upper_barrier_early}
There is a constant $C$ such that for each $\Gamma\ge C$ and for each sufficiently large $N$ we have
\begin{equation}\label{e:upper_barrier_early}
\tilde\PP^{Q_N,0}\left(\bigcup_{x\in Q_N}(\Eup_{\le k_\infty-(\log n)^3,x})^\com\right)\le C\exp\left(-\frac{\Gamma}{C}\right).
\end{equation}
\end{lemma}
\begin{proof}
Note that for $N$ sufficiently large we have $(\log n)^3\ge k_\infty-k_\infty'$ and $3\sqrt{\g}\log (2+k)\le(k\wedge(k_\infty-k))^{2/5}$ for any $k\le k_\infty-(\log n)^3$.

Hence
\begin{equation}\label{e:upper_barrier_early1}
\begin{split}
&\tilde\PP^{Q_N,0}\left(\bigcup_{x\in Q_N}(\Eup_{\le k_\infty-(\log n)^3,x})^\com\right)\\
&\!\!
\le\! \sum_{k=1}^{k_\infty-(\log n)^3}\!\!\tilde\PP^{Q_N,0}\left(\exists x,x'\in Q_N\colon k_0(x)\wedge k_0(x')\ge k, |x-x'|\le \frac{r_k}{k}, \left|\Sav_{k,+}(x,\phi)-\Sav_{k,+}(x',\phi)\right|\ge \frac{\Gamma}{2}\right)\\
&\qquad+ \sum_{k=1}^{k_\infty-(\log n)^3}\tilde\PP^{Q_N,0}\left(\exists x\in Q_N\cap \frac{r_k}{k}\Z^2\colon \Sav_{k,+}(x,\phi)\ge \frac{\sqrt{\g}(m_N-2\Delta)}{k_\infty}k+3\sqrt{\g}\log (2+k)+\frac{\Gamma}{2}\right).
\end{split}
\end{equation}
The summands of the first sum can be bounded using Lemma \ref{l:fluctuations_uncond}. For the summands of the second sum, we first observe that
\[\frac{\sqrt{\g}(m_N-2\Delta)}{k_\infty}k+3\sqrt{\g}\log (2+k)
\ge \sqrt{\g}\left(2k-\frac{3\log n}{4n}k+3\log (2+k)\right)-C
\ge \sqrt{\g}\left(2k+2\log k-C\right),\]
where we used the fact that $t\mapsto\frac{\log t}{t}$ is a decreasing function for $t\ge\e$,  and replace $n$ by $k$ to obtain the last inequality.
Moreover, by \eqref{e:sharpmomentbound_Sk} we have
\[\tilde\E^{Q_N,0}\left(\exp\left(\frac{2}{\sqrt{\g}}\Sav_{k,+}(x,\phi)\right)\right)\le C\exp(2k),\]
and so by the exponential Markov inequality and a union bound
\begin{align*}
&\tilde\PP^{Q_N,0}\left(\exists x\in Q_N\cap \frac{r_k}{k}\Z^2\colon \Sav_{k,+}(x,\phi)\ge \frac{\sqrt{\g}(m_N-2\Delta)}{k_\infty}k+3\sqrt{\g}\log (2+k)+\frac{\Gamma}{2}\right)\\
&\le C\frac{N^2k^2}{r_k^2}\exp(\left(2k-\frac{2}{\sqrt{\g}}\sqrt{\g}\left(2k+2\log k+\frac{\Gamma}{2}-C\right)\right)\\
&\le C\exp(2k+2\log k)\exp\left(-2k-4\log k-\frac{\Gamma}{\sqrt{\g}}\right)\\
&\le \frac{C}{k^2\e^{c\Gamma}}.
\end{align*}
Using this estimate and also Lemma \ref{l:fluctuations_uncond}, we conclude from \eqref{e:upper_barrier_early1} that
\[
\tilde\PP^{Q_N,0}\left((\Eup_{\le k_\infty-(\log n)^3,x})^\com\right)\le C\sum_{k=1}^{k_\infty-(\log n)^3}\exp\left(2k+2\log k-\frac{\Gamma^2 k}{C}\right)+\frac{1}{k^2\e^{c\Gamma}}.
\]
We can estimate $\log k\le k$, and so as soon as $\Gamma$ is large enough so that $4-\frac{\Gamma^2}{C}\le-\Gamma$, the right-hand side is bounded by $\frac{C}{\e^{c\Gamma}}$, as claimed.
\end{proof}

\subsection{Late barrier}

On small scales $k_\infty-(\log n)^3\le k\le k_\infty'$, our general strategy is the same as for the large scales: We consider a set of representatives $B_k$ at some scale slightly smaller than $r_k$, and argue that on the one hand the field $\Sav_{k,+}(x,\phi)$ for $x\in B_k$ is likely to stay well below the barrier, and that on the other the fluctuations are small enough so that actually $\Sav_{k,+}(x,\phi)$ does not cross the barrier for any $x\in Q_N$.

However, we now have much less room to spare below the barrier. In particular, a global estimate for the fluctuations as in Lemma \ref{l:fluctuations_uncond} is already too rough. Instead we will now only estimate the fluctuations near those representatives where the square average process $\Sav_{k',+}(x,\phi)$ for some $k'<k$ is already relatively high. These points are rare, and so we get a much better estimate than in Lemma \ref{l:fluctuations_uncond}.

In fact, we will proceed by induction on $k$, with the cases $k\le n-(\log n)^3$ already established in the previous subsection. By the inductive assumption, we can assume the square average process $\Sav_{k',+}(x,\phi)$ for $k'\le k$ has not crossed the barrier, and this will imply that it is unlikely to end very close to the barrier. So there are only few representatives for which we need a sharp estimate for the fluctuations.

This means, though, that we need to estimate fluctuations conditional on the fact that the square averages on larger scales are high. That is, we need a version of Lemma \ref{l:fluctuations_uncond} that is valid conditionally on the field outside some larger box. We will begin by proving the corresponding estimate. Note that the lengthscale on which we estimate the fluctuations is slightly different than in the previous section. A good choice is $\frac{r_k}{(n-k)^2}$, as it allows us to estimate the error term in \eqref{e:fluctuations_cond_error} below.

\begin{lemma}\label{l:fluctuations_cond}
Let $1\le k'< k\le k_\infty$, $x\in Q_N$ with $\dist(x,\partial^+Q_N)\ge r_{k'}$. Then for any $|\lambda|\le\lambda_*$ we have
\begin{equation}\label{e:fluctuations_cond}
\begin{split}
&\log\tilde\E^{Q_N,0}\left(\exp\left(\lambda\left(\max_{x'\in Q_{r_k/(n-k)^2}(x)} \Biggl(\Sav_{k,+}(x',\phi)-\Sav_{k',-}(x,\phi)\Biggr)\right)\right)\middle|\tilde\F_{Q_{r_{k'}}(x)^\com}\right)\I_{\Ero_{k',x}}\\
&\qquad\le \frac{\lambda^2\g(k-k')}{2}+C_{\lambda_*}.
\end{split}
\end{equation}
Moreover, for any $\lambda\in\R$ we have
\begin{equation}\label{e:fluctuations_cond_BL}
\begin{split}
&\log\tilde\E^{Q_N,0}\left(\exp\left(\lambda\left(\max_{x'\in Q_{r_k/(n-k)^2}(x)}\Biggl(\Sav_{k,+}(x',\phi)-\Sav_{k',-}(x,\phi)\Biggr)\right)\right)\middle|\tilde\F_{Q_{r_{k'}}(x)^\com}\right)\I_{\Ero_{k',x}}\\
&\qquad\le C\lambda^2(k-k')+C.
\end{split}
\end{equation}
\end{lemma}

\begin{proof}
We begin with \eqref{e:fluctuations_cond}. This estimate is similar to \eqref{e:sharpmomentbound_Sk_Sk'}, and in fact we will prove it by a reverse induction on $k'$, just like \eqref{e:sharpmomentbound_Sk_Sk'}. The induction step will be completely analogous, so the only difficulty will be to establish the case $k'=k-1$. That is, we need to prove
\begin{equation}\label{e:fluctuations_cond1}
\log\tilde\E^{Q_N,0}\left(\exp\left(\lambda\left(\max_{x'\in Q_{r_k/(n-k)^2}(x)} \Biggl(\Sav_{k,+}(x',\phi)-\Sav_{k-1,-}(x,\phi)\Biggr)\right)\right)\middle|\tilde\F_{Q_{r_{k-1}}(x)^\com}\right)\I_{\Ero_{k-1,x}}\le C_{\lambda_*}.
\end{equation}
This estimate in turn follows easily from the Gaussian tail bound
\begin{equation}\label{e:fluctuations_cond2}
\tilde\PP^{Q_N,0}\left(\exists x'\in Q_{r_k/(n-k)^2}(x)\colon \left|\Sav_{k,+}(x',\phi)-\Sav_{k-1,-}(x,\phi)\right|\ge t\middle|\tilde\F_{Q_{r_{k-1}}(x)^\com}\right)\I_{\Ero_{k-1,x}}\le C\exp\left(-\frac{t^2}{C}\right).
\end{equation}
The estimate \eqref{e:fluctuations_cond2} is similar to the statement of Lemma \ref{l:fluctuations_uncond}. However here we have the additional difficulty that the conditional expectation of $\Sav_{k,+}(x',\phi)-\Sav_{k-1,-}(x,\phi)$ need not be zero.
In order to control this expectation, we can proceed similarly as in the proof of Lemma \ref{l:badevents}. Namely, on the event $\Ero_{k-1,x}$ we know that
\[\osc_{y\in Q_{r_{k-1}}(x)}\tilde\E^{Q_N,0}\left(\phi(y)\mid\tilde\F_{Q_{r_{k-1}}^\com(x)}\right)\le 4(\log r_{k-1})^2.\]
By Theorem \ref{t:miller_mean}, \[\tilde\E^{Q_N,0}\left(\phi(y)\mid\tilde\F_{Q_{r_{k-1}}^\com(x)}\right)=\hat\phi(y)+O\left(\frac{1}{r_{k-1}^\omega}\right)\text{ on }Q_{r_{k-1,-}}(x),\]
where $\hat\phi$ is the harmonic function that agrees with $\tilde\E^{Q_N,0}\left(\phi(\cdot)\mid\tilde\F_{Q_{r_{k-1}}^\com(x)}\right)$ on $\partial^+Q_{r_{k-1,-}}(x)$.
The maximum principle implies $\osc_{y\in Q_{r_{k'},-}(x)}\hat\phi(y)\le 4(\log r_{k-1})^2$, and by gradient estimates for discrete harmonic functions (cf. \cite[Theorem 6.3.8 (a)]{LL10}), we have that
\begin{equation}\label{e:fluctuations_cond_error}
\left|\Sav_{k,+}(x',\hat\phi)-\Sav_{k-1,-}(x,\hat\phi)\right|\le C\frac{|x-x'|}{r_{k-1}}(\log r_{k-1})^2\le C\frac{r_k(\log r_{k-1})^2}{(n-k)^2r_{k-1}}\le C,
\end{equation}
where in the last step we used that $n-k=\log r_k$.
Hence also
\begin{equation}\label{e:fluctuations_cond3}
\left|\tilde\E^{Q_N,0}\left(\Sav_{k,+}(x',\phi)-\Sav_{k-1,+}(x,\phi)\mid\tilde\F_{Q_{r_{k-1}}^\com(x)}\right)\right|\I_{\Ero_{k-1,x}}\le C+\frac{C}{r_{k-1}^\omega}\le C.
\end{equation}
This is a bound on the conditional expectation of $\Sav_{k,+}(x',\phi)-\Sav_{k,-}(x,\phi)$. We will also need to control the fluctuations. 
To that end, we claim that, with 
\[V_{x,x'}:=\left|\Sav_{k,+}(x',\phi)-\Sav_{k-1,-}(x,\phi)-\tilde\E^{Q_N,0}\left(\Sav_{k,+}(x',\phi)-\Sav_{k-1,-}(x,\phi)\middle|\tilde\F_{Q_{r_{k-1}}^\com(x)}\right)\right|,\]
\begin{equation}\label{e:fluctuations_cond4}
  \tilde\PP^{Q_N,0}\left(\max_{ x'\in Q_{r_k/(n-k)^2}(x)}V_{x,x'}\ge s\middle|\tilde\F_{Q_{r_{k-1}}^\com(x)}\right)
\I_{\Ero_{k-1,x}}
\le C\exp\left(-\frac{(n-k)^2s^2}{C}\right)\le C\exp\left(-\frac{s^2}{C}\right).
\end{equation}
It is clear that the combination of \eqref{e:fluctuations_cond3} and \eqref{e:fluctuations_cond4} implies \eqref{e:fluctuations_cond2} (and then also \eqref{e:fluctuations_cond1}), and so it remains to verify \eqref{e:fluctuations_cond4}.
This is a chaining argument very similar to the one in the proof of Lemma \ref{l:fluctuations_uncond}, using instead the exponential moment estimate \eqref{e:fluctuations_cond3}.

The argument for \eqref{e:fluctuations_cond_BL} is similar, but actually less complicated. Again we proceed by reverse induction on $k'$, with the case $k'=k-1$ still following from \eqref{e:fluctuations_cond2}. For the inductive step we proceed as before, just using the Brascamp-Lieb inequality instead of the sharp moment estimates from Lemma \ref{l:sharpmomentbound}.
\end{proof}

In practice we will use Lemma \ref{l:fluctuations_cond} in the form of a Gaussian tail bound on the fluctuations, that we state separately.
\begin{lemma}\label{l:fluctuations_cond_tail_bound}
In the setting of Lemma \ref{l:fluctuations_cond} we have, for any $|\lambda|\le\lambda_*$ and any $\eps>0,t\geq 0$, the estimate
\begin{equation}\label{e:fluctuations_cond_tail_bound}
\begin{split}
&\tilde\PP^{Q_N,0}\left(\exists x'\in Q_{r_k/(n-k)^2}(x)\colon \left|\Sav_{k,+}(x',\phi)-\Sav_{k',-}(x,\phi)\right|\ge t\middle|\tilde\F_{Q_{r_{k'}}(x)^\com}\right)\I_{\Ero_{k',x}}\\
&\le C_{\lambda_*,\eps}\exp\left(\frac{(\lambda^2+\eps)\g(k-k')}{2}-\lambda t-\frac{t^2}{C(k-k')^3}\right).
\end{split}
\end{equation}
\end{lemma}
\begin{proof}
As a consequence of \eqref{e:fluctuations_cond} and the exponential Markov inequality we have
\begin{equation}\label{e:fluctuations_cond_tail_bound1}
\begin{split}
&\tilde\PP^{Q_N,0}\left(\exists x'\in Q_{r_k/(n-k)^2}(x)\colon \left|\Sav_{k,+}(x',\phi)-\Sav_{k',-}(x,\phi)\right|\ge t\middle|\tilde\F_{Q_{r_{k'}}(x)^\com}\right)\I_{\Ero_{k',x}}\\
&\le C_{\lambda_*}\exp\left(\frac{\lambda^2\g(k-k')}{2}-\lambda t\right).
\end{split}
\end{equation}
On the other hand, \eqref{e:fluctuations_cond_BL} implies that
\begin{equation}\label{e:fluctuations_cond_tail_bound2}
\begin{split}
&\tilde\PP^{Q_N,0}\left(\exists x'\in Q_{r_k/(n-k)^2}(x)\colon \left|\Sav_{k,+}(x',\phi)-\Sav_{k',-}(x,\phi)\right|\ge t\middle|\tilde\F_{Q_{r_{k'}}(x)^\com}\right)\I_{\Ero_{k',x}}\\
&\le \exp\left(-\frac{t^2}{C'(k-k')}\right),
\end{split}
\end{equation}
for some constant $C'$. We choose the constant $C$ in \eqref{e:fluctuations_cond_tail_bound} as $2C'$. Now for $t\le \sqrt{C'\eps\g}(k-k')^2$ we have
\[\frac{\eps\g(k-k')}{2}\ge\frac{t^2}{2C'(k-k')^3}\]
and so
\eqref{e:fluctuations_cond_tail_bound} follows directly from \eqref{e:fluctuations_cond_tail_bound1}. So it remains to consider the case $t> \sqrt{C'\eps\g}(k-k')^2$. If $k-k'\ge2\lambda\sqrt{\frac{C'}{\eps\g}}$ we can estimate
\[\lambda t\le \frac{\lambda t^2}{\sqrt{C'\eps\g}(k-k')^2}\le \frac{t^2}{2C'(k-k')},\]
while if $k-k'\le2\lambda\sqrt{\frac{C'}{\eps\g}}$ we can use that
\[\lambda t\le \frac{t^2}{2C'(k-k')}+\frac{C'\lambda^2(k-k')}{2}\le \frac{t^2}{2C'(k-k')}+\sqrt{\frac{C'^3}{\eps\g}}\lambda_*^3.\]
So in any case, if $t> \sqrt{C'\eps\g}(k-k')^2$  we have
\[\frac{(\lambda^2+\eps)\g(k-k')}{2}-\lambda t-\frac{t^2}{2C'(k-k')^3}\ge -\frac{ t^2}{2C'(k-k')}-C_{\lambda_*,\eps}-\frac{t^2}{2C'(k-k')^3}\ge -\frac{ t^2}{C'(k-k')}-C_{\lambda_*,\eps}, \]
so that \eqref{e:fluctuations_cond_tail_bound} in this remaining case follows from \eqref{e:fluctuations_cond_tail_bound2}.
\end{proof}

We will also rely on a barrier estimate for random walks,  whose proof appears 
 in Section \ref{s:ballot_thms}.
 \begin{lemma}
 \label{l:upper_barrier}
Let $m\in\N$ and $\gamma>0$. Let $(X_j)_{j=0}^m$ be independent random variables possessing a uniformly bounded density,  $C<\infty$ such that $\left|\E_\Q[X_j]\right| \le C\exp(\gamma(j-m))$ and $\left|\E_\Q[X_j^2]-\g\right|\le C\exp(\gamma(j-m))$ for any $j\ge 0$.  We further assume that there exists $\lambda_*> 0$ and $C^*<\infty$  such that $\E\exp(\lambda_* X_j^2)< C^*$ for all $|\lambda|<\lambda_*$ and all $j\ge0$.  Let $a\ge0$, $t\ge0$,  and $0\le k_0\le \frac m2$, and define $\Sigma_k = \sum_{j=k_0}^{k-1} X_j$.

Then there exists some $\ell\ge1$ fixed and a constant $C_{\gamma,\ell}\ge 0$ depending on $\gamma$ and $\ell$ only such that
\begin{align}
\label{eq-230124a}
&\PP \left(-a+\Sigma_j \le (j\wedge (m-j))^{2/5} \ \forall j\in\{k_0\vee\ell,\ldots m-\ell\}, \Sigma_{m-1} \in [ -t-1, -t] \right) \nonumber\\
&\qquad \le \frac{C_{\gamma,\ell}(1+a+k_0^{2/5})(1+t)}{m^{3/2}}.
\end{align}
Moreover for any $k\in\{1,\ldots,m-1\}$ and $a\geq 0$ we have the estimates
\begin{align}
\label{eq-230124b}
&\PP \left(-a+\Sigma_j \le (j\wedge (m-j))^{2/5} \ \forall j\in\{k_0,\ldots k-1\}, \Sigma_{k-1} \in (k\wedge (m-k))^{2/5}+[ -t-1, -t] \right)\nonumber\\
&\qquad \qquad \qquad \le \frac{C_{\gamma,\ell}(1+a+k_0^{2/5})(1+t)}{(k-k_0(x))^{3/2}}.
  \end{align}
 \end{lemma}

As explained in the beginning of the subsection, we will estimate the probability of a late crossing of the barrier by induction on the first time when it is crossed. The following lemma, which is 
the main result of the present section, is an estimate on the probability that our upper barrier is crossed for the first time at some time $k\in\{k_\infty-(\log n)^3,\ldots, k_\infty'\}$. Here we restrict ourselves to points that have at least a small distance to the boundary (which is fine in view of Lemma \ref{l:pointsnearbdry}). We will carefully estimate the probability that the square average process crosses the barrier for the first time at scale $k+1$. The proof will be based on calculations of exponential moments, but we will need to keep track of the correct polynomial prefactors (in particular the one from Lemma \ref{l:upper_barrier}) to get the precise result.

A major challenge in the proof is that the decomposition from Lemma \ref{l:iterated_coupling} involves various good events at different scales, and we do not know if these all occur. We address this by introducing in the proof a (random) scale $\kmax<k$ such that, roughly speaking, all required good events up to that scale occur. Up to scale $\kmax$, Lemma \ref{l:iterated_coupling} allows us to compare our square average process with a random walk with independent increments, and so the barrier estimate from Lemma \ref{l:upper_barrier} (and its tilted version) give us very fine information on the behavior of the process. We give the technical details in Step 3 of the proof below.

From scale $\kmax$ onwards, Lemma \ref{l:upper_barrier} is no longer available, and so the best we can do is to apply the conditional fluctuation estimate from Lemma \ref{l:fluctuations_cond_tail_bound} between $\kmax$ and $k$. However, this fluctuation estimate is less sharp, and we lose a summand $\frac{\eps\g(k-\kmax)}{2}$ in the exponent when applying it. This summand is problematic if $k-\kmax\gg1$, i.e. if an early good event has failed to occur. But the bounds from Lemma \ref{l:badevents} allow us to the estimate the probability that this happens, and it turns out that we gain a summand $-\frac{\omega(n-\kmax)}{2}$ in the exponent, which is dominant once we choose $\eps$ small enough. The lengthy argument is explained in Step 2 of the proof.

Here now is the result.
\begin{lemma}\label{l:upper_barrier_late}
Let $A_{(\log n)^{5/2}}$ be as in \eqref{e:def_Ak}, and let $k\in\{k_\infty-(\log n)^3,\ldots k_\infty'\}$. Then we have
\begin{equation}\label{e:upper_barrier_late}
\tilde\PP^{Q_N,0}\left(\bigcap_{x\in A_{(\log n)^{5/2}}}\Eup_{\le k,x} \cap \bigcup_{x\in A_{(\log n)^{5/2}}}(\Eup_{\le k+1,x})^\com\right)\le C\exp\left(-c(n-k+1)^{2/5}-c\Gamma\right).
\end{equation}
\end{lemma}
Note in particular that the right-hand side of \eqref{e:upper_barrier_late} is summable in $k$. So the lemma together with Lemma \ref{l:upper_barrier_early} implies that with high probability the event $\Eup_x$ occurs for all $x\in A_{(\log n)^{5/2}}$.

\begin{proof}
Fix some $k\in\{k_\infty-(\log n)^3,\ldots k_\infty'\}$. Similarly as in the proof of Lemma \ref{l:upper_barrier_early}, we will work with the set of representatives $\frac{r_k}{(n-k)^2}\Z^2\cap Q_N$, and argue that $\Sav_k(\cdot,\phi)$ is likely to be below the barrier at all representatives, and then use Lemma \ref{l:fluctuations_cond} to conclude on all of $A_{(\log n)^{5/2}}$. Indeed, we can estimate that
\begin{equation}\label{e:upper_barrier_late1}
\begin{split}
&\tilde\PP^{Q_N,0}\left(\bigcap_{x\in A_{(\log n)^{5/2}}}\Eup_{\le k,x} \cap \bigcup_{x\in A_{(\log n)^{5/2}}}(\Eup_{\le k+1,x})^\com\right)\\
&\le\tilde\PP^{Q_N,0}\left(\exists x\in\frac{r_k}{(n-k)^2}\Z^2, x'\in Q_{r_k/(n-k)^2}(x)\colon \phi\in\Eup_{\le k,x} \cap (\Eup_{\le k+1,x'})^\com\right)\\
&\le \sum_{x\in r_k/(n-k)^2\Z^2}\tilde\PP^{Q_N,0}\left(\exists x'\in Q_{r_k/(n-k)^2}(x)\colon \phi\in\Eup_{\le k,x} \cap (\Eup_{\le k+1,x'})^\com\right),
\end{split}
\end{equation}
and we will estimate each of the summands on the right-hand side.

\emph{Step 1: Introduction of the random scale $\kmax$}\\
Fix some $x\in \frac{r_k}{(n-k)^2}\Z^2$. Our goal is to use the barrier estimate from Lemma \ref{l:upper_barrier} on as many scales as possible. 
Toward this end we define the random index $\kmax$ as, roughly speaking, the maximal $k'$ such that $\Ero_{\le k',x}\cap \Ebd_{\le k',x}\cap \Ecp_{\le k'-1,x}$ occur. Up to scale $\kmax$ we will use the estimate from Lemma \ref{l:upper_barrier}, while from scale $\kmax+1$ onwards we will use the estimate from Lemma \ref{l:fluctuations_cond}. Moreover, if $\kmax<k$, we gain an additional small factor from the fact that one of $\Ero_{\kmax+1,x}$, $\Ebd_{\kmax+1,x}$ and $\Ecp_{\kmax,x}$ has failed. Putting these arguments together carefully, we will obtain the desired estimate.

Let us set up this argument in detail. We define
\[\kmax:=(k-1)\wedge\max\left\{k'\in\{k_0(x)-1\ldots, k_\infty'\}\colon \phi\in\Ero_{\le k',x}\cap \Ebd_{\le k',x}\cap \Ecp_{\le k'-1,x} \right\}.\]
Note that the condition $\phi\in\Ero_{\le k',x}\cap \Ebd_{\le k',x}\cap \Ecp_{\le k'-1,x}$ is vacuous for $k'=k_0(x)-1$, and so $\kmax\ge k_0(x)-1$ is well-defined. 
We also have, for any $x\in A_{(\log n)^{5/2}}$,
\[\tilde\PP^{Q_N,0}\left(\kmax=k_0(x)-1\right)\le \tilde\PP^{Q_N,0}\left((\Ero_{k_0(x),x}\cap \Ebd_{k_0(x),x})^\com\right)\le C\exp(-c\log r_{k_0(x)})^3)\le C\exp(-cn^3),\]
because $k_0(x)\le\frac n2$.

We can now decompose
\begin{equation}\label{e:upper_barrier_late2}
\begin{split}
&\tilde\PP^{Q_N,0}\left(\exists x'\in Q_{r_k/(n-k)^2}(x)\colon \phi\in\Eup_{\le k,x} \cap (\Eup_{\le k+1,x'})^\com\right)\\
&\le \tilde\PP^{Q_N,0}(\kmax=k_0(x)-1)\\
&\qquad\qquad +\sum_{k'=k_0(x)-1}^{k-1}\tilde\PP^{Q_N,0}\left(\exists x'\in Q_{r_k/(n-k)^2}(x)\colon \phi\in\Eup_{\le k,x} \cap (\Eup_{\le k+1,x'})^\com ,\kmax=k'\right)\\
&\le C\exp(-cn^3)+\sum_{k'=k_0(x)}^{k-1}\tilde\PP^{Q_N,0}\left(\exists x'\in Q_{r_k/(n-k)^2}(x)\colon \phi\in\Eup_{\le k,x} \cap (\Eup_{\le k+1,x'})^\com ,\kmax=k'\right).
\end{split}
\end{equation}

We will bound each of the summands separately. The dominant contribution will come from the term with $k'=k-1$, and at some point later on we will distinguish the cases of whether
 $k'<k-1$ or $k'=k-1$. For now let us consider an arbitrary $k'$, and simplify the probability in question.

On the event $\{\kmax=k'\}\subset\Eup_{\le k',x} \cap \Ero_{\le k',x}\cap \Ecp_{\le k'-1,x}$ we can use the representation \eqref{e:telescope2} to see that
\[\Sav_{j,-}(x,\phi)=\Sav_{k_0(x),+}(x,\phi)+\Ind_{k_0(x)}(x,\phi)+\Ind_{k_0(x)+1}(x,\phi)+\ldots + \Ind_{j-1}(x,\phi)+O(1)\]
for all $j\le k'$. Moreover, suppose that we know that $\Sav_{k_0(x),+}(x,\phi)\in[s,s+1]$ and $\Sav_{k',-}(x,\phi)\in\Bup(k')+\Gamma-[t,t+1]$ for some $s\in\Z$, $t\in\N$. Then $\Eup_{\le k,x}\subset \Eup_{\le k',x}\subset \widetilde{\Eup_{\le k',x}}$, where
\[\widetilde{\Eup_{\le k',x}}:=\left\{\phi\colon s+\Ind_{k_0(x)}(x,\phi)+\ldots + \Ind_{j-1}(x,\phi)\le\Bup(j)+\Gamma+C\ \forall j\in\{k_0(x),\ldots,k'\}\right\}.\]
Similarly, we see that if there is some $x'$ for which $(\Eup_{\le k+1,x'})^\com$ occurs, then
\[\max_{x'\in Q_{r_k/(n-k)^2}(x)}\biggl( \Sav_{k,+}(x',\phi)-\Sav_{k',-}(x,\phi) \biggr) \ge \Bup(k+1)-\Bup(k')+t-C.\]
So we can rewrite the summands in \eqref{e:upper_barrier_late2}, using
\[W_{x,x'}:=\Sav_{k,+}(x',\phi)-\Sav_{k',-}(x,\phi),\]
as
\begin{equation}\label{e:upper_barrier_late3}
\begin{split}
&\tilde\PP^{Q_N,0}\left(\exists x'\in Q_{r_k/(n-k)^2}(x)\colon \phi\in\Eup_{\le k,x} \cap (\Eup_{\le k+1,x'})^\com ,\kmax=k'\right)\\
&\le \sum_{s=-\infty}^{\Bup(k_0(x))}\sum_{t=0}^\infty\tilde\PP^{Q_N,0}\Bigg(\Sav_{k_0(x),+}(x,\phi)\in[s,s+1],\Sav_{k',-}(x,\phi)\in\Bup(k')+\Gamma-[t,t+1],\widetilde{\Eup_{\le k',x}},\\
&\qquad\max_{x'\in Q_{r_k/(n-k)^2}(x)} W_{x,x'}
\ge \Bup(k+1)-\Bup(k')+t-C,\kmax=k'\Bigg)\\
&\le\sum_{s=-\infty}^{\Bup(k_0(x))}\sum_{t=0}^\infty\tilde\E^{Q_N,0}\Bigg(\tilde\E^{Q_N,0}\Bigg(\\
&\qquad\tilde\PP^{Q_N,0}\Bigg(\max_{x'\in Q_{r_k/(n-k)^2}(x)} W_{x,x'}
\ge \Bup(k+1)-\Bup(k')+t-C,\\
&\qquad\qquad\Sav_{k',-}(x,\phi)\in\Bup(k')+\Gamma-[t,t+1], \kmax=k'\Bigg|F_{Q_{r_{k'}}^\com(x)}\Bigg)\I_{\Ero_{k',x}}\I_{\widetilde{\Eup_{\le k',x}}}\\
&\qquad\qquad\I_{s+\Ind_{k_0(x)}(x,\phi)+\ldots + \Ind_{k'-1}(x,\phi)\in \Bup(k')-[t-C,t+C]}
\Bigg| F_{Q_{r_{k_0(x)}}^\com(x)}\Bigg)\I_{\Ero_{k_0(x),x}}\I_{\Sav_{k_0(x),+}(x,\phi)\in[s,s+1]}\Bigg).
\end{split}
\end{equation}
We will compute this tower of conditional expectations from the inside out. 

\emph{Step 2: Estimate of the contribution between $\kmax$ and $k$}\\
Regarding the innermost conditional probability, we claim that
\begin{equation}\label{e:upper_barrier_late4}
\begin{split}
  &\tilde\PP^{Q_N,0}\Bigg(\max_{x'\in Q_{r_k/(n-k)^2}(x)} W_{x,x'}
  \ge \Bup(k+1)-\Bup(k')+t-C,\\
&\qquad\qquad\Sav_{k',-}(x,\phi)\in\Bup(k')+\Gamma-[t,t+1], \kmax=k'\Bigg|F_{Q_{r_{k'}}^\com(x)}\Bigg)\I_{\Ero_{k',x}}\\
&\le C\exp\left(\frac{\bar\lambda^2\g(k-k')}{2}-\bar\lambda(\Bup(k+1)-\Bup(k')-t)-\frac{t^2}{C(k-k')^3}-c(k-k')\right),
\end{split}
\end{equation}
where
\begin{equation}\label{e:def_lambda}
\bar\lambda:=\frac{\sqrt{\g}(m_N-2\Delta)}{\g k_\infty}=\frac{m_N-2\Delta}{\sqrt{\g} (n-\Delta)}=\frac{2}{\sqrt{\g}}-\frac{3\log n}{4\sqrt{\g}(n-\Delta)}.
\end{equation}

To show this claim, we distinguish several cases. The simplest case is if $k'=k-1$. In that case we can just drop the condition that $\kmax=k'$, and estimate that
\begin{align*}
  &\tilde\PP^{Q_N,0}\Bigg(\max_{x'\in Q_{r_k/(n-k)^2}(x)} W_{x,x'}
  \ge \Bup(k+1)-\Bup(k')+t-C,\\
&\qquad\qquad\Sav_{k',-}(x,\phi)\in\Bup(k')+\Gamma-[t,t+1], \kmax=k'\Bigg|\tilde\F_{Q_{r_{k'}}^\com(x)}\Bigg)\I_{\Ero_{k',x}}\\
&\le \tilde\PP^{Q_N,0}\Bigg(\max_{x'\in Q_{r_k/(n-k)^2}(x)} W_{x,x'}
\ge \Bup(k+1)-\Bup(k-1)+t-C\Bigg)\I_{\Ero_{k',x}}\\
&\le \tilde\PP^{Q_N,0}\Bigg(\max_{x'\in Q_{r_k/(n-k)^2}(x)} W_{x,x'}
\ge t-C\Bigg)\I_{\Ero_{k',x}},
\end{align*}
so that \eqref{e:upper_barrier_late4} immediately follows from Lemma \ref{l:fluctuations_cond_tail_bound}. If, on the other hand, $k'<k-1$, then $\kmax=k'$ implies that at least one of $(\Ero_{k'+1,x})^\com$, $(\Ebd_{k'+1,x})^\com$ or $(\Ecp_{k',x})^\com$ occurs. We will distinguish the two subcases that $(\Ero_{k'+1,x})^\com\cup(\Ebd_{k'+1,x})^\com$ or that $\Ero_{k'+1,x}\cap\Ebd_{k'+1,x}\cap (\Ecp_{k',x})^\com$ occur.

Let us begin with the first subcase, that $(\Ero_{k'+1,x})^\com\cup(\Ebd_{k'+1,x})^\com$ occurs. Then, using Cauchy-Schwarz and Lemma \ref{l:fluctuations_cond_tail_bound} (with $\lambda=2\bar\lambda$ and $\eps=1$, say), we can estimate that
\begin{equation}\label{e:upper_barrier_late5}
\begin{split}
  &\tilde\PP^{Q_N,0}\Bigg(\max_{x'\in Q_{r_k/(n-k)^2}(x)} W_{x,x'}
  \ge \Bup(k+1)-\Bup(k')+t-C,\\
&\qquad\qquad\Sav_{k',-}(x,\phi)\in\Bup(k')+\Gamma-[t,t+1], (\Ero_{k'+1,x})^\com\cup(\Ebd_{k'+1,x})^\com\Bigg|\tilde\F_{Q_{r_{k'}}^\com(x)}\Bigg)\I_{\Ero_{k',x}}\\
&\le \left(\tilde\PP^{Q_N,0}\Bigg(\max_{x'\in Q_{r_k/(n-k)^2}(x)} W_{x,x'}
\ge \Bup(k+1)-\Bup(k')+t-C\Bigg)\right)^{1/2}\\
&\qquad\times\left(\tilde\PP^{Q_N,0}\left((\Ero_{k'+1,x})^\com\cup(\Ebd_{k'+1,x})^\com\right)\right)^{1/2}\I_{\Ero_{k',x}}\\
&\le C\left(\exp\left(\frac{((2\bar\lambda)^2+1)\g(k-k')}{2}-2\bar\lambda(\Bup(k+1)-\Bup(k')-t)-\frac{t^2}{C(k-k')^3}\right)\right)^{1/2}\\
&\qquad \qquad \times \left(\exp\left(-c(\log r_{k'+1})^3\right)\right)^{1/2}\\
&\le C\exp\left(\left(\bar\lambda^2+\frac14\right)\g(k-k')-\bar\lambda(\Bup(k+1)-\Bup(k')-t)-\frac{t^2}{C(k-k')^3}-c(n-k')^3\right)\\
&\le C\exp\left(\frac{\bar\lambda^2\g(k-k')}{2}-\bar\lambda(\Bup(k+1)-\Bup(k')-t)-\frac{t^2}{C(k-k')^3}-c(n-k')^3\right),
\end{split}
\end{equation}
as the term $(n-k')^3\ge (k-k')^3$ easily beats $C(k-k')$.

In the second subcase that $\Ero_{k'+1,x}\cap\Ebd_{k'+1,x}\cap (\Ecp_{k',x})^\com$ occurs, it would be too wasteful to use Cauchy-Schwary from scale $k'$ onwards as in the first subcase. Instead we will use Cauchy-Schwarz only on scale $k'$, while returning to the sharp asymptotics from Lemma \ref{l:fluctuations_cond_tail_bound} from scale $k'+1$ onwards. To do so, observe that on $\Ebd_{k'+1,x}$ we have $\Sav_{k'+1,-}(x,\phi)=\Sav_{k'+1,+}(x,\phi)+O(1)$. So, if we know that $\Sav_{k'+1,+}(x,\phi)\in \Bup(k'+1)+\Gamma-[\tilde t,\tilde t+1]$, then necessarily $\Sav_{k'+1,-}(x,\phi)\le \Bup(k'+1)+\Gamma-\tilde t+C$.
We can now condition on $\tilde\F_{Q_{r_{k'+1}}^\com(x)}$ and see that
\begin{align*}
  &\tilde\PP^{Q_N,0}\Bigg(\max_{x'\in Q_{r_k/(n-k)^2}(x)} W_{x,x'}
  \ge \Bup(k+1)-\Bup(k')+t-C,\\
&\qquad\qquad\Sav_{k',-}(x,\phi)\in\Bup(k')+\Gamma-[t,t+1], \Ero_{k'+1,x}\cap\Ebd_{k'+1,x}\cap (\Ecp_{k',x})^\com\Bigg|\tilde\F_{Q_{r_{k'}}^\com(x)}\Bigg)\I_{\Ero_{k',x}}\\
&\le \sum_{\tilde t=-\infty}^{\infty}\tilde\PP^{Q_N,0}\Bigg(\max_{x'\in Q_{r_k/(n-k)^2}(x)} W_{x,x'}
\ge \Bup(k+1)-\Bup(k'+1)+\tilde t-C,\\
&\qquad\qquad\Sav_{k',-}(x,\phi)\in\Bup(k')+\Gamma-[t,t+1], \Sav_{k'+1,+}(x,\phi)\in \Bup(k'+1)+\Gamma-[\tilde t,\tilde t+1],\\
&\qquad\qquad\Ero_{k'+1,x}\cap\Ebd_{k'+1,x}\cap (\Ecp_{k',x})^\com\Bigg|\tilde\F_{Q_{r_{k'}}^\com(x)}\Bigg)\I_{\Ero_{k',x}}\\
&\le \sum_{\tilde t=-\infty}^{\infty}\tilde\E^{Q_N,0}\Bigg(\tilde\PP^{Q_N,0}\Bigg(\max_{x'\in Q_{r_k/(n-k)^2}(x)}W_{x,x'}
\ge \Bup(k+1)-\Bup(k'+1)+\tilde t-C\Bigg|\tilde\F_{Q_{r_{k'+1}^\com}(x)}\Bigg)\\
&\qquad\qquad\I_{\Ero_{k'+1,x}}\I_{|\Sav_{k',-}(x,\phi)-\Sav_{k'+1,-}(x,\phi)|\ge \tilde t-t-C}\I_{(\Ecp_{k',x})^\com}\Bigg|\tilde\F_{Q_{r_{k'}}^\com(x)}\Bigg)\I_{\Ero_{k',x}}\\
&\le \sum_{\tilde t=-\infty}^{\infty}C_{\eps}\exp\left(\frac{(\bar\lambda^2+\eps)\g(k-k'-1)}{2}-\bar\lambda(\Bup(k+1)-\Bup(k'+1)-\tilde t)-\frac{\tilde t^2}{C(k-k'-1)^3}\right)\\
&\qquad\qquad\tilde\PP^{Q_N,0}\Bigg(|\Sav_{k',-}(x,\phi)-\Sav_{k'+1,-}(x,\phi)|\ge \tilde t-t-C,(\Ecp_{k',x})^\com\Bigg|\tilde\F_{Q_{r_{k'}}^\com(x)}\Bigg)\I_{\Ero_{k',x}},
\end{align*}
where we used Lemma \ref{l:fluctuations_cond_tail_bound} (for some $
\eps$ to be fixed shortly) in the last step. The remaining probability is easily bounded using Brascamp-Lieb and Cauchy-Schwarz, and we obtain
\begin{align*}
&\tilde\PP^{Q_N,0}\Bigg(\max_{x'\in Q_{r_k/(n-k)^2}(x)} W_{x,x'}\ge \Bup(k+1)-\Bup(k')+t-C,\\
&\qquad\qquad\Sav_{k',-}(x,\phi)\in\Bup(k')+\Gamma-[t,t+1], \Ero_{k'+1,x}\cap\Ebd_{k'+1,x}\cap (\Ecp_{k',x})^\com\Bigg|\tilde\F_{Q_{r_{k'}}^\com(x)}\Bigg)\I_{\Ero_{k',x}}\\
&\le\sum_{\tilde t=-\infty}^{\infty}C_{\eps}\exp\left(\frac{(\bar\lambda^2+\eps)\g(k-k'-1)}{2}-\bar\lambda(\Bup(k+1)-\Bup(k'+1)-\tilde t)-\frac{\tilde t^2}{C(k-k'-1)^3}\right)\\
&\qquad\qquad\left(\exp\left(-\frac{(\tilde t-t-C)^2}{C}\right)\right)^{1/2}\left(\frac{1}{r_{k'}^\omega}\right)^{1/2}\\
&\le C_{\eps}\exp\left(\frac{(\bar\lambda^2+\eps)\g(k-k')}{2}-\bar\lambda(\Bup(k+1)-\Bup(k')-t)-\frac{t^2}{C(k-k')^3}-\frac{\omega(n-k')}{2}\right),
\end{align*}
where in the last step we used that because of the factor $\exp\left(-\frac{(\tilde t-t-C)^2}{C}\right)$ we essentially only need to consider $\tilde t=t+O(1)$. We now choose $\eps=\frac{\omega}{4\g}$ and conclude that
\begin{equation}\label{e:upper_barrier_late6}
\begin{split}
&\tilde\PP^{Q_N,0}\Bigg(\max_{x'\in Q_{r_k/(n-k)^2}(x)} W_{x,x'}\ge \Bup(k+1)-\Bup(k')+t-C,\\
&\qquad\qquad\Sav_{k',-}(x,\phi)\in\Bup(k')+\Gamma-[t,t+1], \Ero_{k'+1,x}\cap\Ebd_{k'+1,x}\cap (\Ecp_{k',x})^\com\Bigg|\tilde\F_{Q_{r_{k'}}^\com(x)}\Bigg)\I_{\Ero_{k',x}}\\
&\le C\exp\left(\frac{\bar\lambda^2\g(k-k')}{2}-\bar\lambda(\Bup(k+1)-\Bup(k')-t)-\frac{t^2}{C(k-k')^3}-\frac{\omega(n-k')}{4}\right).
\end{split}
\end{equation}
If we now combine \eqref{e:upper_barrier_late5} and \eqref{e:upper_barrier_late6}, we obtain that
\begin{align*}
&\tilde\PP^{Q_N,0}\Bigg(\max_{x'\in Q_{r_k/(n-k)^2}(x)} \Sav_{k,+}(x',\phi)-\Sav_{k',-}(x,\phi)\ge \Bup(k+1)-\Bup(k')+t-C,\\
&\qquad\qquad\Sav_{k',-}(x,\phi)\in\Bup(k')+\Gamma-[t,t+1], \kmax=k'\Bigg|F_{Q_{r_{k'}^\com}(x)}\Bigg)\I_{\Ero_{k',x}}\\
&\le C\exp\left(\frac{\bar\lambda^2\g(k-k')}{2}-\bar\lambda(\Bup(k+1)-\Bup(k')-t)-\frac{t^2}{C(k-k')^3}\right)\\
& \qquad \times \left(\exp\left(-c(n-k')^3\right)+\exp\left(-\frac{\omega(n-k')}{4}\right)\right),
\end{align*}
which implies \eqref{e:upper_barrier_late4} in the remaining case that $k'<k-1$.

Thus we now have established \eqref{e:upper_barrier_late4} in all cases. We can now insert this estimate into \eqref{e:upper_barrier_late3} and arrive at
\begin{equation}\label{e:upper_barrier_late7}
\begin{split}
&\tilde\PP^{Q_N,0}\left(\exists x'\in Q_{r_k/(n-k)^2}(x)\colon \phi\in\Eup_{\le k,x} \cap (\Eup_{\le k+1,x'})^\com ,\kmax=k'\right)\\
&\le C\sum_{s=-\infty}^{\Bup(k_0(x))}\sum_{t=0}^\infty \exp\left(\frac{\bar\lambda^2\g(k-k')}{2}-\bar\lambda(\Bup(k+1)-\Bup(k')-t)-\frac{t^2}{C(k-k')^3}-c(k-k')\right)\\
&\tilde\E^{Q_N,0}\Bigg(\tilde\PP^{Q_N,0}\Bigg(\widetilde{\Eup_{\le k',x}},s+\Ind_{k_0(x)}(x,\phi)+\ldots + \Ind_{k'-1}(x,\phi)\in \Bup(k')-[t-C,t+C]
\Bigg| F_{Q_{r_{k_0(x)}}^\com(x)}\Bigg)\\
&\qquad\qquad\I_{\Ero_{k_0(x),x}}\I_{\Sav_{k_0(x),+}(x,\phi)\in[s,s+1]}\Bigg).
\end{split}
\end{equation}

\emph{Step 3: Estimate of the contribution between $k_0$ and $\kmax$}\\
In order to deal with the new innermost probability, we will use the barrier estimate from Lemma \ref{l:upper_barrier}.
We consider the tilted measure
\[
\frac{d\Q}{d\tilde\PP^{Q_N,0}\left(\cdot\middle|\tilde\F_{Q_{r_{k_0(x)}}^\com(x)}\right)} :=\frac{ \exp\left(\bar\lambda \left(\sum_{k=k_0}^{k'-1} \Ind_k (x,\phi)\right) \right)}{\tilde\E^{Q_N,0}\left(\exp\left(\bar\lambda \left(\sum_{k=k_0}^{k'-1} \Ind_k (x,\phi)\right) \right) \middle|\tilde\F_{Q_{r_{k_0(x)}}^\com(x)}\right)},\]
with $\bar\lambda$ as in \eqref{e:def_lambda}.
As we verify in Lemma \ref{l:tilting}, under $\Q$ the random variables $\Ind_k(x)-\g\bar\lambda$ satisfy the assumptions of Lemma \ref{l:upper_barrier}. Using that the interval $[t-C,t+C]$ is the union of $2C$ intervals of length $1$, we thus obtain from  from Lemma \ref{l:upper_barrier}  that
\begin{align*}
&\Q\left(\widetilde{\Eup_{\le k',x}},s+\Ind_{k_0(x)}(x,\phi)+\ldots + \Ind_{k'-1}(x,\phi)\in \Bup(k')-[t-C,t+C]\right)\\
&\le \sum_{i=-C}^{C-1}\Q\left(\widetilde{\Eup_{\le k',x}},s+\Ind_{k_0(x)}(x,\phi)+\ldots + \Ind_{k'-1}(x,\phi)\in \Bup(k')-[t+i,t+i+1]\right)\\
&\le C\frac{\left(1+\Bup(k_0(x))-s\right)\left(1+t+\Bup(k')-\g\bar\lambda(k'-k_0(x))\right)}{(k'-k_0(x))^{3/2}}.
\end{align*}
(The terms $\Bup(k_0(x))=\g\bar\lambda k_0(x)+k_0(x)^{2/5}$ and $\Bup(k')-\g\bar\lambda(k'-k_0(x))=\g\bar\lambda k_0(x)+(k'\vee (k_\infty-k'))^{2/5}$   appear because of the fact that our random walk starts at $k_0(x)$ instead of 0, which effectively shifts the barrier upwards by $\bar\lambda k_0(x)$.)
Now
\begin{equation}\label{e:upper_barrier_late8}
\begin{split}
&\tilde\PP^{Q_N,0}\left(\widetilde{\Eup_{\le k',x}},s+\Ind_{k_0(x)}(x,\phi)+\ldots + \Ind_{k'-1}(x,\phi)\in \Bup(k')-[t-C,t+C]\right)\\
&\le\E\Bigg(\I_{\widetilde{\Eup_{\le k',x}}}\I_{s+\Ind_{k_0(x)}(x,\phi)+\ldots + \Ind_{k'-1}(x,\phi)\in \Bup(k')+\Gamma-[t-C,t+C]}\\
&\qquad\qquad\exp\left(\bar\lambda \left(\sum_{k=k_0(x)}^{k'-1} \Ind_k (x,\phi)-\left(\Bup(k')+\Gamma-s-t-C\right)\right)\right)\Bigg|\tilde\F_{Q_{r_{k_0(x)}}^\com(x)}\Bigg)\\
\end{split} \end{equation}
which, undoing the tilt, is in turn bounded by
\begin{equation*}
\begin{split}
&\le \exp\left(-\bar\lambda\left(\Bup(k')+\Gamma-s-t-C\right)\right)\tilde\E^{Q_N,0}\left(\exp\left(\bar\lambda \left(\sum_{k=k_0}^{k'-1} \Ind_k (x,\phi)\right) \right) \middle|\tilde\F_{Q_{r_{k_0(x)}}^\com(x)}\right)\\
&\qquad\qquad\Q\left(\widetilde{\Eup_{\le k,x}},s+\Ind_{k_0(x)}(x,\phi)+\ldots + \Ind_{k-2}(x,\phi)\in \Bup(k-1)+\Gamma-[t-C,t+C]\right)\\
&\le C\frac{\left(1+\Bup(k_0(x))-s\right)\left(1+t+\Bup(k')-\g\bar\lambda(k'-k_0(x))\right)}{(k'-k_0(x))^{3/2}}\\
&\qquad\exp\left(\frac{\bar\lambda^2\g(k'-k_0(x))}{2}-\bar\lambda(\left(\Bup(k')+\Gamma-s-t-C\right)\right)\\
&\le C\frac{\left(1+\Bup(k_0(x))-s\right)\left(1+t+\Bup(k')-\g\bar\lambda(k'-k_0(x))\right)}{(k'-k_0(x))^{3/2}}\\
&\qquad\exp\left(\frac{\bar\lambda^2\g(k'-k_0(x))}{2}-\bar\lambda(\left(\Bup(k')+\Gamma-s-t-C\right)\right).
\end{split}
\end{equation*}

If we combine \eqref{e:upper_barrier_late8} with \eqref{e:upper_barrier_late7}, all that is left to estimate is the outer expectation. For it we use
\[\tilde\PP^{Q_N,0}\left(\Sav_{k_0(x),+}(x,\phi)\in[s,s+1]\right)\le \exp\left(-\frac{s^2}{C}\right),\]
as follows from Lemma \ref{l:sharpmomentbound}. 

\emph{Step 4: Completion of the proof}\\
Putting things together, we arrive at
\begin{equation}\label{e:upper_barrier_late9}
\begin{split}
&\tilde\PP^{Q_N,0}\left(\exists x'\in Q_{r_k/(n-k)^2}(x)\colon \phi\in\Eup_{\le k,x} \cap (\Eup_{\le k+1,x'})^\com ,\kmax=k'\right)\\
&\le C\sum_{s=-\infty}^{\Bup(k_0(x))+1}\sum_{t=0}^\infty \exp\left(-\frac{s^2}{C}\right) \frac{\left(1+\Bup(k_0(x))-s\right)\left(1+t+\Bup(k')-\g\bar\lambda(k'-k_0(x))\right)}{(k'-k_0(x))^{3/2}}\\
&\qquad\exp\left(\frac{\bar\lambda^2\g(k'-k_0(x))}{2}-\bar\lambda(\left(\Bup(k')+\Gamma-s-t-C\right)\right)\\
&\qquad \exp\left(\frac{\bar\lambda^2\g(k-k')}{2}-\bar\lambda(\Bup(k+1)-\Bup(k')-t)-\frac{t^2}{C(k-k')^3}-c(k-k')\right)\\
&\le C\sum_{s=-\infty}^{\Bup(k_0(x))+1}\sum_{t=0}^\infty \frac{\left(1+\g\bar\lambda k_0(x)+k_0(x)^{2/5}-s\right)\left(1+t+\g\bar\lambda k_0(x)+(n-k')^{2/5}\right)}{(k'-k_0(x))^{3/2}}\\
&\qquad\exp\left(\frac{\bar\lambda^2\g(k-k_0(x))}{2}-\bar\lambda\Bup(k+1)-\bar\lambda(\Gamma-s-t)-\frac{s^2}{C}-\frac{t^2}{C(k-k')^3}-c(k-k')\right).
\end{split}
\end{equation}
Consider the factor
\[\frac{\left(1+\g\bar\lambda k_0(x)+k_0(x)^{2/5}-s\right)\left(1+t+\g\bar\lambda k_0(x)+(k_\infty-k')^{2/5}\right)}{(k'-k_0(x))^{3/2}}\exp\left(-\frac{\bar\lambda^2\g k_0(x)}{2}\right).\]
We can expand the numerator in powers of $k_0(x)$. As a function of $k_0(x)$, the exponentially decaying term dominates all polynomial terms $k_0(x)$. Thus, effectively, we can replace $k_0(x)$ by $0$ at the cost of possibly increasing $C$. This leads to the bound
\begin{align*}
&\frac{\left(1+\g\bar\lambda k_0(x)+k_0(x)^{2/5}-s\right)\left(1+t+\g\bar\lambda k_0(x)+(n-k')^{2/5}\right)}{(k'-k_0(x))^{3/2}}\exp\left(-\frac{\bar\lambda^2\g k_0(x)}{2}\right)\\
&\le C\frac{((C-s)\vee 1)\left(1+t+(n-k')^{2/5}\right)}{k'^{3/2}}.
\end{align*}
Inserting this bound into \eqref{e:upper_barrier_late9} and recalling the definition of $\bar\lambda$ from \eqref{e:def_lambda}, we obtain
\begin{align*}
&\tilde\PP^{Q_N,0}\left(\exists x'\in Q_{r_k/(n-k)^2}(x)\colon \phi\in\Eup_{\le k,x} \cap (\Eup_{\le k+1,x'})^\com ,\kmax=k'\right)\\
&\le C\sum_{s=-\infty}^{\Bup(k_0(x))+1}\sum_{t=0}^\infty\frac{((C-s)\vee 1)\left(1+t+(n-k')^{2/5}\right)}{k'^{3/2}}\\
&\qquad\exp\left(\frac{\bar\lambda^2\g k}{2}-\bar\lambda\Bup(k+1)-\bar\lambda(\Gamma-s-t)-\frac{s^2}{C}-\frac{t^2}{C(k-k')^3}-c(k-k')\right)\\
&\le C\sum_{s=-\infty}^\infty\sum_{t=0}^\infty(\frac{((C-s)\vee 1)\left(1+t+(n-k')^{2/5}\right)}{k'^{3/2}}\\
&\qquad\times \exp\left(-2k+\frac{3k\log n}{2n}-\frac{2}{\sqrt{\g}}(n-k-1)^{2/5}-c\Gamma-O\left(\frac{k(\log n)^2}{n^2}\right)-c(k-k')\right)\\
&\qquad \times \exp\left(\frac{2}{\sqrt{\g}}(s+t)-\frac{s^2}{C}-\frac{t^2}{C(k-k')^3}\right).
\end{align*}
As $a\mapsto\frac{\log a}{a}$ is a decreasing function of $a$ for $a\ge\e$, we have $\frac{k\log n}{n}\le \log k$, if $n$ (and thus $k$) is sufficiently large. This allows us to bound
\[
\frac{1}{k'^{3/2}}\exp\left(\frac{3k\log n}{2n}-c(k-k')\right)\le\exp\left(\frac32(\log k-\log k')-c(k-k')\right)\le C\exp\left(-c'(k-k')\right),
\]
so that we obtain the estimate
\begin{align*}
&\tilde\PP^{Q_N,0}\left(\exists x'\in Q_{r_k/(n-k)^2}(x)\colon \phi\in\Eup_{\le k,x} \cap (\Eup_{\le k+1,x'})^\com ,\kmax=k'\right)\\
&\le C\sum_{s=-\infty}^\infty\sum_{t=0}^\infty\left((C-s)\vee 1\right)\left(1+t+(n-k')^{2/5}\right)\\
&\qquad\exp\left(-2k-\frac{2}{\sqrt{\g}}(n-k-1)^{2/5}-c\Gamma-c'(k-k')+\frac{2}{\sqrt{\g}}(s+t)-\frac{s^2}{C}-\frac{t^2}{C(k-k')^3}\right).
\end{align*}
Next, note that the sums over $s$ and $t$ are dominated by the first $O(1)$ and $O((k-k')^{3/2})$ terms, respectively, so that actually
\begin{equation}\label{e:upper_barrier_late10}
\begin{split}
&\tilde\PP^{Q_N,0}\left(\exists x'\in Q_{r_k/(n-k)^2}(x)\colon \phi\in\Eup_{\le k,x} \cap (\Eup_{\le k+1,x'})^\com ,\kmax=k'\right)\\
&\le C\left((k-k')^{3/2}+(n-k')^{2/5}\right)\exp\left(-2k-\frac{2}{\sqrt{\g}}(n-k-1)^{2/5}-c\Gamma-c'(k-k')\right)\\
&\le C\exp\left(-2k-c''(n-k-1)^{2/5}-c''(k-k')-c\Gamma\right).
\end{split}
\end{equation}
The estimate \eqref{e:upper_barrier_late10} allows to bound the summands of the right-hand side of \eqref{e:upper_barrier_late2}. We obtain
\begin{equation}\label{e:upper_barrier_late11}
\begin{split}
&\tilde\PP^{Q_N,0}\left(\exists x'\in Q_{r_k/(n-k)^2}(x)\colon \phi\in\Eup_{\le k,x} \cap (\Eup_{\le k+1,x'})^\com\right)\\
&\le C\exp(-cn^3)+C\sum_{k'=k_0(x)}^{k-1}\exp\left(-2k-c(n-k-1)^{2/5}-c\Gamma-c(k-k')\right)\\
&\le C\exp\left(-2k-c(n-k-1)^{2/5}-c\Gamma\right),
\end{split}
\end{equation}
for sufficiently large $n$ (depending on $\Gamma$).

Now from \eqref{e:upper_barrier_late1} and \eqref{e:upper_barrier_late11} we directly obtain
\begin{align*}
&\tilde\PP^{Q_N,0}\left(\bigcap_{x\in A_{(\log n)^{5/2}}}\Eup_{\le k,x} \cap \bigcup_{x\in A_{(\log n)^{5/2}}}(\Eup_{\le k+1,x})^\com\right)\\
&\le C\left(\frac{N(n-k)}{r_k}\right)^2\exp\left(-2k-c(n-k-1)^{2/5}-c\Gamma\right)\\
&\le C\exp\left(2\log (n-k)-c(n-k-1)^{2/5}-c\Gamma\right),
\end{align*}
and the lemma follows once we note that $2\log (n-k)$ can be absorbed into $c(n-k-1)^{2/5}$.
\end{proof}
For the very last step of our inductive argument, we need the following variant of Lemma \ref{l:upper_barrier_late}.
\begin{lemma}\label{l:upper_barrier_laststep}
Let $A_{(\log n)^{5/2}}$ be as in \eqref{e:def_Ak}. We have
\begin{equation}\label{e:upper_barrier_laststep}
\tilde\PP^{Q_N,0}\left(\bigcap_{x\in A_{(\log n)^{5/2}}}\Eup_{\le k_{\infty'},x} \cap \{\phi\colon \phi(x)\ge m_N+\Gamma\}\right)\le C\exp\left(-c\Gamma\right).
\end{equation}
\end{lemma}
\begin{proof}
The proof is similar to the proof of Lemma \ref{l:upper_barrier_late}, with the main difference being that we use the tail bounds following from Lemma \ref{l:sharpmomentbound} instead of Lemma \ref{l:fluctuations_cond_tail_bound}.

Let us give some details on this. We define
\[\kmax:=\max\left\{k'\in\{k_0(x)-1\ldots, k_\infty'\}\colon \phi\in\Ero_{\le k',x}\cap \Ebd_{\le k',x}\cap \Ecp_{\le k'-1,x} \right\}\]
and write
\begin{equation}\label{e:upper_barrier_laststep1}
\begin{split}
&\tilde\PP^{Q_N,0}\left(\bigcap_{x\in A_{(\log n)^{5/2}}}\Eup_{\le k_{\infty'},x} \cap \{\phi\colon \phi(x)\ge m_N+\Gamma\}\right)\\
&\le \PP(\kmax=k_0(x)-1)+\sum_{k'=k_0(x)-1}^{k-1}\tilde\PP^{Q_N,0}\left(\phi(x)\ge m_N+\Gamma , \Eup_{\le k_\infty',x} ,\kmax=k'\right).
\end{split}
\end{equation}
For each individual summand in \eqref{e:upper_barrier_laststep1}, we can use a decomposition analogous to \eqref{e:upper_barrier_late3} and write
\begin{equation}\label{e:upper_barrier_laststep2}
\begin{split}
&\tilde\PP^{Q_N,0}\left(\phi(x)\ge m_N+\Gamma , \Eup_{\le k_\infty',x} ,\kmax=k'\right)\\
&\le \sum_{s=-\infty}^{\Bup(k_0(x))}\sum_{t=0}^\infty\tilde\PP^{Q_N,0}\Bigg(\Sav_{k_0(x),+}(x,\phi)\in[s,s+1],\Sav_{k',-}(x,\phi)\in\Bup(k')+\Gamma-[t,t+1],\\
&\qquad\qquad \widetilde{\Eup_{\le k_\infty',x}},\phi(x)-\Sav_{k',-}(x,\phi)\ge m_N-\Bup(k')-\Gamma+t-C\Bigg)\\
&\le\sum_{s=-\infty}^{\Bup(k_0(x))}\sum_{t=0}^\infty\tilde\E^{Q_N,0}\Bigg(\tilde\E^{Q_N,0}\Bigg(\tilde\PP^{Q_N,0}\Bigg(\phi(x)-\Bup(k')+t+1\ge m_N+\Gamma-C,\\
&\qquad\qquad\Sav_{k',-}(x,\phi)\in\Bup(k')+\Gamma-[t,t+1]\Bigg|\F_{Q_{r_{k'}^\com}(x)}\Bigg)\I_{\Ero_{k',x}}\I_{\widetilde{\Eup_{\le k',x}}}\\
&\qquad\qquad\I_{s+\Ind_{k_0(x)}(x)+\ldots + \Ind_{k'-1}(x)\in \Bup(k')-[t-C,t+C]}
\Bigg|\F_{Q_{r_{k_0(x)}^\com}(x)}\Bigg)\I_{\Ero_{k_0(x),x}}\I_{\Sav_{k_0(x),+}(x,\phi)\in[s,s+1]}\Bigg).
\end{split}
\end{equation}
(This is the analogue of estimate \eqref{e:upper_barrier_late3}.) Again we estimate this tower of conditional expectations from the inside out. For that purpose, we need a replacement for \eqref{e:upper_barrier_late4}. However, this is immediate from the exponential moment estimate in Lemma \ref{l:sharpmomentbound} (valid for all $|\lambda|\le\lambda_*$) and the Brascamp-Lieb inequality (valid for all $\lambda$). Just like we deduced \eqref{e:fluctuations_cond_tail_bound} from Lemma \ref{l:fluctuations_cond}, we deduce from these results that
\begin{equation}\label{e:upper_barrier_laststep3}
\tilde\PP^{Q_N,0}\left(\left|\phi(x)-\Sav_{k,-}(x,\phi)\right|\ge t\middle|\tilde\F_{Q_{r_{k}}(x)^\com}\right)\I_{\Ero_{k',x}}\le C_{\lambda_*,\eps}\exp\left(\frac{(\lambda^2+\eps)(n-k)}{2}-\lambda t-\frac{t^2}{C_{\lambda_*,\eps}(n-k)^3}\right).
\end{equation}
Using \eqref{e:upper_barrier_laststep3} (with $k=k'$) to estimate the innermost probability in \eqref{e:upper_barrier_laststep2}, we can now continue just like in the proof of Lemma \ref{l:upper_barrier_late}. We omit further details.
\end{proof}

In the proof of Lemma \ref{l:upper_barrier_late}, we used the following characterisation of the tilted $\Ind_k$.
  \begin{lemma}\label{l:tilting}
Let $|\lambda|\le\lambda_*$, and consider the tilted measure
\[\frac{d\Q}{d\tilde\PP^{Q_N,0}} :=\frac{ \exp\left(\lambda \Ind_k (x,\phi)\right) }{\tilde\E^{Q_N,0}\left(\exp\left(\lambda \Ind_k (x,\phi)\right) \right)}\]
Then under $\Q$  the random variable $\Ind_k(x)$ has mean $\g\lambda+O_{\lambda_*}\left(\frac{1}{r_k^\gamma}\right)$ and variance $\g+O_{\lambda_*}\left(\frac{1}{r_k^\gamma}\right)$.

Moreover, the law of $\Ind_k(x)$ under $\Q$ is absolutely continuous with respect to Lebesgue measure, with a density that is bounded above by some constant $C$ (independent of $k$).
   \end{lemma}
 \begin{proof}
  Denote by $U_i = X_i - \g\lambda$.  We first notice that $|\E_\Q [U_k]| \le Ce^{-\gamma (n-k) }$,  where $\gamma >0$ is the constant from Theorem \ref{t:qclt}.  Indeed,  by the definition of $\Q$,
 $$
 \E_\Q [U_k] = \frac {\tilde\E^{Q_N,0}[(X_k -  \g \lambda)e^{\lambda X_k}]}{\tilde\E^{Q_N,0}[e^{\lambda X_k}]}.
 $$
 It follows from Lemma \ref{l.meantilt} that
 \begin{align*}
 \frac{d}{d\lambda} \E^{Q_N,0}[e^{\lambda X_k}]
 =
 \g \lambda \E^{Q_N,0}[e^{\lambda X_k}] + O(e^{-\gamma (n-k)}) \lambda.
  \end{align*}
  Therefore we conclude the estimate for the mean.  For the variance,  we notice that for any $s\in\R$, it follows  from Theorem \ref{t:qclt} that
  \begin{align*}
  \E_\Q [e^{s X_k}] = \frac{\E^{Q_N,0} [e^{(s+\lambda) X_k}]}{\E^{Q_N,0} [e^{\lambda X_k}]}
  =
  e^{(\lambda s + \frac 12 \lambda^2)\g} (1+ O(e^{-\gamma (n-k) }(s^2+\lambda^2)),
   \end{align*}
   from which we conclude that
   \[
   \Var_\Q [U_k] =   \Var_\Q [X_k] = \g + O(e^{-\gamma (n-k) }).
   \]

It remains to show the boundedness of the density of $X_k$.  It is known that both marginals and linear functions of a log-concave distribution are log-concave \cite{P73}, \cite{BL76}, and therefore $X_k$ has a log-concave density, and the latter is bounded by $1/\sqrt{\Var(X_k)}$ by, e.g., the upper bound in \cite[(1.4)]{BC15}.
\end{proof}

\subsection{Completion of the proof}
After the hard work in the previous subsections, we can now establish the desired upper bound on the maximum.

\begin{theorem}\label{t:upper_bound_max}
The random variable $\max_{x\in Q_N}\phi(x)-\sqrt{\g}m_N$ has a tight right tail. That is, for each $\eps>0$ there is $C_\eps>0$ such that
\begin{equation}\label{e:upper_bound_max}
\PP^{Q_N,0}\left(\max_{x\in Q_N}\phi(x)-\sqrt{\g}m_N\le C_\eps\right)=\tilde\PP^{Q_N,0}\left(\max_{x\in Q_N}\phi(x)-\sqrt{\g}m_N\le C_\eps\right)\ge 1-\eps
\end{equation}
for all $N\in\N$.
\end{theorem}
\begin{proof}
Let $A_{(\log n)^{5/2}}$ be as in \eqref{e:def_Ak}. By Lemma \ref{l:pointsnearbdry}, with probability at least $1-\frac1n$ the maximum of $\phi$ over points in $Q_N\setminus A_{(\log n)^{5/2}}$ is at most $\sqrt{\g}m_N$. So it suffices to consider the maximum of the field over $A_{(\log n)^{5/2}}$.

From Lemma \ref{l:upper_barrier_early} we know that
\[\tilde\PP^{Q_N,0}\left(\bigcap_{x\in Q_N}\Eup_{\le k_\infty-(\log n)^3,x}\right)\ge 1-C\exp\left(-\frac{\Gamma}{C}\right).\]
Applying Lemma \ref{l:upper_barrier_late} iteratively for $k\in\left\{k_\infty-(\log n)^3,\ldots,k_\infty'\right\}$ we we can conclude from this that even
\[\tilde\PP^{Q_N,0}\left(\bigcap_{x\in Q_N}\Eup_{\le k_\infty',x}\right)\ge 1-C\exp\left(-\frac{\Gamma}{C}\right)-C\sum_{k=n/2}^{k_\infty'}\exp\left(-\frac{(n-k+1)^{2/5}+\Gamma}{C}\right)\ge 1-C'\exp\left(-\frac{\Gamma}{C}\right).\]
Finally, using Lemma \ref{l:upper_barrier_laststep}, we see that 
\[\tilde\PP^{Q_N,0}\left(\bigcap_{x\in A_{(\log n)^{5/2}}}\Eup_{\le k_{\infty'},x} \cap \{\phi\colon \phi(x)\le \sqrt{\g}m_N+\Gamma\}\right)\ge 1-C\exp\left(-c\Gamma\right).\]
Now for a given $\eps$ we can choose $\Gamma$ large enough so that $C\exp\left(-c\Gamma\right)\le\eps$, and pick $C_\eps=\Gamma$. Then \eqref{e:upper_bound_max} follows.
\end{proof}

\section{Lower bound on the maximum}\label{s:lower}

In order to establish the lower bound for the maximum of the field, we will apply a truncated second-moment method. To do so, we need to introduce  some definitions. We try to follow the notation used in the proof of the upper bound in Section \ref{s:upper}.

As we are only interested in a lower bound on the maximum, we directly restrict attention to a subset of $Q_N$ with particularly nice properties. We define $k_0=\Delta$ for some constant $\Delta$ to be chosen later, and
\begin{align*}
I_{k_0}&=10r_{k_0}\Z^2\cap Q_{N/2},\\
A_{k_0}&=\bigcup_{y\in I_{k_0}}Q_{r_{k_0}}(y).
\end{align*}
Thus $A_{k_0}$ is the union of $\ge c\left(\frac{N}{r_{k_0}}\right)^2$ squares of sidelength $2r_{k_0}$ that are well-separated and far away from the boundary of $Q_N$. We will show that already the maximum of $\phi$ over $A_{k_0}$ is with high probability of the right order.

It will be easy to see that the maximum of $\left|\Sav_{k_0,+}(x, \phi)\right|$ over $x\in A_{k_0}$ is bounded by $Ck_0=C\Delta$ with high probability, and so it is negliglible for the purpose of proving tightness from below.

The main challenge will be to show that the probability that there is a point in $x\in A_{k_0}$ where $\phi(x)-\Sav_{k_0,-}(x, \phi)$ is of order $\sqrt{\g}m_N-C(\Delta+1)$ tends to 1 as $\Delta\to\infty$.
To do so, we will apply the second-moment method to estimate the number of points $x$ where $\phi(x)-\Sav_{k_0,+}(x, \phi)$ is large, while the square average process stays below a certain barrier that is curved downwards (and some other conditions are satisfied).

To introduce the corresponding barrier event, we introduce more parameters. We let $\ell$ be a parameter that will later be chosen as a large constant (independent of $\Delta,N$) and let $k_\infty=n-\Delta$. We also let $k_0'=k_0+\ell$ and $k_\infty'=k_\infty-\ell$ and consider the barrier functions
\begin{align*}
\Blwp(j)&:=\frac{\sqrt{\g}(m_N-4\Delta)}{k_\infty-k_0}(k-k_0)-(k\wedge(k_\infty-k))^{2/5},\\
\Blwm(j)&:=\frac{\sqrt{\g}(m_N-4\Delta)}{k_\infty-k_0}(k-k_0)-(k\wedge(k_\infty-k))^{3/5},
\end{align*}
and the event
\begin{equation}
\begin{split}
  \label{eq-Elw}
&\Elw_x=\\
&\vast\{\phi\colon \sum_{i=k_0}^{j-1}\Inc_i(x,\phi)\in\left[\Blwm(j),\Blwp(j)\right]\forall j\in\{k_0',\ldots k_\infty'\},\;
\sum_{i=k_0}^{k_\infty-1}\Inc_i(x,\phi)\in\sqrt{\g}(m_N-4\Delta)+[-1,0]\vast\}.
\end{split}
\end{equation}
For later use, let us also define truncations of this event. For $k\in\{k_0,\ldots k_\infty-1\}$ we define
\[
\Elw_{\le k,x}=\vast\{\phi\colon \sum_{i=k_0}^{j-1}\Inc_i(x,\phi)\in\left[\Blwm(j),\Blwp(j)\right]\forall j\in\{k_0',\ldots k_\infty'\wedge k\}\vast\}
\]
and
\begin{align*}
&\Elw_{\ge k,x}=\\
&\vast\{\phi\colon \exists a\in\sqrt{\g}(m_N-4\Delta)+[-1,0]\text{ s.t. } a-\sum_{i=j}^{k_\infty-1}\Inc_i(x,\phi)\in\left[\Blwm(j),\Blwp(j)\right]\forall j\in\{k_0\vee k,\ldots k_\infty'\}\vast\}.
\end{align*}

We also define the events
\begin{align*}
\Ecp_x&=\bigcap_{k_0\le k\le k_\infty-1}\Ecp_{k,x},\qquad
\Ero_x=\bigcap_{k_0\le k\le k_\infty}\Ero_{k,x},\qquad
\Ebd_x=\bigcap_{k_0\le k\le k_\infty}\Ebd_{k,x}
\end{align*}
and, for $k\in\{k_0,\ldots k_\infty-1\}$, the events
\begin{align*}
\Ecp_{\le k,x}&=\bigcap_{k_0\le j\le k}\Ecp_{k,x},\qquad
\Ero_{\le k,x}=\bigcap_{k_0\le j\le k}\Ero_{k,x},\qquad
\Ebd_{\le k,x}=\bigcap_{k_0\le j\le k}\Ebd_{k,x}.
\end{align*}

We point out a conceptual difference between the definitions of our upper and lower barrier events. We defined the upper barrier event $\Eup_{k,x}$ in terms of the square averages $\Sav_k(x,\phi)$, while the lower barrier event $\Elw_{k,x}$ is defined using the sum of increments $\sum_{i=k_0}^{k_\infty-1}\Inc_i(x,\phi)$. There are two reasons for making these choices. The first is that when proving a lower bound, we can assume that certain good events (such as the event $\Ebd_x$ that all boundary layer terms are small) occur. This allows to work with the sum of increments (that has a nice structure as a sum of near-independent increments). In the proof of the upper bound, this is not easily possible. The second reason is that in the proof of the lower bound there is a need to ``cut the tree'', i.e. to make sure that the barrier events for two points at macroscopic distance are independent. This we achieve by considering the sum of increments only starting at $j=k_0=\Delta$, and eventually choosing $\Delta$ large.

Note in particular that on the event $\Ebd_{x,\le k}$ the random variable $\Sav_{k,+}(x,\phi)-\Sav_{k_0,+}(x,\phi)$ is well approximated by $\sum_{j=k_0}^{k-1}\Inc_j(x,\phi)$, while on the event $\Ero_{x,\le k}$ we can couple the $\Inc_j$ to independent $\Ind_j$ as in Lemma \ref{l:iterated_coupling}. Finally on the event $\Ecp_{x,\le k}$ that all these couplings succeed, the definition of $\Elw_{x,\le k}$ simplifies as it can now be written in terms of the $\Inc_j$. In particular on the event $\Ecp_x\cap\Ero_x\cap\Ero_{k,x}$ the barrier event $\Elw_x$ can be phrased in terms of the (genuinely independent) increments $\Inc_j$.

We can now consider the counting random variable
\[\mathcal{N}_{\GGamma,N}=\sum_{x\in A_{k_0}}\I_{\Elw_x}\I_{\Ero_x}\I_{\Ebd_x}\I_{\Ecp_x}\I_{\phi(x)-\Sav_{k_0,-}(x,\phi)\ge \sqrt{\g}m_N-\GGamma}\]
Here $\GGamma$ is another constant (that we will eventually choose as $C\Delta$ for an absolute constant $C$).

We will show that there the probability that $\mathcal{N}_{\GGamma,N}$ is non-zero tends to 1 as $\Delta\to\infty$. In fact, we will show that this is the case conditionally on the field outside of $A_{k_0}$. We set
\[\tilde\F:=\tilde\F_{Q_N\setminus\bigcup_{y\in I_{k_0}}Q_{r_{k_0-1}}(y)}\]
and
\[\Ero_{A_{k_0}}:=\bigcap_{y\in I_{k_0}}\Ero_{k_0-1,y}.\]
To show that $\mathcal{N}_{\GGamma,N}$ is non-zero, we need a lower bound on its expectation and an upper bound on its variance, conditional on $\F$. This is the subject of the following two subsections.

Before turning to this, let us state two technical estimates.
\begin{lemma}\label{l:roughnessAk0}
We have
\begin{equation}\label{e:roughnessAk0}
\tilde\PP^{Q_N,0}\left((\Ero_{A_{k_0}})^\com\right)\le C\exp\left(-c(\log r_{k_0})^{3}\right).
\end{equation}
Moreover, for any $x\in A_{k_0}$ we have
\begin{equation}\label{e:roughnessAk0cond}
\tilde\PP^{Q_N,0}\left((\Ero_{k_0,x})^\com\mid\tilde\F\right)I_{\Ero_{A_{k_0}}}\le C\exp\left(-c(\log r_{k_0})^{3}\right).
\end{equation}
\end{lemma}
\begin{proof}
The estimate \eqref{e:roughnessAk0} follows from \eqref{e:badevents1} and a union bound.

Regarding \eqref{e:roughnessAk0cond}, it suffices to prove that for $y\in I_{k_0}$ such that $x\in Q_{r_{k_0}}(y)$,
\begin{equation}\label{e:roughnessAk01}
\tilde\PP^{Q_N,0}\left((\Ero_{k_0,x})^\com\mid\tilde\F_{Q_{r_{k_0}}^\com(x)}\right)I_{\Ero_{k_0-1,y}}\le C\exp\left(-c(\log r_{k_0})^{3}\right).
\end{equation}
To see \eqref{e:roughnessAk01}, we argue very similarly as in the proof of estimate \eqref{e:badevents3} in Lemma \ref{l:badevents}. The only difference is that $Q_{r_{k_0}}(x)$ and $Q_{r_{k_0-1}}(y)$ need not be concentric. However we do have
\[Q_{r_{k_0}}(x)\subset Q_{r_{2k_0}}(y)\subset Q_{r_{k_0-1}}(y),\]
and so we can still use the discrete Harnack inequality to obtain decay of the oscillation of the harmonic extension of $\phi\restriction_{\partial^+Q_{r_{k_0-1}}(y)}$.
\end{proof}

\subsection{Lower bound on the first moment}
We begin by giving a lower bound on the expectation of $\mathcal{N}_{\GGamma,N}$.

\begin{lemma}\label{l:lower_bound_first}
There is a constant $\ell>0$ such that for any $\Delta>0$ and for any $N$ sufficiently large (depending on $\Delta,\ell$) we have for $x\in A_{k_0}$ the estimate
\begin{equation}\label{e:lower_bound_first}
\tilde\PP^{Q_N,0}\left(\Elw_x,\Ero_x,\Ebd_x,\Ecp_x,\phi(x)-\Sav_{k_0,-}(x,\phi)\ge \sqrt{\g}(m_N-C(\Delta+1))\middle|\tilde\F\right)\I_{\Ero_{A_{k_0}}}\ge c\frac{\e^{4\Delta}}{N^2}\I_{\Ero_{A_{k_0}}}.
\end{equation}
In particular, there is a constant $C$, such that if $\GGamma\ge C(\Delta+1)$ then for any $N$ sufficiently large,
\begin{equation}\label{e:lower_bound_firstmoment}
\tilde\E^{Q_N,0}(\mathcal{N}_{\GGamma,N}|\tilde\F)\I_{\Ero_{A_{k_0}}}\ge c\e^{4\Delta}\I_{\Ero_{A_{k_0}}}.
\end{equation}
\end{lemma}

Even though this is a lower, not an upper bound, the proof is quite similar to the proof of Lemma \ref{l:upper_barrier_late}. 
The main ingredient is the following barrier estimate for the tilted process, whose proof 
appears  in  Section \ref{s:ballot_thms}.
 \begin{lemma}
 \label{l:lower_barrier}
Let $m\in\N$ and $\gamma>0$. Let $(X_j)_{j=0}^m$ be independent random variables possessing a uniformly bounded density,  $C<\infty$ such that $\left|\E_\Q[X_j]\right| \le C\exp(\gamma(j-m))$ and $\left|\E_\Q[X_j^2]-\g\right|\le C\exp(\gamma(j-m))$ for any $j\ge 0$.  We further assume that there exists $\lambda_*> 0$ and $C^*<\infty$  such that $\E\exp(\lambda_* X_j^2)< C^*$ for all $|\lambda|<\lambda_*$ and all $j\ge0$.  Let $\Sigma_k = \sum_{j=0}^{k-1} X_j$.

Then there exists some $\ell\ge1$ fixed and constants $c_{\gamma,\ell},C_{\gamma,\ell}\ge 0$ depending on $\gamma$ and $\ell$ only such that
\begin{align*}
\frac{ c_{\gamma,\ell}}{m^{3/2}}&\le \PP \left(\Sigma_j \in\left[-(j\wedge (m-j))^{3/5}, -(j\wedge (m-j))^{2/5}\right] \ \forall j\in\{\ell,\ldots m-\ell\}, \Sigma_{m-1} \in [ -1, 0] \right)\\
&\qquad\qquad\qquad\qquad\qquad  \le \frac{C_{\gamma,\ell}}{m^{3/2}}.
\end{align*}
  Moreover for any $k\in\{\ell,\ldots,m-\ell\}$ and any $t>0$ we have the estimates
  \begin{align*}
 \PP \left(\Sigma_j\le 0\ \forall j\in\{\ell,\ldots,k-1\},\Sigma_{k}\in [-t-1,-t]\right)&\le \frac{C_{\gamma,\ell}(1+t)}{k^{3/2}},\\
  	\PP \left(\Sigma_{m-1}-\Sigma_j\le 0\ \forall j\in\{k,\ldots,m-\ell\},\Sigma_{m-1}-\Sigma_{k}\in [-t-1,-t]\right)&\le \frac{C_{\gamma,\ell}(1+t)}{(m-k)^{3/2}}.
  \end{align*}
 \end{lemma}

\begin{proof}[Proof of Lemma \ref{l:lower_bound_first}]
We first note that \eqref{e:lower_bound_firstmoment} is an easy consequence of \eqref{e:lower_bound_first} and Lemma \ref{l:badevents}, and so the main challenge is to establish \eqref{e:lower_bound_first}.
We fix the constant $\ell$ so that Lemma \ref{l:lower_barrier} can be applied.

Note that on the event $\Ecp_x$, we have that $\Inc_j(x,\phi)=\Ind_j(x)$, and so on $\Ecp_x$ we have $\Elw_x=\widetilde{\Elw_x}$, where $\widetilde{\Elw_x}$ is the same event, just with $\Ind_j$ instead of $\Inc_j$. That is
\begin{equation}\label{e:def_tildeElw}
\begin{split}
&\widetilde{\Elw_x}\\
&=\vast\{\phi\colon \sum_{j=k_0}^{k-1}\Ind_j(x)\in\frac{\sqrt{\g}(m_N-4\Delta)}{k_\infty-k_0}(k-k_0)+\left[-((k-k_0)\wedge (k_\infty-k))^{3/5},-((k-k_0)\wedge (k_\infty-k))^{2/5}\right]\\
&\qquad\qquad\forall k\in\{k_0',\ldots k_\infty'\},\quad\sum_{j=k_0}^{k_\infty-1}\Ind_j(x)\in\sqrt{\g}(m_N-4\Delta)+[-1,0]\vast\}.
\end{split}
\end{equation}

Consider the tilted measure
\begin{equation}\label{e:def_Q}
\frac{d\Q}{d\PP(\cdot\mid\tilde\F)} :=\frac{ \exp\left(\lambda \left(\sum_{k=k_0}^{k_\infty-1} \Ind_k (x,\phi)\right) \right)}{\E\left[\exp\left(\lambda \left(\sum_{k=k_0}^{k_\infty-1} \Ind_k (x,\phi)\right) \right) \middle|\tilde\F\right]}\end{equation}
where
\begin{equation}\label{e:def_lambda2}
\lambda:=\frac{\sqrt{\g}(m_N-4\Delta)}{\g (k_\infty-k_0)}=\frac{m_N-4\Delta}{\sqrt{\g} (n-2\Delta)}=\frac{2}{\sqrt{\g}}-\frac{3\log n}{4\sqrt{\g}(n-2\Delta)}.
\end{equation}

As in the proof of Lemma \ref{l:tilting}, we see that under $\Q$ the $\Ind_k-\lambda$ are independent random variables with expectation $O\left(\frac{1}{r_k^\gamma}\right)$ and variance $\g+O\left(\frac{1}{r_k^\gamma}\right)$. 
Lemma \ref{l:lower_barrier} with $m=k_\infty-k_0\ge\frac n2$ implies that
\[\Q\left(\widetilde{\Elw_x}\right)\ge \frac{c}{(k_\infty-k_0)^{3/2}}\ge \frac{c}{n^{3/2}},\]
with the constants independent of $\Delta$ for large enough $n$.
As on the event $\widetilde{\Elw_x}$ we know in particular that $\sum_{k=k_0}^{k_\infty-1}\Ind_k(x)\le\sqrt{\g}(m_N-4\Delta)$, we can calculate
\begin{align*}
&\tilde\PP^{Q_N,0}\left(\widetilde{\Elw_x}\mid\tilde\F\right)\\
&\ge\tilde\E^{Q_N,0}\left(\I_{\widetilde{\Elw_x}}\exp\left(\lambda \left(\sum_{j=k_0}^{k_\infty-1} \Ind_j (x)-\sqrt{\g}(m_N-4\Delta)\right)\right)\middle|\tilde\F\right)\\
&\ge\exp\left(-\lambda \sqrt{\g}(m_N-4\Delta)\right)\tilde\E^{Q_N,0}\left(\I_{\widetilde{\Elw_x}}\exp\left(\lambda \sum_{j=k_0}^{k_\infty-1} \Ind_j (x)\right)\middle|\tilde\F\right)\\
&=\exp\left(-\lambda \sqrt{\g}(m_N-4\Delta)\right)\tilde\E^{Q_N,0}\left(\exp\left(\lambda \sum_{j=k_0}^{k_\infty-1} \Ind_j (x)\right)\middle|\tilde\F\right)\Q\left(\widetilde{\Elw_x}\right).
\end{align*}
We can compute the expectation here using Lemma \ref{l:exponential_moment_bounds} and obtain that
\begin{equation}\label{e:lower_bound_first1}
\begin{split}
&\tilde\PP^{Q_N,0}\left(\widetilde{\Elw_x}\mid\tilde\F\right)\\
&\ge \frac{c_\ell}{n^{3/2}}\exp\left(-\lambda \sqrt{\g}(m_N-4\Delta)+\frac{\lambda^2\g}{2}(k_\infty-k_0)-C\sum_{j=k_0}^{k_\infty-1}\frac{1}{r_k^\omega}\right)\\
&\ge c\exp\left(-\frac32\log n-\frac{(m_N-4\Delta)^2}{2(k_\infty-k_0)}\right)\\
&\ge c\exp\left(-\frac32\log n-2(n-4\Delta)+\frac32\log n-O_{\Delta}\left(\frac{\log n}{n}\right)\right)\\
&\ge c\exp\left(-2n+8\Delta\right),
\end{split}
\end{equation}
for any $N$ large enough.

In the following we will try to argue that on the event $\Ero_{A_{k_0}}$ we can add in the events $\Ero_x, \Ebd_x, \Ecp_x$ only changing the probability by a constant factor. As Lemma \ref{l:badevents} and Lemma \ref{l:exponential_moment_bounds} give of bounds for the conditional probabilities of $\Ebd_{k,x}$ and $\Ecp_{k,x}$ only on the event $\Ero_{k,x}$, it is important to start with $\Ero_x$. Again, we have a bound for the conditional probability of $\Ero_{k,x}$ only on the event $\Ero_{k-1,x}$, and so we need to do this inductively. We have that
\begin{equation}\label{e:lower_bound_first2}
\begin{split}
&\tilde\PP^{Q_N,0}\left(\widetilde{\Elw_x},(\Ero_x)^\com\middle|\tilde\F\right)\I_{\Ero_{A_{k_0}}}\\
&=\tilde\PP^{Q_N,0}\left(\widetilde{\Elw_x},(\Ero_{x,k_0})^\com\middle|\tilde\F\right)\I_{\Ero_{A_{k_0}}}+\sum_{k=k_0}^{k_\infty-1}\tilde\PP^{Q_N,0}\left(\widetilde{\Elw_x},\Ero_{\le k,x},(\Ero_{x,k+1})^\com\middle|\tilde\F\right)\I_{\Ero_{A_{k_0}}}.
\end{split}
\end{equation}
Here the first summand on the right-hand side is controlled by Lemma \ref{l:roughnessAk0}. For the other summands, we use that on the event $\widetilde{\Elw_x}$ we have that $\sum_{k=k_0}^{k_\infty-1} \Ind_k \ge\sqrt{\g}(m_N-4\Delta)-1$. We also drop the barrier starting from $k$, i.e. we use that $\widetilde{\Elw_x}\subset \widetilde{\Elw_{\le k-1,x}}$, and choose $\lambda$ as in \eqref{e:def_lambda2}.

Now by iterative conditioning (relying on Lemma \ref{l:iterated_coupling}) we find that with
\begin{align*}
&\tilde\PP^{Q_N,0}\left(\widetilde{\Elw_x},\Ero_{\le k,x},(\Ero_{k+1,x})^\com\middle|\tilde\F\right)\I_{\Ero_{A_{k_0}}}\\
&\le \tilde\E^{Q_N,0}\left(\exp\left(\lambda\left(\sum_{k=k_0}^{k_\infty-1}\Ind_k(x)-\sqrt{\g}(m_N-4\Delta)+1\right)\right)\I_{\widetilde{\Elw_{\le k-1,x}}}\I_{\Ero_{\le k,x}}\I_{(\Ero_{x,k+1})^\com}\middle|\tilde\F\right)\I_{\Ero_{A_{k_0}}}\\
&\le\exp\left(-\lambda(\sqrt{\g}(m_N-4\Delta)+1)\right)\\
&\qquad\times\E\Bigg(
\tilde\E^{Q_N,0}\left(
\tilde\E^{Q_N,0}\left(\exp\left(\lambda\sum_{j=k+1}^{k_\infty-1}\Ind_j(x)\right)\middle|\tilde\F_{Q_{r_{k+1}}^\com(x)}\right)
\exp(\lambda \Ind_{k}(x,\phi))\I_{(\Ero_{k+1,x})^\com}
\middle|\tilde\F_{Q_{r_{k}}^\com(x)}\right)\\
&\qquad\qquad\exp\left(\lambda\sum_{j=k_0}^{k-1}\Ind_j(x)\right)\I_{\widetilde{\Elw_{\le k-1,x}}}\I_{\Ero_{\le k,x}}\Bigg|\tilde\F\Bigg)\I_{\Ero_{A_{k_0}}}.
\end{align*}
This tower of conditional expectations can be evaluated from the inside out. First of all we have
\[\tilde\E^{Q_N,0}\left(\exp\left(\lambda\sum_{j=k+1}^{k_\infty-1}\Ind_j(x)\right)\middle|\tilde\F_{Q_{r_{k+1}}^\com(x)}\right)\le C\exp\left(\frac{\lambda^2\g}{2}(k_\infty-k+1)\right)\]
by Lemma \ref{l:exponential_moment_bounds}. By the same lemma and the Cauchy-Schwarz inequality we have
\[\tilde\E^{Q_N,0}\left(\exp(\lambda \Ind_k(x,\phi))\I_{(\Ero_{k+1,x})^\com}
\middle|\tilde\F_{Q_{r_{k}}^\com(x)}\right)\I_{\Ero_{k,x}}\le C\exp(-c(\log r_{k+1})^3)\]
and so we obtain
\begin{equation}\label{e:lower_bound_first3}
\begin{split}
&\tilde\PP^{Q_N,0}\left(\widetilde{\Elw_x},\Ero_{\le k,x},(\Ero_{k+1,x})^\com\middle|\tilde\F\right)\I_{\Ero_{A_{k_0}}}\\
&\le C\exp\left(-\lambda(\sqrt{\g}(m_N-4\Delta)+1)\right)\exp\left(\frac{\lambda^2\g}{2}(k_\infty-k+1)-c(\log r_{k+1})^3\right)\\
&\qquad\tilde\E^{Q_N,0}\left(\exp\left(\lambda\sum_{j=k_0}^{k-1}\Ind_k(x)\right)\I_{\widetilde{\Elw_{\le k-1,x}}}\middle|\tilde\F\right).
\end{split}
\end{equation}
To evaluate the remaining expectation, we distinguish two cases. In the easy case $k<k_0'$ we can ignore $\I_{\widetilde{\Elw_{\le k-1,x}}}$ and just use Lemma \ref{l:exponential_moment_bounds} once more to obtain
\begin{equation}\label{e:lower_bound_first4}
\tilde\E^{Q_N,0}\left(\exp\left(\lambda\sum_{j=k_0}^{k-1}\Ind_j(x)\right)\I_{\widetilde{\Elw_{\le k-1,x}}}\middle|\tilde\F\right)\le C\exp\left(\frac{\lambda^2\g}{2}(k-k_0)\right).
\end{equation}
In case $k\ge k_0'$ it would be to wasteful to ignore $\I_{\widetilde{\Elw_{\le k-1,x}}}$. Instead we use a tilting argument very similar as the one that lead to \eqref{e:lower_bound_first1}. Namely with $\Q$ defined as before, Lemma \ref{l:lower_barrier} with $k\wedge k_\infty'$ in place of $k$ implies that
\[\Q\left(\widetilde{\Elw_{\le k-1,x}}\right)\le C_\ell\frac{(k\wedge k_\infty'\wedge(k_\infty-(k\wedge k_\infty')))^{3/10}}{(k\wedge k_\infty')^{3/2}}\]
and therefore
\begin{equation}\label{e:lower_bound_first5}
\begin{split}
\tilde\E^{Q_N,0}\left(\exp\left(\lambda\sum_{j=k_0}^{k-1}\Ind_k(x)\right)\I_{\widetilde{\Elw_{\le k-1,x}}}\middle|\tilde\F\right)&\le \tilde\E^{Q_N,0}\left(\exp\left(\lambda \sum_{j=k_0}^{k-1} \Ind_k (x,\phi)\right)\middle|\tilde\F\right)\Q\left(\widetilde{\Elw_{\le k-1,x}}\right)\\
&\le C\exp\left(\frac{\lambda^2\g}{2}(k-k_0)\right)\frac{(k\wedge k_\infty'\wedge(k_\infty-(k\wedge k_\infty')))^{3/10}}{(k\wedge k_\infty')^{3/2}}.
\end{split}
\end{equation}
We can now use \eqref{e:lower_bound_first3} together with \eqref{e:lower_bound_first4} (for $k<k_0'$) or \eqref{e:lower_bound_first5} (for $k\ge k_0'$), and obtain from \eqref{e:lower_bound_first2} that
\begin{align*}
&\tilde\PP^{Q_N,0}\left(\widetilde{\Elw_x},(\Ero_x)^\com\middle|\tilde\F\right)\I_{\Ero_{A_{k_0}}}\\
&\le C\exp(-c(\log r_{k_0})^3)\\
&+C\exp\left(\frac{\lambda^2\g}{2}(k_\infty-k_0)-\lambda\sqrt{\g}(m_N-4\Delta)\right)\\
&\qquad\times\sum_{k=k_0}^{k_\infty-1}\exp(-c(\log r_{k+1})^3)\left(C\I_{k<k_0'}+\frac{(k\wedge k_\infty'\wedge(k_\infty-(k\wedge k_\infty')))^{3/10}}{(k\wedge k_\infty')^{3/2}}\I_{k\ge k_0'}\right).
\end{align*}
The sum on the right-hand side is dominated by its last terms and we obtain
\begin{equation}\label{e:lower_bound_first6}
\begin{split}
\tilde\PP^{Q_N,0}\left(\widetilde{\Elw_x},(\Ero_x)^\com\middle|\tilde\F\right)\I_{\Ero_{A_{k_0}}}&\le C\exp(-100n)+C\exp\left(-\frac{(m_N-4\Delta)^2}{2(k_\infty-k_0)}\right)\frac{\exp(-c\Delta^3)}{n^{3/2}}\\
&\le C\exp(-2n+8\Delta)\exp(-c\Delta^3).
\end{split}
\end{equation}
If $\Delta$ is large enough, then the right-hand side of \eqref{e:lower_bound_first6} is less than half the right-hand side of \eqref{e:lower_bound_first1}. So by combining these two estimates we find that
\begin{equation}\label{e:lower_bound_first7}
\tilde\PP^{Q_N,0}\left(\widetilde{\Elw_x},\Ero_x\middle|\tilde\F\right)\I_{\Ero_{A_{k_0}}}
\ge c\exp(-2n+8\Delta)
\end{equation}
for any $N$ large enough.

Now that we have succeeded in introducing the event $\Ero_x$, we can repeat the process twice more to introduce $\Ebd_x$ and $\Ecp_x$ as well. In fact, the argument is slightly simpler, as we now can assume that $\Ero_x$ holds, so that the conditional estimates from Lemma \ref{l:badevents} and Lemma \ref{l:iterated_coupling} are directly applicable. In this manner one can obtain the estimates
\begin{align*}
\tilde\PP^{Q_N,0}\left(\widetilde{\Elw_x},\Ero_x,(\Ebd_x)^\com\middle|\tilde\F\right)\I_{\Ero_{A_{k_0}}}
&\le C\exp(-2n+8\Delta)\exp\left(-ce^{\Delta\omega}\right),\\
\tilde\PP^{Q_N,0}\left(\widetilde{\Elw_x},\Ero_x,(\Ecp_x)^\com\middle|\tilde\F\right)\I_{\Ero_{A_{k_0}}}
&\le C\exp(-2n+8\Delta)\exp(-\Delta\omega).
\end{align*}
If we combine these two estimates with \eqref{e:lower_bound_first7} and enlarge $\Delta$ if necessary, we conclude
\begin{equation}\label{e:lower_bound_first8}
\tilde\PP^{Q_N,0}\left(\widetilde{\Elw_x},\Ero_x,\Ebd_x,\Ecp_x\middle|\tilde\F\right)\I_{\Ero_{A_{k_0}}}
\ge c\exp(-2n+8\Delta).
\end{equation}

From here on it is straightforward to complete the proof. Already at the beginning of the proof we had observed that $\widetilde{\Elw_x}\cap \Ecp_x=\Elw_x\cap\Ecp_x$. Moreover, on the event $\Ebd_x$ we have
\[\Ind_{k_0}(x)+\ldots+\Ind_{k_\infty-1}(x)=\Sav_{k_\infty,+}(x,\phi)-\Sav_{k_0,-}(x,\phi)+O(1)\]
and
\[\Sav_{k_\infty,-}(x,\phi)-\Sav_{k_\infty,+}(x,\phi)\ge-C.\]
So on the event $\widetilde{\Elw_x}\cap\Ero_x\cap\Ebd_x\cap\Ecp_x=\Elw_x\cap\Ero_x\cap\Ebd_x\cap\Ecp_x$ we know that
\[\Sav_{k_\infty,-}(x,\phi)-\Sav_{k_0,+}(x,\phi)\ge \sqrt{\g}(m_N-4\Delta)-C',\]
for some constant $C'$.

Additionally note that by \eqref{e:exponential_moment_bounds3} and the exponential Markov inequality there is a constant $C''$ such that
\[\tilde\PP^{Q_N,0}\left(\phi(x)-\Sav_{k_\infty,-}(x,\phi)\ge -C''\Delta\mid\tilde\F_{Q_{r_{k_\infty}}^\com(x)}\right)\I_{\Ero_{k_\infty,x}}\ge \frac34.\]
As Lemma \ref{l:iterated_coupling} implies that (for $\Delta$ large enough) we have
\[\tilde\PP^{Q_N,0}\left(\Ebd_{k_\infty,x}\middle|\tilde\F_{Q_{r_{k_\infty}}^\com(x)}\right)\I_{\Ero_{k_\infty,x}}\ge \frac34,\]
we know that also
\[\tilde\PP^{Q_N,0}\left(\phi(x)-\Sav_{k_\infty,-}(x,\phi)\ge -C''\Delta,\Ebd_{k_\infty,x}\middle|\tilde\F_{Q_{r_{k_\infty}}^\com(x)}\right)\I_{\Ero_{k_\infty,x}}\ge \frac12.\]

So from \eqref{e:lower_bound_first8} we conclude that
\begin{align*}
&\tilde\PP^{Q_N,0}\left(\Elw_x,\Ero_x,\Ebd_x,\Ecp_x,\phi(x)-\Sav_{k_0,-}(x,\phi)\ge \sqrt{\g}(m_N-4\Delta)-C'-C''\Delta\mid\tilde\F\right)\I_{\Ero_{A_{k_0}}}\\
&\ge\tilde\PP^{Q_N,0}\Bigg(\Elw_x,\Ero_x,\Ebd_x,\Ecp_x,\phi(x) - \Sav_{k_\infty,-}(x,\phi)\ge -C''\Delta,\\
&\qquad\qquad\Sav_{k_\infty,-}(x,\phi)-\Sav_{k_0,+}(x,\phi)\ge \sqrt{\g}(m_N-4\Delta)-C'\Bigg|\tilde\F\Bigg)\I_{\Ero_{A_{k_0}}}\\
&=\tilde\PP^{Q_N,0}\Bigg(\Elw_x,\Ero_x,\Ebd_x,\Ecp_x,\phi(x) - \Sav_{k_\infty,-}(x,\phi)\ge -C''\Delta\Bigg|\tilde\F\Bigg)\I_{\Ero_{A_{k_0}}}\\
&\ge\tilde\E^{Q_N,0}\Bigg(\tilde\PP^{Q_N,0}\left(\phi(x)-\Sav_{k_\infty,-}(x,\phi)\ge -C''\Delta,\Ebd_{k_\infty,x}\middle|\tilde\F_{Q_{r_{k_\infty}}^\com(x)}\right)\I_{\Elw_x}\I_{\Ero_x}\I_{\Ebd_{\le k_\infty-1,x}}\I_{\Ecp_x}\Bigg|\tilde\F\Bigg)\I_{\Ero_{A_{k_0}}}\\
&\ge\frac12\tilde\PP^{Q_N,0}\Bigg(\Elw_x,\Ero_x,\Ebd_{\le k_\infty-1,x},\Ecp_x\Bigg|\tilde\F\Bigg)\I_{\Ero_{A_{k_0}}}\\
&\ge c\exp(-2n+4\Delta)\I_{\Ero_{A_{k_0}}}.
\end{align*}
This is what we wanted to show.
\end{proof}

\subsection{Upper bound on the second moment}
Throughout this section, $\ell$ is fixed to be the same value as in Lemma \ref{l:lower_bound_first} (so that in particular Lemma \ref{l:lower_barrier} applies for $\ell$).

In this section we will show that
\[\E(\mathcal{N}_{\GGamma,N}^2\mid\tilde\F)\I_{\Ero_{A_{k_0}}}=\left(\E(\mathcal{N}_{\GGamma,N}\mid\tilde\F)\right)^2\I_{\Ero_{A_{k_0}}}(1+o_\Delta(1))\]
for any $\Delta\ge0$ and any $\GGamma\ge C(\Delta+1)$.
Once we know this, the Paley-Zygmund inequality will imply that $\mathcal{N}_{C(\Delta+1),N}\ge1$ with probability tending to 1 as $\Delta\to\infty$.

Let us abbreviate
\[\Ev_x=\Elw_x\cap\Ero_x\cap\Ebd_x\cap\Ecp_x\cap\left\{\phi\colon\phi(x)-\Sav_{k_0,-}(x,\phi)\ge \sqrt{\g}(m_N-C(\Delta+1))\right\},\]
so that
\begin{equation}\label{e:twopointfct}\E(\mathcal{N}_{\GGamma,N}^2\mid\tilde\F)=\sum_{x,x'\in A_{k_0}}\tilde\PP^{Q_N,0}(\Ev_x\cap\Ev_{x'}\mid\tilde\F),
\end{equation}
and so we will have to estimate the two-site probabilities on the right-hand side. To begin, note that the fields in different boxes $Q_{r_{k_0}}(x)$ for $x\in I_{k_0}$ are conditionally independent given $\F$, and so
\begin{equation}\label{e:twopointdiffboxes}
\tilde\PP^{Q_N,0}(\Ev_x\cap\Ev_{x'}\mid\tilde\F)=\tilde\PP^{Q_N,0}(\Ev_x\mid\tilde\F)\tilde\PP^{Q_N,0}(\Ev_{x'}\mid\tilde\F)
\end{equation}
if $x,x'$ are in different subboxes of $A_{k_0}$. This means that we can restrict attention in \eqref{e:twopointfct} to points $x,x'$ in the same subbox $Q_{r_{k_0}}(y)$ for some $y\in I_{k_0}$.

In any case, we have the trivial bound
\begin{equation}\label{e:twopointtrivial}\tilde\PP^{Q_N,0}(\Ev_x\cap\Ev_{x'}\mid\tilde\F)\le\tilde\PP^{Q_N,0}(\Ev_x\mid\tilde\F)\wedge\tilde\PP^{Q_N,0}(\Ev_{x'}\mid\tilde\F)
\end{equation}
that we will employ when $x,x'$ are very close. To make it useful, we need an upper bound for $\tilde\PP^{Q_N,0}(\Ev_x\mid\tilde\F)$ restricted on the nice events.  This is the subject of the following lemma.

\begin{lemma}\label{l:late_branching}
For any $\Delta>0,\GGamma\ge0$ and for any $N$ sufficiently large (depending on $\Delta$) we have for $x\in A_{k_0}$ the estimate
\begin{equation}\label{e:late_branching}
\tilde\PP^{Q_N,0}\left(\Elw_x,\Ero_x,\Ebd_x,\Ecp_x,\phi(x)-\Sav_{k_0,-}(x,\phi)\ge \sqrt{\g}(m_N-\GGamma)\middle|\tilde\F\right)\I_{\Ero_{A_{k_0}}}\le C\frac{\e^{4\Delta}}{N^2}.
\end{equation}
In particular,
\begin{equation}\label{e:upper_bound_first}
\E(\mathcal{N}_{\GGamma,N}\mid\tilde\F)\I_{\Ero_{A_{k_0}}}\le C\e^{4\Delta}.
\end{equation}
Moreover, for any $x,x'\in A_{k_0}$ we have
\begin{equation}\label{e:late_branching_twopoint}
\tilde\PP^{Q_N,0}(\Ev_x\cap\Ev_{x'}\mid\tilde\F)\I_{\Ero_{A_{k_0}}}\le C\frac{\e^{4\Delta}}{N^2}.
\end{equation}
\end{lemma}

The main challenge is to obtain good estimates in the intermediate regime. There we prove the following estimate.
\begin{lemma}\label{l:intermediate_branching}
Let $x,x'\in Q_{r_{k_0}}(y)$ for some $y\in I_{k_0}$, and suppose that $2r_{k+1}\le |x-x'|\le 2r_{k}$ for some $k$ with $k_0'\le k< k_\infty'$. Let also $\GGamma\ge0$ be arbitrary. Then we have
\begin{equation}\label{e:intermediate_branching}
\tilde\PP^{Q_N,0}(\Ev_x\cap\Ev_{x'}\mid\tilde\F)\I_{\Ero_{A_{k_0}}}\le C\frac{\e^{6\Delta}}{N^2r_{k}^2}\cdot\frac{((k-k_0)\wedge(k_\infty-k)))^{3/10}}{\e^{2((k-k_0)\wedge(k_\infty-k))^{2/5}/\sqrt{\g}}}.
\end{equation}
\end{lemma}

In the early branching regime we have a similar, but easier estimate.
\begin{lemma}\label{l:early_branching}
Let $x,x'\in Q_{r_{k_0}}(y)$ for some $y\in I_{k_0}$, and suppose that $2r_{k+1}\le |x-x'|\le 2r_{k}$ for some $k$ with $k_0\le k< k_0'$. Let also $\GGamma\ge0$ be arbitrary. Then we have
\begin{equation}\label{e:early_branching}
\tilde\PP^{Q_N,0}(\Ev_x\cap\Ev_{x'}\mid\tilde\F)\I_{\Ero_{A_{k_0}}}\le C\frac{\e^{8\Delta}}{N^4}.
\end{equation}
\end{lemma}

Before we prove these lemmas, let us show how to combine them into the desired upper bound on the second moment of $\mathcal{N}_{\GGamma,N}$.

\begin{lemma}\label{l:upper_bound_second}
For any $\Delta$ and for any $N$ sufficiently large,
\begin{equation}\label{e:upper_bound_second}
\left(\tilde\E^{Q_N,0}(\mathcal{N}_{\GGamma,N}^2\mid\tilde\F)-\left(\tilde\E^{Q_N,0}(\mathcal{N}_{\GGamma,N}\mid\tilde\F)\right)^2\right)\I_{\Ero_{A_{k_0}}}\le C\e^{4\Delta}.
\end{equation}
\end{lemma}
\begin{proof}
Using \eqref{e:twopointfct} and \eqref{e:twopointdiffboxes} we compute
\begin{align*}
&\tilde\E^{Q_N,0}(\mathcal{N}_{\GGamma,N}^2\mid\tilde\F)\I_{\Ero_{A_{k_0}}}\\
&=\sum_{x,x'\in A_{k_0}}\tilde\PP^{Q_N,0}(\Ev_x\cap\Ev_{x'}\mid\tilde\F)\I_{\Ero_{A_{k_0}}}\\
&=\sum_{y\in I_{k_0}}\sum_{x,x'\in Q_{r_{k_0}}(y)}\tilde\PP^{Q_N,0}(\Ev_x\cap\Ev_{x'}\mid\tilde\F)\I_{\Ero_{A_{k_0}}}+\sum_{\substack{y,y'\in I_{k_0}\\y\neq y'}}\sum_{\substack{x\in Q_{r_{k_0}}(y)\\x'\in Q_{r_{k_0}}(y')}}\tilde\PP^{Q_N,0}(\Ev_x\mid\tilde\F)\tilde\PP^{Q_N,0}(\Ev_{x'}\mid\tilde\F)\I_{\Ero_{A_{k_0}}}.
\end{align*}
We can bound the second summand by
\[\sum_{x,x'\in A_{k_0}}\tilde\PP^{Q_N,0}(\Ev_x\mid\tilde\F)\tilde\PP^{Q_N,0}(\Ev_{x'}\mid\tilde\F)\I_{\Ero_{A_{k_0}}}=\left(\tilde\E^{Q_N,0}(\mathcal{N}_{\GGamma,N}\mid\tilde\F)\right)^2,\]
and so \eqref{e:upper_bound_second} follows once we prove that
\begin{equation}\label{e:upper_bound_second1}
\sum_{y\in I_{k_0}}\sum_{x,x'\in Q_{r_{k_0}}(y)}\tilde\PP^{Q_N,0}(\Ev_x\cap\Ev_{x'}\mid\tilde\F)\I_{\Ero_{A_{k_0}}}\le C\e^{4\Delta}.
\end{equation}
To show this, we fix some $y\in I_{k_0}$ and split the inner sum according to the largest $k$ such that $|x-x'|\le 2r_{k}$ and use Lemmas \ref{l:early_branching}, \ref{l:intermediate_branching} and \ref{l:late_branching}, respectively. We find that
\begin{align*}
&\sum_{x,x'\in Q_{r_{k_0}}(y)}\tilde\PP^{Q_N,0}(\Ev_x\cap\Ev_{x'}\mid\tilde\F)\I_{\Ero_{A_{k_0}}}\\
&=\sum_{\substack{x,x'\in Q_{r_{k_0}}(y)\\|x-x'|\ge 2r_{k_0'}}}\tilde\PP^{Q_N,0}(\Ev_x\cap\Ev_{x'}\mid\tilde\F)\I_{\Ero_{A_{k_0}}}+\sum_{k=k_0'}^{k_\infty'-1}\sum_{\substack{x,x'\in Q_{r_{k_0}}(y)\\2r_{k+1}\le |x-x'|\le 2r_{k}}}\tilde\PP^{Q_N,0}(\Ev_x\cap\Ev_{x'}\mid\tilde\F)\I_{\Ero_{A_{k_0}}}\\
&\qquad\qquad+\sum_{\substack{x,x'\in Q_{r_{k_0}}(y)\\|x-x'|\le 2r_{k_\infty'}}}\tilde\PP^{Q_N,0}(\Ev_x\cap\Ev_{x'}\mid\tilde\F)\I_{\Ero_{A_{k_0}}}\\
&\le Cr_{k_0}^4\frac{\e^{8\Delta}}{N^4}+C\sum_{k=k_0'}^{k_\infty'-1}r_{k_0}^2r_{k}^2\frac{\e^{6\Delta}}{N^2r_{k}^2}\cdot\frac{((k-k_0)\wedge(k_\infty-k)))^{3/10}}{\e^{2((k-k_0)\wedge(k_\infty-k))^{2/5}/\sqrt{\g}}}+Cr_{k_0}^2r_{k_\infty'}^2\frac{\e^{4\Delta}}{N^2}\\
&\le C\e^{4\Delta}+C\e^{4\Delta}\sum_{k=k_0'}^{k_\infty'-1}\frac{((k-k_0)\wedge(k_\infty-k)))^{3/10}}{\e^{2((k-k_0)\wedge(k_\infty-k))^{2/5}/\sqrt{\g}}}+C\e^{4\Delta}.
\end{align*}
The sum over $k$ here is easily seen to be bounded uniformly, and so we obtain \eqref{e:upper_bound_second1}, which finishes the proof.
\end{proof}

Let us now give the proofs of the two crucial estimates, Lemma \ref{l:late_branching} and Lemma \ref{l:intermediate_branching}. The proofs are quite similar to the proof of Lemma \ref{l:upper_barrier_late}, but actually simpler, as we can directly work with the independent increments $\Ind_k$ instead of the $\Inc_k$. In fact, because of the latter simplification, the argument closely resembles the lower bound in Lemma \ref{l:lower_bound_first}.

\begin{proof}[Proof of Lemma \ref{l:late_branching}]
Note first that \eqref{e:upper_bound_first} and \eqref{e:late_branching_twopoint} follow immediately from  \eqref{e:late_branching} (the latter in combination with \eqref{e:twopointtrivial}). Therefore it suffices to prove \eqref{e:late_branching}.

As in the proof of Lemma \ref{l:lower_bound_first}, on the event $\Ecp_x$ we have $\Inc_k(x,\phi)=\Ind_k(x)$ for all $k_0\le k\le k_\infty$, and so $\Elw_x=\widetilde{\Elw_x}$ with $\widetilde{\Elw_x}$ as defined in \eqref{e:def_tildeElw}. So \eqref{e:late_branching} follows from
\begin{equation}\label{e:late_branching1}
\tilde\PP^{Q_N,0}\left(\widetilde{\Elw_x}\mid\tilde\F\right)\le C\frac{\e^{4\Delta}}{N^2}.
\end{equation}
In \eqref{e:lower_bound_first1} we showed the corresponding lower bound. As it turns out, the upper bound follows from basically the same argument. Namely with $\Q$ as in \eqref{e:def_Q} and $\lambda$ as in \eqref{e:def_lambda2} Lemma \ref{l:lower_barrier} implies that
\[\Q\left(\widetilde{\Elw_x}\right)\le \frac{C}{n^{3/2}}\]
On the event $\widetilde{\Elw_x}$ we have that $\sum_{j=k_0}^{k_\infty-1}\Ind_j(x)\ge\sqrt{\g}(m_N-4\Delta)-1$, and hence
\begin{align*}
&\tilde\PP^{Q_N,0}\left(\widetilde{\Elw_x}\middle|\tilde\F\right)\\
&\le\tilde\E^{Q_N,0}\left(\I_{\widetilde{\Elw_x}}\exp\left(\lambda \left(\sum_{j=k_0}^{k_\infty-1} \Ind_j (x,\phi)-\sqrt{\g}(m_N-4\Delta)+1\right)\right)\middle|\tilde\F\right)\\
&\le\exp\left(-\lambda (\sqrt{\g}(m_N-4\Delta)-1)\right)\tilde\E^{Q_N,0}\left(\I_{\widetilde{\Elw_x}}\exp\left(\lambda \sum_{j=k_0}^{k_\infty-1} \Ind_k (x,\phi)\right)\middle|\tilde\F\right)\\
&=\exp\left(-\lambda (\sqrt{\g}(m_N-4\Delta)-1)\right)\tilde\E^{Q_N,0}\left(\exp\left(\lambda \sum_{j=k_0}^{k_\infty-1} \Ind_j (x,\phi)\right)\middle|\tilde\F\right)\Q\left(\widetilde{\Elw_x}\right).
\end{align*}
After a short calculation just like the one that lead to \eqref{e:lower_bound_first1} we obtain \eqref{e:late_branching1}, as desired.
\end{proof}

\begin{proof}[Proof of Lemma \ref{l:intermediate_branching}]
As before, observe that on $\Ecp_x$ we have $\Elw_x=\widetilde{\Elw_x}$ with $\widetilde{\Elw_x}$ as in \eqref{e:def_tildeElw}, and the analogous statement holds with $x'$ in place of $x$.
This means that
\[\Ev_x\cap\Ev_{x'}\subset\widetilde{\Elw_x}\cap\widetilde{\Elw_{x'}}\subset \widetilde{\Elw_x}\cap\widetilde{\Elw_{\ge k,x'}},\]
and so it suffices to prove
\begin{equation}\label{e:intermediate_branching1}
\tilde\PP^{Q_N,0}\left(\widetilde{\Elw_x}\cap\widetilde{\Elw_{\ge k,x'}}\middle|\tilde\F\right)\le C\frac{\e^{6\Delta}}{N^2r_{k}^2}\cdot\frac{((k-k_0)\wedge(k_\infty-k)))^{3/10}}{\e^{2((k-k_0)\wedge(k_\infty-k))^{2/5}/\sqrt{\g}}}.
\end{equation}
The benefit of replacing $\widetilde{\Elw_{x'}}$ by $\widetilde{\Elw_{\ge k,x'}}$ is that $\widetilde{\Elw_x}$  and $\widetilde{\Elw_{\ge k,x'}}$ are independent. Indeed, by our assumption on $k$, $Q_{r_{k}}(x)$ and $Q_{r_{k}}(x')$ are disjoint. According to Lemma \ref{l:iterated_coupling} this means that $(\Ind_k(x))_{k=k_0}^{k_\infty-1}$ and $(\Ind_k(x',\phi))_{j=k}^{k_\infty-1}$ are all independent (conditionally on $\tilde \F$). For us this means that $\widetilde{\Elw_x}$  and $\widetilde{\Elw_{\ge k,x'}}$ are indeed conditionally  independent, so that
\begin{equation}\label{e:intermediate_branching2}
\tilde\PP^{Q_N,0}\left(\widetilde{\Elw_x}\cap\widetilde{\Elw_{\ge k,x'}}\middle|\tilde\F\right)=\tilde\PP^{Q_N,0}\left(\widetilde{\Elw_x}\middle|\tilde\F\right)\tilde\PP^{Q_N,0}\left(\widetilde{\Elw_{\ge k,x'}}\middle|\tilde\F\right).
\end{equation}
The first factor was already estimated in \eqref{e:late_branching1}. For the second factor we use they same tilting argument that we have seen a few times already. We pick $\Q$ as in \eqref{e:def_Q} and $\lambda$ as in \eqref{e:def_lambda2}. Then Lemma \ref{l:lower_barrier} implies that
\[\Q\left(\widetilde{\Elw_{\ge k,x'}}\middle|\tilde\F\right)\le C\frac{((k-k_0)\wedge(k_\infty-k)))^{3/10}}{(k_\infty-k)^{3/2}}.\]
Moreover, on the event $\widetilde{\Elw_{\ge k,x'}}$ we have
\[\sum_{j=k}^{k_\infty-1}\Ind_j(x)\ge\frac{\sqrt{\g}(m_N-4\Delta)}{k_\infty-k_0}(k_\infty-k)+((k-k_0)\wedge(k_\infty-k))^{2/5}-1,\]
and thus
\begin{align*}
&\tilde\PP^{Q_N,0}\left(\widetilde{\Elw_{\ge k,x'}}\middle|\tilde\F\right)\\
&\le\tilde\E^{Q_N,0}\left(\I_{\widetilde{\Elw_x}}\exp\left(\lambda \left(\sum_{j=k}^{k_\infty-1} \Ind_j (x,\phi)-\frac{\sqrt{\g}(m_N-4\Delta)}{k_\infty-k_0}(k_\infty-k)-((k-k_0)\wedge(k_\infty-k))^{2/5}+1\right)\right)\middle|\tilde\F\right)\\
&\le\exp\left(-\lambda \left(\frac{\sqrt{\g}(m_N-4\Delta)}{k_\infty-k_0}(k_\infty-k)+((k-k_0)\wedge(k_\infty-k))^{2/5}-1\right)\right)\\
&\qquad \qquad \times \tilde\E^{Q_N,0}\left(\I_{\widetilde{\Elw_x}}\exp\left(\lambda \sum_{j=k}^{k_\infty-1} \Ind_k (x,\phi)\right)\middle|\tilde\F\right)\\
&=\exp\left(-\lambda \left(\frac{\sqrt{\g}(m_N-4\Delta)}{k_\infty-k_0}(k_\infty-k)+((k-k_0)\wedge(k_\infty-k))^{2/5}-1\right)\right)\\
&\qquad\qquad\times\tilde\E^{Q_N,0}\left(\exp\left(\lambda \sum_{j=k}^{k_\infty-1} \Ind_k (x,\phi)\right)\middle|\tilde\F\right)\Q\left(\widetilde{\Elw_{\ge k,x'}}\middle|\tilde\F\right)\\
&\le C\exp\left(-\frac{(m_N-4\Delta)^2}{2(k_\infty-k_0)^2}(k_\infty-k)-\frac{m_N-4\Delta}{\sqrt{\g} (n-2\Delta)}((k-k_0)\wedge(k_\infty-k))^{2/5}\right)\\
& \qquad\qquad\qquad\qquad \qquad \qquad\qquad \qquad \qquad \qquad \qquad \times \frac{((k-k_0)\wedge(k_\infty-k)))^{3/10}}{(k_\infty-k)^{3/2}}\\
&\le C\exp\left(-2(k_\infty-k)+\frac{3\log n}{2n}(k_\infty-k)-\frac{2}{\sqrt{\g}}((k-k_0)\wedge(k_\infty-k))^{2/5}\right.\\
& \qquad \qquad \qquad\qquad \qquad \left. -\frac{3\log(k_\infty-k)}{2}+O_\Delta\left(\frac{\log n}{n^{3/5}}\right)\right)
\times((k-k_0)\wedge(k_\infty-k)))^{3/10}.
\end{align*}
As $(\log t)/{t}$ is a decreasing function of $t$ for $t\ge\e$, we have that $\frac{\log n}{n}(k_\infty-k)\le\log(k_\infty-k)$ (at least when $\ell\ge3$, say). So we can simplify further and obtain
\begin{equation}\label{e:intermediate_branching3}
\begin{split}
&\tilde\PP^{Q_N,0}\left(\widetilde{\Elw_{\ge k,x'}}\middle|\tilde\F\right)\\
&\le C\exp\left(-2(k_\infty-k)-\frac{2}{\sqrt{\g}}((k-k_0)\wedge(k_\infty-k))^{2/5}\right)\times((k-k_0)\wedge(k_\infty-k)))^{3/10}\\
&\le \frac{C}{r_{k}^2}\exp\left(2\Delta-\frac{2}{\sqrt{\g}}((k-k_0)\wedge(k_\infty-k))^{2/5}\right)\times((k-k_0)\wedge(k_\infty-k)))^{3/10}
\end{split}
\end{equation}
for $N$ large enough. Combining \eqref{e:late_branching1} and \eqref{e:intermediate_branching3} with \eqref{e:intermediate_branching2}, we easily obtain \eqref{e:intermediate_branching1}.

\end{proof}
\begin{proof}[Proof of Lemma \ref{l:early_branching}]
The argument is very similar to the one for Lemma \ref{l:intermediate_branching}. The main difference is that this time we modify both $\widetilde{\Elw_x}$ and $\widetilde{\Elw_{x'}}$ to gain independence. That is, we write
\[\Ev_x\cap\Ev_{x'}\subset\widetilde{\Elw_x}\cap\widetilde{\Elw_{x'}}\subset \widetilde{\Elw_{\ge k_0',x}}\cap\widetilde{\Elw_{\ge k_0',x'}}.\]
As in Lemma \ref{l:intermediate_branching}, the two events on the right-hand side are independent conditional on $\F$, and hence
\[\tilde\PP^{Q_N,0}\left(\widetilde{\Elw_x}\cap\widetilde{\Elw_{\ge k,x'}}\middle|\tilde\F\right)=\tilde\PP^{Q_N,0}\left(\widetilde{\Elw_{\ge k_0',x}}\mid\tilde\F\right)\tilde\PP^{Q_N,0}\left(\widetilde{\Elw_{\ge k_0',x'}}\middle|\tilde\F\right).\]
Both factors can be estimated in the same manner as in Lemma \ref{l:intermediate_branching}, and we obtain
\[\tilde\PP^{Q_N,0}\left(\widetilde{\Elw_{\ge k_0',x}}\middle|\tilde\F\right)\le C\frac{\e^{4\Delta}}{N^2},\]
as well as the analogous estimate for $x'$. The lemma follows immediately.
\end{proof}

\subsection{Completion of the proof}

As promised, we can now combine Lemma \ref{l:lower_bound_first} and Lemma \ref{l:upper_bound_second} with the Paley-Zygmund inequality to obtain a lower bound on the probability that $\mathcal{N}_{\GGamma,N}$ is positive.
\begin{lemma}\label{l:fine_field}
There is a constant $C$ such that for any $\Delta$ and for any sufficiently large $N$ (depending on $\Delta$) we have
\begin{equation}\label{e:fine_field}
\tilde\PP^{Q_N,0}\left(\exists x\in A_{k_0}\colon \phi(x)-\Sav_{k_0,-}(x,\phi)\ge \sqrt{\g}(m_N-C(\Delta+1)) \right)\ge 1-C\e^{-4\Delta}.
\end{equation}
\end{lemma}
\begin{proof}
In Lemma \ref{l:upper_bound_second} we have shown that for any $\GGamma\ge0$ we have
\[\left(\tilde\E^{Q_N,0}(\mathcal{N}_{\GGamma,N}^2\mid\tilde\F)-\left(\tilde\E^{Q_N,0}(\mathcal{N}_{\GGamma,N}\mid\tilde\F)\right)^2\right)\I_{\Ero_{A_{k_0}}}\le C\e^{4\Delta}\]
with the constant indepent of $\GGamma$. Meanwhile in Lemma \ref{l:lower_bound_first} we proved that for $\GGamma=C'(\Delta+1)$ with a sufficiently large $C'$ we have
\[\tilde\E^{Q_N,0}(\mathcal{N}_{\GGamma,N}\mid\tilde\F)\I_{\Ero_{A_{k_0}}}\ge c\e^{4\Delta}\I_{\Ero_{A_{k_0}}}.\]
The combination of these two estimates leads to
\[
\frac{\tilde\E^{Q_N,0}(\mathcal{N}_{C'(\Delta+1),N}^2\mid\tilde\F)}{\left(\tilde\E^{Q_N,0}(\mathcal{N}_{C'(\Delta+1),N}\mid\tilde\F)\right)^2}\I_{\Ero_{A_{k_0}}}\le 1+C\e^{-4\Delta},
\]
and hence by the Paley-Zygmund inequality
\begin{equation}\label{e:fine_field1}
\tilde\PP^{Q_N,0}(\mathcal{N}_{C'\Delta,N}>0\mid\tilde\F)\I_{\Ero_{A_{k_0}}}\ge\frac{\left(\tilde\E^{Q_N,0}(\mathcal{N}_{C'(\Delta+1),N}\mid\tilde\F)\right)^2}{\tilde\E^{Q_N,0}(\mathcal{N}_{C'(\Delta+1),N}^2\mid\tilde\F)}\I_{\Ero_{A_{k_0}}}\ge (1-C\e^{-4\Delta})\I_{\Ero_{A_{k_0}}}.
\end{equation}
In order to establish \eqref{e:fine_field}, we need to integrate this. We have
\begin{align*}
&\tilde\PP^{Q_N,0}\left(\exists x\in A_{k_0}\colon \phi(x)-\Sav_{k_0,-}(x,\phi)\ge \sqrt{\g}(m_N-C(\Delta+1))\right)\ge\tilde\PP^{Q_N,0}\left(\mathcal{N}_{C'(\Delta+1),N}>0\right)\\
&\ge\tilde\E^{Q_N,0}\left(\tilde\PP^{Q_N,0}(\mathcal{N}_{C'(\Delta+1),N}>0\middle|\tilde\F)\I_{\Ero_{A_{k_0}}}\right)
\ge(1-C\e^{-4\Delta})\tilde\PP^{Q_N,0}\left(\Ero_{A_{k_0}}\right),
\end{align*}
and \eqref{e:fine_field} now follows from \eqref{e:roughnessAk0}.
\end{proof}

Next, we show that the maximum of $S_{k_0,+}(x,\phi)$ over $x\in A_{k_0}$ is negligible. This is a chaining estimate similar to Lemma \ref{l:fluctuations_uncond}.

\begin{lemma}\label{l:coarse_field}
There is $C>0$ such that for any $\Delta>0$ and any sufficiently large $N$ we have
\begin{equation}\label{e:coarse_field}
\tilde\PP^{Q_N,0}\left(\exists x\in A_{k_0}\colon \left|\Sav_{k_0,+}(x,\phi)\right|\ge C\Delta\right)\le\frac{1}{n}.
\end{equation}
\end{lemma}
\begin{proof}
Heuristically, the field $\Sav_{k_0,+}(\cdot,\phi)$ has variance of order $k_0$, and typically has oscillations of order 1 in each box $Q_{r_{k_0}}(y)$ for $y\in I_{k_0}$. Thus, we expect the maximum of $\Sav_{k_0,+}(x,\phi)$ to be of order $C\log\left(\frac{N}{r_{k_0}}\right)\cdot k_0=Ck_0$, which is in line with \eqref{e:coarse_field}.

This can be made rigorous with a chaining argument similar to the proof of Lemma \ref{l:fluctuations_uncond},
using instead the exponential moment bound
\[
\tilde\E^{Q_N,0}\left(\exp\left(\lambda\Sav_{k_0,+}(x,\phi)\right)\right)\le \exp\left(C\lambda^2k_0\right).
\]
\end{proof}

We can now complete the proof of the lower bound on the maximum.
\begin{theorem}\label{t:lower_bound_max}
The random variable $\max_{x\in Q_N}\phi(x)-\sqrt{\g}m_N$ has a tight left tail. That is, for each $\eps>0$ there is $C_\eps>0$ such that for all $N\in\N$,
\begin{equation}\label{e:lower_bound_max}
\PP^{Q_N,0}\left(\max_{x\in Q_N}\phi(x)-\sqrt{\g}m_N\ge -C_\eps\right)=\tilde\PP^{Q_N,0}\left(\max_{x\in Q_N}\phi(x)-\sqrt{\g}m_N\ge -C_\eps\right)\ge 1-\eps.
\end{equation}
\end{theorem}
\begin{proof}
On the event $\Ebd_{k_0,x}$ we have $\left|\Sav_{k_0,-}(x,\phi)-\Sav_{k_0,+}(x,\phi)\right|\le 1$. Moreover, \eqref{e:badevents2} and a union bound imply that
\begin{equation}\label{e:lower_bound_max1}
\tilde\PP^{Q_N,0}\left(\bigcup_{x\in A_{k_0}}\Ebd_{k_0,x}\right)\le CN^2\exp(-cr_{k_0}^{\omega/2})\le \frac{C}{N}.
\end{equation}
We can now write
\begin{equation}\label{e:lower_bound_max2}
\begin{split}
&\tilde\PP^{Q_N,0}\left(\exists x\in A_{k_0}\colon\phi(x)\ge m_N-(2C+1)\Delta-1\right)\\
&\ge \tilde\PP^{Q_N,0}\left(\exists x\in A_{k_0}\colon\phi(x)-\Sav_{k_0,-}(x,\phi)\ge \sqrt{\g}(m_N-C(\Delta+1)\right)\\
&\qquad \qquad -\tilde\PP^{Q_N,0}\left(\exists x\in A_{k_0}\colon \left|\Sav_{k_0,+}(x,\phi)\right|\ge C\Delta\right)
-\tilde\PP^{Q_N,0}\left(\bigcup_{x\in A_{k_0}}\Ebd_{k_0,x}\right)\\
&\ge 1-C\e^{-4\Delta}-\frac{C}{\log N},
\end{split}
\end{equation}
where we used Lemma \ref{l:fine_field}, Lemma \ref{l:coarse_field} and \eqref{e:lower_bound_max1}.

Now, given $\eps>0$, we can choose $\Delta>0$ such that $1-C\e^{-4\Delta}\ge 1-\frac{\eps}{2}$. Then the right-hand side of \eqref{e:lower_bound_max2} is at least $1-\eps$ for all $N$ sufficiently large (depending on $\Delta,\eps$). So if we choose $C_\eps=2C+1$, then \eqref{e:lower_bound_max} holds for all sufficiently large $N$. But by increasing $C_\eps$, we can ensure that it actually holds for all $N$.
\end{proof}

\section{Proofs of the ballot estimates}\label{s:ballot_thms}
The goal of this section is to prove the ballot estimates in Lemma \ref{l:upper_barrier} and Lemma \ref{l:lower_barrier}. We will in fact prove a common generalization of these two lemmas.

 \begin{lemma}\label{l:upperlower_barrier}
 Let $m \in\N$ and $\gamma>0$. Then there is a constant $\ell>0$ with the following property. Let $(X_j)_{j=0}^{m}$ be independent random variables possessing a uniformly bounded density,  $C<\infty$, such that $\left|\E(X_j)\right| \le C\exp(\gamma(j-m))$ and $\left|\E(X_j^2)-\g\right|\le C\exp(\gamma(j-m))$ for any $j\ge 0$. We further assume that there exists $\lambda_*> 0$ and $C^*<\infty$  such that $\E\exp(\lambda_* X_j^2)< C^*$ for all $|\lambda|<\lambda_*$ and all $j\ge0$.  Let $\Sigma_k = \sum_{j=0}^{k-1} X_k$.

 Then there are $c_{\gamma,\ell},C_{\gamma,\ell}\ge 0$ depending on $\gamma,\ell$ as well as $\lambda_*,C^*$ such that
\begin{equation}
  \label{eq:100923a}
  \frac{ c_{\gamma,\ell}}{m^{3/2}}\le \PP \left(\Sigma_j \in\left[-(j\wedge (m-j))^{3/5}, -(j\wedge (m-j))^{2/5}\right] \ \forall j\in\{\ell,\ldots m-\ell\}, \Sigma_{m-1} \in [ -1, 0]\right)
  \le \frac{C_{\gamma,\ell}}{m^{3/2}}
\end{equation}
Moreover for any $k\in\{\ell,\ldots,m\}$, $h\geq 0$ and any $t>0$ we have 
  \begin{align}  
  &\begin{split}\label{eq:100923b}
    &\PP \left(-h+\Sigma_j\le  (j\wedge (m-j))^{2/5}\,\, \forall j\in\{\ell,\ldots,k-1\},\Sigma_{k}\in (k\wedge (m-k))^{2/5}+[-t-1,-t]\right)\\
    &\quad\le \frac{C_{\gamma,\ell}(1+h+\ell^{1/2})(1+t)}{(k-\ell)^{3/2}},
  \end{split}\\
 &\begin{split}\label{eq:100923b'}
    &\PP \left(-h+\Sigma_j\le  -(j\wedge (m-j))^{2/5}\,\, \forall j\in\{\ell,\ldots,k-1\},\Sigma_{k}\in -(k\wedge (m-k))^{2/5}+[-t-1,-t]\right)\\
    &\quad\le \frac{C_{\gamma,\ell}(1+h+\ell^{1/2})(1+t)}{(k-\ell)^{3/2}}.
  \end{split}  
  \end{align}
 \end{lemma}
 We note that 
 the upper bound in \eqref{eq:100923a}
 is a consequence of \eqref{eq:100923b'} with $k=m-\ell$, and Gaussian estimates for
 $\Sigma_{m-1}-\Sigma_{m-\ell}$. For this reason, we prove the
 upper bounds in \eqref{eq:100923b} and \eqref{eq:100923b'} and the lower bound
 in \eqref{eq:100923a}.

 The proof of Lemma
 \ref{l:upperlower_barrier} is based on Skorokhod embedding, which we recall in
 Section \ref{subsec-skorokhod}. We then prove the upper bounds in Section
 \ref{subsec-UBbar} and the lower bound in Section \ref{subsec-LBbar}.

 \subsection{The Skorokhod embedding}
 \label{subsec-skorokhod}
We begin with recalling a variant of the Skorokhod coupling that is constructed in such a way as to minimize the fluctuations of the Brownian motion $W_t$ for $t\le\tau$.
\begin{lemma}
  \label{lem-basicSkor}
  Let $X$ be a zero mean random variable, of unit variance and law
  $\mu$ which is absolutely continuous with respect to Lebesgue measure. Assume that $\mu$ possesses sub-Gaussian tails, that is, for some constant $C_\mu$,
  \[ \mu (|X|>t)\leq e^{-C_\mu t^2}.\]
  Then, there exist a Brownian motion $W_\cdot$,
  a stopping time $\tau=\tau_\mu$ with $\E\tau=1$ and constants $c_\mu$, $x_\mu$
  depending on $C_\mu$ only,
   so that $W_\tau=X$ in distribution and, for all $x>x_\mu$,
  \begin{equation}
    \label{eq-2907a}
    \PP\left(\sup_{t\leq \tau} |W_t|>x\right)\leq \e^{-c_\mu x^2}.
  \end{equation}
  Further,
  \begin{equation}
    \label{eq-2907b}
    \PP(\tau>t) \leq \e^{-c_\mu \sqrt{t}}.
  \end{equation}
\end{lemma}

\begin{proof}
  A good reference for such results is \cite{P86}. Let $\tau=T_{S}$ with
  $T_{S}$ the Skorokhod embedding
  defined in \cite[(2.2)]{P86}. Let $M=\max_{t\leq \tau} W_t$. Then, by
  item (a) in the proof of \cite[Proposition 2.1]{P86}, and with $\rho(x)$ defined in \cite[Page 176]{P86}, we have
  that for $\lambda>0$,
  \[ \PP(M\geq \lambda)\leq \int_\lambda^\infty \frac{x+\rho(x)}{\lambda+\rho(x)}
  d\mu(x)\leq \int_\lambda^\infty \Big(1+\frac{x}{\lambda}\Big) d\mu(x)
  \leq \e^{-c_\mu \lambda^2}.\]
  Applying this to $-X$ and adjusting $c_\mu$ if necessary gives
  \eqref{eq-2907a}.

  To see \eqref{eq-2907b}, write $M^*=\max_{s\leq \tau} |W_s|$. Then,
  \begin{equation}
  \begin{split}
    \label{eq-2907c}
    \PP(\tau>t)&=\PP(M^*\geq \lambda,\tau>t)+P(M^*<\lambda,\tau>t)\\
    &\leq
    \PP(M^*\geq \lambda)+P(\max_{s\leq t} |W_s|<\lambda) \\
    &\leq e^{-c_\mu \lambda^2} +e^{-c_w t/\lambda^2},
\end{split}  
  \end{equation}
  where $c_w>0$ is an absolute constant (that can be taken as half the 
 principal eigenvalue of the Dirichlet Laplacian on $[-1,1]$), and the last inequality follows from
  \eqref{eq-2907a}. Substituting $\lambda^2=\sqrt{t/c_\mu}$ and reducing $c_\mu$ if necessary  gives
  \eqref{eq-2907b}.
\end{proof}

\subsection{A basic one sided estimate}
\label{sec-onesided}
In this section we let $X_i$ denote a sequence of independent
random variables satisfying the assumptions of Lemma \ref{lem-basicSkor},
\emph{with the same constant $C_\mu$}.
We set $S_n=\sum_{i=1}^n X_i$, and use $\PP$ to denote probabilities for $X_i$ and $S_n$.

Let $\Phi(t)= t^{\theta}$ with $\theta<1/2$ and $\Psi(t)=t^{\iota}$ with $\iota>1/2$. Fix $\delta<1/2$. Introduce the events
\[ \mathcal{A}_{h,\ell, n,+}=\{-h+S_k\leq  \Phi(k), \; k=\ell, \ldots, n\},\qquad \mathcal{A}_{h,\ell, n,0}=\{-h+S_k\leq 0, \; k=\ell, \ldots, n\}
\]
and
\[\mathcal{A}_{h,\ell,n,-}=\{-\Psi(k)\leq -h+S_k\leq - \Phi(k), \; k=\ell, \ldots, n, S_n\in -n^{1/2}[(1-\delta),(1+\delta)]\}.\]
\begin{lemma}
\label{lem-onesided}
There exists a constant $c$ depending on $C_\mu$, $\theta$  only, so that for $h\geq 0$,
\begin{equation}
\label{eq-060823a}
\PP( \mathcal{A}_{h,\ell, n,+})\leq c\frac{(h+1+\ell^{1/2})}{\sqrt{n-\ell}}.
\end{equation}
Further, there exists a constant $c_-$ depending on $C_\mu$, $\theta,\iota,\delta$ only and a constant $\ell_0$ so that for all $\ell\geq \ell_0$, and all $h\in[0,1]$,
\begin{equation}
\label{eq-070823aa}
\PP( \mathcal{A}_{h,\ell, n,-})\geq \frac{c_-}{\sqrt{n}}.
\end{equation}
\end{lemma}
Note that the upper bound \eqref{eq-060823a} obviously
implies a similar bound with  $\mathcal{A}_{h,\ell, n,+}$ replaced by
$ \mathcal{A}_{h,\ell, n,0}$.
\begin{proof} 
The idea of the proof is to use the Skorokhod embedding from Lemma \ref{lem-basicSkor} to relate the barrier events for $S_k$ to barrier events for Brownian motion. This works well if the stopping times in the Skorokhod embedding are well behaved. In the unlikely event that at some time $\gg1$ a bad event occurs, we can consider the last such time. Then we gain a small factor from the probability of that event, and moreover, we can still use that the Brownian motion at even larger times will be below the barrier.

Throughout the proofs, $\theta,\iota,\delta$ are considered fixed and all constants may depend on them.
  Also, note that we may assume that $\ell<0.05 n$, because otherwise 
  \eqref{eq-060823a} is trivial whereas the left hand side of
  \eqref{eq-070823aa} is monotone increasing in $\ell$.
  For the upper bound, write $\mathcal{A}_{h,\ell, n}:= \mathcal{A}_{h,\ell, n,+}$. By Gaussian bounds, it is
  enough to consider $\ell>\ell_0$ for some fixed (i.e. independent from $n$) $\ell_0$.
  Let $\tau_i$ denote the stopping times given by Lemma \ref{lem-basicSkor}
  for $X_i$, and let $W_\cdot$ denote the concatenation of the Brownian motions
  given there, so that $\Sigma_k=W_{T_k}$, with $T_k=\sum_{i=1}^k \tau_i$. Note that
\[ \mathcal{A}_{h,\ell, n}=\{-h+W_{T_k}\leq \Phi(k), \; k=\ell, \ldots, n\}.\]

Call an index $k$ \emph{bad} if either $|T_k-k|>k/20$ or there is a $t\in (T_k,T_{k+1})$ so that $|W_t-W_{T_k}|>k^{\theta}/2$. Let $J_n-1$ denote the \emph{last} bad index
(with $J_n=\ell-1$ if no such bad time exists).
  Introduce the events
\begin{align}
\label{eq-190224a}
  \mathcal{B}_1&=\{ W_t\leq  2 \Phi(t)+h\ \forall 1.1 \ell  \leq t\leq 0.9 n\},\nonumber\\
\mathcal{B}_2&=\{ J_n\geq (\log n)^{1/\theta}\},\\
\mathcal{B}_3&= \{(\log n)^{1/\theta}\geq J_n> \ell-1,W_t\leq  4\Phi(t)+h\ \forall  0.9 n>t>T_{J_n}\}.\nonumber 
\end{align}
Note that for $\ell$ large enough we have
\begin{equation}
\label{eq-0708230}
 \mathcal{A}_{h,\ell, n}\subset \bigcup_{i=1}^3 \mathcal{B}_i.\end{equation}
 Indeed, on
 $ \mathcal{A}_{h,\ell, n}$, if there are no bad indices larger than or equal to
 $\ell$ (and hence $J_n=\ell-1$), then $\mathcal{B}_1$ must occur, while otherwise,
 note that
 for $k\geq J_n$, $T_k$ is comparable to $k$ and the path $W_t$
 does not stray, for $t\in [T_k,T_{k+1}]$, too much from $W_{T_k}$; thus, on
 $ \mathcal{A}_{h,\ell, n}$ either $\mathcal{B}_2$ occurs or
 the continuous path $W_t, t\in (T_{J_n},0.9 n)$ satisfies $\mathcal{B}_3$.

By standard Brownian motion estimates (compare e.g. with \cite[Lemma 2.7]{B83}, or see \cite[Proposition 2.2]{CHL19}) we have that
\begin{equation}
\label{eq-B2070823}
\PP(\mathcal{B}_1) \leq \frac{C(h+1+\ell^{1/2})}{\sqrt{n-\ell}}.
\end{equation}

  It is standard (and follows from Chebyshev's exponential inequality
  for the random variables $\tau_i\wedge c k$ for $c$ small, together with \eqref{eq-2907b}), that
  \begin{equation}
\label{eq-070823}
 \PP(|T_k-k|>k/20) \leq \e^{-c(\mu)\sqrt{k}},\end{equation}
and also, using \eqref{eq-2907a} and
modifying $c(\mu)$ if necessary,
\begin{equation}
\PP(
\exists t\in (T_k,T_{k+1}): \,
W_t-W_{T_k}>k^{\theta}/2)
\label{eq-B1060823}
\leq \e^{-c(\mu) k^{2\theta}}.
\end{equation}
Therefore, we obtain (again modifying $c(\mu)$) that
\begin{equation}
\label{eq-B3070823}
\PP(\mathcal{B}_2)
\leq  \e^{-c_\mu (\log n)^{1/2\theta}} + \e^{-c_\mu(\log n)^2}=o(1/\sqrt{n}).
\end{equation}
Concerning $\mathcal{B}_3$, we have that on that event, $|T_k-k|<k/20$ for all $k\geq J_n$ and hence for all such $k$,
$0.9\leq \Phi(T_k)/\Phi(k)\leq 1.1$. We then have that
\begin{align*}
\PP(\mathcal{B}_3) &\le \sum_{j=\ell}^{(\log n)^{1/\theta}}\PP\left(J_n=j,W_t\leq  4\Phi(t)+h\ \forall  0.9 n>t>T_{j}\right)\\
&\le \sum_{j=\ell}^{(\log n)^{1/\theta}}\PP\left(j-1\ \text{is bad},j\ \text{is good},W_t\leq  4\Phi(t)+h\ \forall  0.9 n>t>T_{j}\right)
\end{align*}
We can resample here $W_t$ for $t\ge T_j$ and replace it by an independent Brownian motion $\bar W_t$ (whose law when started from $y+1$ we denote by $\PP^{y+1}$). This allows us to estimate
\begin{equation}\label{eq-070823a}
\begin{split}
\PP(\mathcal{B}_3)&\leq \sum_{j=\ell}^{(\log n)^{1/\theta}} \sum_{y\in\Z}
\E \I_{\{|T_{j-1}-(j-1)|<(j-1)/20,|T_{j}-j|\leq j/20,\exists t\in (T_{j-1},T_j): |W_t-W_{T_j}|\geq (j-1)^{\theta}/2, W_{T_j}\in [y,y+1]\}} \\
& \qquad \qquad \times\PP^{y+1}(\bar  W_t\leq 4\Phi(t+T_{j})+h, t\in [0,n/2])\\
&+\sum_{j=\ell}^{(\log n)^{1/\theta}} \sum_{y\in\Z} \sum_{s=1}^\infty \E \I_{\{|T_{j-1}-(j-1)|\in [s(j-1)/20, (s+1) (j-1)/20] , |T_j-j|\leq j/20\}}
 \I_{ \{W_{T_{j}}\in [y,y+1]\} } \\
& \qquad \qquad  \times\PP^{y+1}( \bar W_t\leq 4\Phi(t+T_{j})+h, t\in [0,n/2]) \\
&=: P_1+P_2.
\end{split}
\end{equation}

By standard Brownian motion barrier estimates from \cite{B83} and on the event  $|T_j-j|<j/20$, we have that 
\[ \PP^{y+1}( \bar W_t\leq 4\Phi(t+T_{j})+h, t\in [0,n/2]) \leq C \frac{(y+h+\Phi(1.05j))\vee1}{\sqrt{n}}.\]
We then obtain using
\eqref{eq-2907a}, \eqref{eq-070823}, \eqref{eq-B1060823} and Gaussian estimates for $W_{1.05 j}$ together with the Cauchy-Schwarz inequality  that
\begin{equation}
\label{eq-B4070823}
P_1
\leq \frac{C}{\sqrt{n}} \sum_{j=\ell}^{(\log n)^{1/\theta}}\sum_{y\in\Z}  \e^{-c_\mu j^{2\theta}}  \e^{-cy^2/j} (h+y+(1.05 j)^\theta)\vee1\leq C\e^{-c_\mu\ell^{2\theta}} \frac{h+1+\ell^{1/2}}{\sqrt{n}}.
\end{equation}
Similarly,
\begin{equation}
\label{eq-B4070823a}
P_2
\leq \frac{C}{\sqrt{n}} \sum_{j=\ell}^{(\log n)^{1/\theta}}\sum_{y\in\Z} \sum_{s=1}^\infty \e^{-c_\mu\sqrt{sj}}   \e^{-cy^2/j} (h+y+(1.05 j)^\theta)\vee1\leq C\e^{-c_\mu\ell^{2\theta}} \frac{h+1+\ell^{1/2}}{\sqrt{n}}.
\end{equation}
 Altogether we obtain that
\begin{equation}
\label{eq-070823c}
\PP(\mathcal{B}_3)\leq C\e^{-c_\mu\ell^{2\theta}} \frac{h+1+\ell^{1/2}}{\sqrt{n}},
\end{equation}
which together with \eqref{eq-0708230}, \eqref{eq-B2070823} and \eqref{eq-B3070823},
completes the proof of \eqref{eq-060823a}.

We next turn to the proof of the lower bound \eqref{eq-070823aa}, which
uses the same ideas but needs some modifications that we spell  out.
Let $\tau_i, T_k$ and $W_\cdot$ be as in the proof of the upper bound,
and let $J_n^--1$ denote the last index $k$ so that $|T_k-k|>k/20$, with $J_n^-=\ell-1$ if no such index exists.
Similarly to \eqref{eq-190224a}, introduce the events
\begin{align*}
\mathcal{B}_1^-&=  \{-0.5\Psi(t)\leq  -h+W_t\leq  -2 \Phi(t), 0.9 \ell  \leq t\leq 0.9 n,  W_t\in -n^{1/2} [(1-\delta),(1+\delta)], t\in [0.9n,1.1n] \},\\
\mathcal{B}_2^-&=\{ J_n\geq (\log n)^{1/\theta}\},\\
\mathcal{B}_3^-&= \{(\log n)^{1/\theta}\geq J_n> \ell-1,
-0.5\Psi(t)\leq  -h+W_t\leq  -2 \Phi(t), J_n  \leq t\leq 0.9 n, \\
&\qquad \qquad \qquad \qquad \qquad \qquad \qquad \qquad
W_t\in -n^{1/2} [(1-\delta),(1+\delta)], t\in [0.9n,1.1n]\}.
\end{align*}
(Note that $\mathcal{B}_2=\mathcal{B}_2^-$, see \eqref{eq-190224a}.)

We now have, similarly to the argument for the upper bound, that
\begin{equation}
  \label{eq:130923}
  \mathcal{A}_{\theta,\ell, n,-}\supset \mathcal{B}_1^-\cap (\mathcal{B}_2^-)^\complement\cap (\mathcal{B}_3^-)^\complement.
\end{equation}
Therefore,
\begin{equation}
  \label{eq-130923aa}
  \PP( \mathcal{A}_{\theta,h,\ell, n,-})\geq P(\mathcal{B}_1^-)-
  \sum_{i=2}^3 P(\mathcal{B}_i^-).\end{equation}
By standard barrier estimates for Brownian motion, see e.g.
\cite{B83}, we have that
\begin{equation}
  \label{eq:130923b}
  P(\mathcal{B}_1^-)\geq \frac{2c_-}{\sqrt{h}}.
\end{equation}
Recall that in \eqref{eq-070823c},
$c(\ell)\to 0$ as $\ell\to\infty$, and therefore we can choose $\ell_0$ so that for $\ell>\ell_0$, $c(\ell)\leq c_-$.
With such choice, the lower bound \eqref{eq-070823aa} for all $\ell>\ell_0$ follows from \eqref{eq-130923aa}, \eqref{eq:130923b}, \eqref{eq-B3070823} 
and the analogue of
 \eqref{eq-070823c}.
\end{proof}
\subsection{The upper bound}
\label{subsec-UBbar}
We can now provide the proof of the upper bounds in Lemma \ref{l:upperlower_barrier}.
\begin{proof}[Proof of \eqref{eq:100923b}.]
  Set $\bar X_i=(X_i-EX_i)/\sqrt{\mbox{\rm Var} X_i}$ and $S_j=\sum_{i=1}^j \bar X_i$.
Our first goal is to transfer \eqref{eq-060823a} to $\Sigma_j$.
Note that
$\Delta_j:=|\Sigma_j-S_j|\leq C+C\sum_{i=1}^j |X_i| e^{-\gamma i}$  and that
$\PP(\Delta_j\geq j^{0.1})\leq Ce^{-c j^{0.1}}$. Now, using a decomposition according to the last time
where $\Delta_j>j^{0.1}$,
\begin{align*}
  &\PP(-h+\Sigma_j\leq j^{2/5}, j=\ell,\ldots, m/3)- P(-h+S_j\leq 2j^{2/5}, j=\ell,\ldots,m/3)\\
&\leq \sum_{l=\ell}^{m/3} \PP(\Delta_l\geq l^{0.1}, -h+S_j\leq 2j^{2/5}, j=l+1, \ldots, m/3)\\
&\leq c m^{-2}+\sum_{l=\ell}^{(\log m)^{20} }\PP(\Delta_l\geq l^{0.1}, -h+S_j\leq 2j^{2/5}, j=l+1, \ldots, m/3) \\
&\leq  cm^{-2}+\sum_{l=\ell}^{(\log m)^{20} }\PP(\Delta_l\geq l^{0.1}, S_{l+1}\geq -l^{1/2}, -h+S_j\leq 2j^{2/5}, j=l+1, \ldots, m/3) \\
&\quad +\sum_{l=\ell}^{(\log m)^{20} }\sum_{s=1}^\infty \PP(\Delta_l\geq l^{0.1}, S_{l+1}\in [-(s+1)l^{1/2},-s l^{1/2}],  -h+S_j\leq 2j^{2/5}, j=l+1, \ldots, m/3)
\\
&\leq  cm^{-2}+\sum_{l=\ell}^{(\log m)^{20} } \e^{-cl^{0.1}} \frac{h+l^{1/2}}{\sqrt{m-\ell}}+\sum_{l=\ell}^{(\log m)^{20} }\sum_{s=1}^\infty \e^{-cl^{0.1}} \e^{-cs^2}
\frac{h+1+sl^{1/2}}{\sqrt{m-\ell}}\leq c\frac{h+1+\ell^{1/2}}{\sqrt{m}}.
\end{align*}
This yields
 \eqref{eq-060823a} with $S_.$ replaced by $\Sigma_.$.

We now turn to the proof of the upper bound in
\eqref{eq:100923b}. We first assume that $k\ge m/2$.
Introduce
\[ p_{m,\ell,t,k,h}:=\PP(-h+\Sigma_j\leq (j\wedge (m-j))^{2/5}\,\,, j=\ell,\ldots, k, \Sigma_k\in [-(t+1),-t]+(m-k)^{2/5}).\]
Set $\bar \Sigma_j= -\sum_{i=k-j+1}^{k-1} X_i$. Introduce the events 
\begin{align*}
{\mathcal B}_1&=\{-h+\Sigma_j\leq j^{2/5}, j=\ell,\ldots, m/6\}\\
{\mathcal B}_2&=\{-(t+1)+\bar \Sigma_j\leq 
((m-k+j)\wedge (k-j))^{2/5}-(m-k)^{2/5}, 
j=1,\ldots,m/6\}\\
{\mathcal B}_3&=\left\{ \sum_{j=m/6+1}^{k-m/6} X_j \in \Sigma_{m/6}-\bar\Sigma_{m/6}+t+[-2,2]\right\}.
\end{align*}  
(Do not confuse these with the events in Lemma \ref{lem-onesided}.)

Note that, by monotonicity,
\[p_{m,\ell,t,k,h}\leq \PP\left(\bigcap_{i=1}^3 {\mathcal B}_i\right)
\leq \prod_{i=1}^2 \PP({\mathcal B}_i) \times \sup_x \PP\left(\sum_{j=m/6+1}^{k-m/6} X_j\in x+[-2,2]\right),\]
where in the second inequality we used the Markov property and independence of
increments. By \eqref{eq-060823a} applied to $\Sigma_.$,
we have that
\[\PP({\mathcal B}_1)\leq C(1+h+\ell^{1/2})/\sqrt{m-\ell}\leq C\frac{1+\ell^{1/2}}
{\sqrt{m}}.\]
Similarly, using that 
\[((m-k+j)\wedge (k-j))^{2/5}-(m-k)^{2/5}\le (m-k+j)^{2/5}-(m-k)^{2/5}\leq j^{2/5},\]
we obtain from \eqref{eq-060823a} applied to $\bar \Sigma_\cdot$ that
\[\PP({\mathcal B}_2)\leq C_{\ell}(t+1)/\sqrt{m}\leq C_{\ell}\frac{t+1}{\sqrt{m}}.\]
On the other hand, by the local CLT as in e.g. \cite[Corollary 1.3]{D16}
\[\sup_x \PP\left(\sum_{j=m/6+1}^{k-m/6} X_j\in x+[-2,2]\right)\leq \frac{C}{\sqrt{m-k/3}}\le \frac{C}{\sqrt{m}}\]
(note that this conclusion holds for either of the alternatives in \cite[Corollary 1.3]{D16}).
Combining these estimates yields
\eqref{eq:100923b}.

The same argument works for $3\ell\le k<m/2$. The only difference is that we 
divide the interval $[\ell, k]$ in three equal pieces, and use the estimate for 
$\mathcal{B}_1$ with $m/6$ replaced by $\ell+(k-\ell)/3$,
and, using that $j^{2/5}-k^{2/5}<0$, replace $\mathcal{B}_2$ by 
$\{-(t+1)+\bar \Sigma_j\leq 0, j=1,\ldots, (k-\ell)/3\}$.

Finally, the case $k<3\ell$ is trivial, provided we choose $C_{\gamma,\ell}$ large enough.
\end{proof}

\begin{proof}[Proof of \eqref{eq:100923b'}.]
This is very similar to the proof of \eqref{eq:100923b}. The fact that the barrier is curved downwards only makes the arguments easier.
\end{proof}

\begin{proof}[Proof of Lemma \ref{l:upper_barrier}]
 \eqref{eq-230124b} follows from  \eqref{eq:100923b} upon taking $\tilde\ell=k_0\vee\ell$ and $h=a$. \eqref{eq-230124a} follows from \eqref{eq:100923b} upon taking $\tilde\ell=k_0\vee\ell$, $k=m-\ell$, $t=s$  and $h=a$ in the latter, and using the local CLT
to estimate $\PP(\Sigma_m-\Sigma_{m-\ell}\in [s-t-2,s-t])$ together with summation on $s$ integer. 
\end{proof}

\subsection{The lower bound}
\label{subsec-LBbar}
 The proof of the lower bound in Lemma \ref{l:upperlower_barrier} follows a pattern similar to the proof of the upper bound, with one important difference in that the division of $m$ to three roughly equal intervals is replaced by a division where the middle interval is much smaller. 

In more details, recall that we may and will assume in the lower bound that $\ell$ is fixed independent of $m$. We  prove the lower bound in \eqref{eq:100923a}
for $S_\cdot$ instead of $\Sigma_\cdot$. As in the proof of the upper bound,
one may then transfer the
estimate to $\Sigma_\cdot$.

Recall the variables $\tau_i$, $T_i$. Introduce the event
\[\mathcal{D}=\{\exists k>0.1 m: |T_k-k|>k^{3/4}\}.\]
By \eqref{eq-2907b} and Chebycheff's inequality, we obtain that
$P(\mathcal{D})\le 0.9m \e^{-cm^{3/8}}=o(1/m^2)$. Set $m_-=m/2-\eps m$.
Introduce also the event
\[\mathcal{D}'=\{\exists k\in [m_-,m-m_-]: |S_k-S_{m_-}|>m^{1/2}/4\}.\]
By
\eqref{eq-2907a} and Chebycheff's inequality,
 we also have that $P(\mathcal{D}')\leq c_\epsilon m^{-1/2}$, with $c_\epsilon\to 0$ as $\epsilon \to 0$.

Define now the events
\begin{align*}
\mathcal{B}_{1}&=\{-k^{3/5}\leq S_k\leq - k^{2/5}, \; k=\ell, \ldots, m_-, S_{m_-}\in -m^{1/2}[(1-\delta),(1+\delta)]\},\\
\mathcal{B}_{2}&=\{-(m-k)^{3/5}\leq S_m-S_k\leq - (m-k)^{2/5}, \; k=m-\ell, \ldots, m-m_-, \\
& \qquad \qquad \qquad\qquad \qquad \qquad S_m-S_{m_-}\in -m^{1/2}[(1-\delta),(1+\delta)]\},\\
\mathcal{B}_3&=\{|S_{m-m_-}-S_{m_-}|\leq 1/2\}.
\end{align*}
With $\mathcal{C}'$ as in the statement of Lemma \ref{l:lower_barrier} but with $S_\cdot$ replacing $\Sigma_\cdot$, we have that
\[\mathcal{C}'\supset \cap_{i=1}^3 \mathcal{B}_i \cap \mathcal{D}\cap \mathcal{D}'.\]
Therefore, using the Markov property for $S_\cdot$,
\begin{align*}
  \PP(\mathcal{C'})&\geq
  \PP(\mathcal{B}_3\cap \mathcal{D}') \prod_{i=1}^2 \PP(\mathcal{B}_i)-o(1/m^2)\\
  &\stackrel{\eqref{eq-070823aa}}{\geq}
  \frac{c_-^2}{m} \PP(\mathcal{B}_3\cap \mathcal{D}')-o(1/m^2)\\
  &\geq
  \frac{c_-^2}{m} (\PP(\mathcal{B}_3)-\frac{c_\epsilon}{m^{1/2}}) -o(1/m^2),
\end{align*}
where $c_->0$ depends on $\delta$ and $\ell$ only.
Choosing $\delta=\delta(\epsilon)$ small enough, we obtain that
$\PP(\mathcal{B}_3)\geq C_\epsilon/\sqrt{m}$ by the local CLT in \cite{M78} or \cite[Corollary 3.3]{D16} (where in the latter, case (i) holds by virtue of the assumption of uniformly bounded density
or existence of a component with bounded density, as follows e.g. from \cite[Corollary 3.3]{D16}),
where now $C_\epsilon$ is bounded below uniformly in $\epsilon$ (and in fact,
with e.g. the choice $\delta=\epsilon^4$, one obtains $C_\epsilon\to\infty$ as $\epsilon\to 0$).
Choosing $\epsilon$ small enough so that $C_\epsilon>2c_\epsilon$ completes the proof.\qed

\begin{proof}[Proof of Lemma \ref{l:lower_barrier}]
  The lemma follows from \eqref{eq:100923a}.
\end{proof}

\bibliographystyle{alpha_edited2}
\bibliography{MaximumGLSecondOrder}

\end{document}